\documentclass[11pt,oneside,reqno]{article}
\pdfoutput=1

\usepackage{amsmath}
\usepackage{amsfonts}
\usepackage{amssymb}

\usepackage{graphicx}
\usepackage{subfig}
\usepackage{placeins}

\usepackage{color}
\usepackage{authblk}

\usepackage{hyperref}  
\hypersetup{
    linkcolor=blue,    
    citecolor=red,     
    colorlinks=true
}

\usepackage[x11names]{xcolor}               
\usepackage{todonotes}						

\newcommand\pare[1]{\left(#1\right)}                     
\newcommand\crochets[1]{\left[#1\right]}                 
\newcommand\abs[1]{\left\lvert#1\right\rvert}            

\newtheorem{remark}{Remark}
\newcommand{\ds}{\displaystyle}
\newcommand{\ud}{\,{\mathrm{d}}}

\font\bba=msbm10
\font\bbb=msbm8
\font\bbc=msbm6
\newfam\bbfam
\textfont\bbfam=\bba
\scriptfont\bbfam=\bbb
\scriptscriptfont\bbfam=\bbc
\def\bb{\fam\bbfam\bba}

\def\N{{\bb N}}
\def\Z{{\bb Z}}
\def\r{{\bb R}}

\def\ve{\varepsilon}

\makeatletter
\newcommand{\subjclass}[2][2010]{%
  \let\@oldtitle\@title%
  \gdef\@title{\@oldtitle\footnotetext{#1 \emph{Mathematics subject classification.} #2}}%
}
\makeatother

\title
{Asymptotic problems and numerical schemes for traffic flows with unilateral constraints describing the formation of jams}
\author[1]{F. Berthelin\thanks{{\tt Florent.Berthelin@unice.fr}}}
\author[1]{T. Goudon\thanks{ {\tt thierry.goudon@inria.fr}}}
\author[1,2]{B. Polizzi\thanks{ {\tt bastien.polizzi@imft.fr}}}
\author[3]{M. Ribot\thanks{{\tt magali.ribot@univ-orleans.fr}}}
\affil[1]{Universit\'e C\^ote d'Azur, Inria, CNRS, LJAD}
\affil[2]{Institut de M\'ecanique des Fluides de Toulouse, CNRS UMR~5502}
\affil[3]{Universit\'e d'Orl\'eans, MAPMO, UMR CNRS 7349}

\begin{document}
\maketitle

\begin{abstract}
We discuss numerical strategies to deal with PDE systems describing traffic flows, taking into account a density threshold, which restricts the vehicles density in the situation of congestion. 
These models are obtained through asymptotic arguments.
Hence, we are interested in the simulation of approached models that contain stiff terms and large speeds of propagation.
We design schemes intended to apply with relaxed stability conditions.
\end{abstract}


\section{Introduction}

In order to describe traffic flows and to reproduce the formation of congestions, several models based either on Ordinary Differential Equations (ODE) or Partial Differential Equations (PDE) have been proposed. Starting from individual--based "Follow-the-Leader" models \cite{gazis}, a very active stream in the traffic community considers now PDE models. A first example dates back to Lighthill and Whitham in the 50's \cite{LW}: the evolution of the density of cars is described by means of a mass conservation equation, where the flux is defined by a prescribed function of the density. In these so-called first-order models, the relation between flux and density is referred to as the fundamental diagram in the traffic flows community. A more accurate description can be expected by considering second-order models where a system of PDE governs the evolution of the density and the speed of cars. 
A first attempt in this direction is due to Payne \cite{payne}, strongly inspired by the principles of fluid mechanics. However, Daganzo \cite{daganzo} pointed out the drawbacks of this approach: the Payne-Whitham model may lead to inconsistent behaviors for the flow, such as vehicles going backwards. 
The model introduced independently by Aw and Rascle \cite{AwR1} and by
Zhang \cite{Zhang}, which still has the form of a $2\times 2$ system of conservation laws, is intended to correct these inconsistencies.
In \cite{AwR2}, a derivation of the system is proposed from a Follow-the-Leader model. We can also mention that some kinetic models \cite{BBNS, Nel, PF, PSTV, WK} are 
 under consideration, after the pioneering work \cite{Prig}. Further details and references can be found in the survey~\cite{BeDo}.
\medskip

This work is concerned with the numerical simulation of certain variants of the Aw-Rascle-Zhang  model.
Let $\rho(x,t)$ and $v(x,t)$ be the density and the velocity of cars at position $x \in \r$ and time $t>0$, respectively. The Aw-Rascle-Zhang model writes 
\begin{align}
\label{AW1}
	\left\{
	\begin{array}{l}
		\partial_t \rho+\partial_x(\rho v)=0,     \\
		\partial_t \big(v+p(\rho)\big)+v\partial_x \big(v+p(\rho)\big)=0,
	\end{array}
	\right.
\end{align}
where $\rho\mapsto p(\rho)$ plays the role of the pressure in the gas dynamic equations. 
In fact, the quantity $w=v+p(\rho)$ describes \ the desired velocity of the drivers, whereas $v$ corresponds to the actual velocity of the cars. Therefore, the term $p(\rho)$, the velocity offset, stands for the difference between these two velocities, reflecting the fact that the drivers slow down because of the density of cars.
It is convenient to rewrite the equations in the more convenient form of a \emph{conservative} system; namely \eqref{AW1} is, at least formally, equivalent to
\begin{align}
\label{AW1b}
	\left\{
	\begin{array}{l}
		\partial_t \rho+\partial_x(\rho v)=0,		\\
		\partial_t \big(\rho(v+p(\rho))\big)+\partial_x \big(\rho v(v+p(\rho))\big)=0.
	\end{array}
	\right.
\end{align}

Of course, a crucial modeling issue relies on the expression of the velocity offset $p(\rho)$.
At first glance, again inspired from gas dynamics, we can set $p(\rho)=\rho^\gamma$ for some $\gamma> 1$. However, such a model does not permit to impose a priori a limitation to the cars density. Consequently, Berthelin, Degond, Delitala and Rascle proposed in \cite{BeDe} to define the velocity offset as follows: given $0<\rho_\star<\infty$, 
\begin{align*}
\rho\in [0,\rho_\star)\longmapsto p(\rho)
=\left(\ds\frac{\rho_\star\rho}{\rho_\star-\rho}\right)^{\gamma}, \qquad \gamma> 1,
\end{align*}
which can be rewritten, for $\rho \neq 0$, by
$p(\rho)=\left(\ds\frac{1}{\rho}-\ds\frac{1}{\rho_\star}\right)^{-\gamma}$,
where $\rho_{\star}$ denotes a maximal value for the density. 
The velocity offset tends to infinity when $\rho \to \rho_{\star}$ while we get the classical expression $\ds p(\rho) \sim \rho^\gamma$ when $\rho \to 0$.
Such a pressure law also arises in gas dynamics, where it is referred to as the Bethe--Weyl law \cite{BGS}; for instance it is used to model close--packing effects in multi-fluid flows, see \cite{BGM} and the references therein.
This expression for $p$ has the role of enforcing that $\rho$ satisfies the constraint $\ds 0 \leq \rho(x,t) \leq \rho_\star$ for all $x \in \r$ and $t>0$, as it can be seen from the bounds on the Riemann invariants of the system \cite{BeDe}.
Moreover, \cite{BeDe} points out that drivers do not reduce significantly their speed unless they reach a congested region. Accordingly the velocity offset is appropriately rescaled with a small parameter $0<\ve\ll 1$, and we use 
\begin{align}
\label{pr1}\tag{{\bf{VO1}}}
	p^{\ve}(\rho)=\ve\left(\ds\frac{\rho_\star\rho}{\rho_\star-\rho}\right)^{\gamma}, \quad
	\textrm{ for } 0 \leq \rho < \rho_\star,
\end{align}
in the Aw-Rascle-Zhang model. We are thus led to the Rescaled Modified Aw-Rascle (RMAR) system
\begin{align}
\label{RMAR}
\begin{cases}
\partial_t \rho^{\ve}+\partial_x\pare{\rho^{\ve} v^{\ve}}=0,
\\
\partial_t \big(\rho^{\ve}(v^{\ve}+p^{\ve}(\rho^{\ve}))\big)+
\partial_x \big(\rho^{\ve} v^{\ve} (v^{\ve}+p^{\ve}(\rho^{\ve}))\big)=0.
\end{cases}
\end{align}
In this model, the velocity offset is small unless the density is getting close to the threshold $\rho_\star$.
Finally, \cite{BeDe} studies the limit when $\ve \to 0$ in \eqref{RMAR}.
In this regime we obtain (at least formally) the constrained system 
\begin{align}\label{PGCD}
\left\{
\begin{array}{l}
\partial_t \rho+\partial_x(\rho v)=0,
\\
\partial_t(\rho(v+\pi))+\partial_x(\rho v(v+\pi))=0,
\\
0 \leq \rho \leq \rho_{\star}, \, \pi\geq 0, \, (\rho_\star-\rho)\pi=0.
\end{array}
\right.
\end{align}
In \eqref{PGCD}, the limit ``pressure'' 
 $\pi=\lim_{\ve \to 0 }\ds p^\ve (\rho^{\ve})$ appears as the Lagrange multiplier associated to the unilateral constraint $0 \leq \rho \leq \rho_{\star}$. 
In particular, $\pi$ becomes active only in the congested regions, where the density reaches the threshold $\rho_\star$. 
Otherwise, in absence of congestion, the system is reduced to the pressure-less gas dynamic model \cite{BrGr,Gr}
\begin{align}\label{PGD}
\left\{
\begin{array}{l}
\partial_t \rho+\partial_x(\rho v)=0,	\\
\partial_t (\rho v)+\partial_x (\rho v^2)=0.
\end{array}
\right.
\end{align}
The asymptotic model is further investigated in \cite{BeDe}, 
exhibiting the formation of clusters, and proving the existence of weak solutions to the system \eqref{PGCD} through the stability analysis of ``sticky blocks''
dynamics.
It is also worth pointing out the original numerical approach developed in \cite{Maury} for \eqref{PGCD}
which uses ideas from the modelling of crowd motion and  includes a fine description of the non elastic collision processes.

\medskip

The asymptotic system \eqref{PGCD} is thus specifically intended to describe the formation and the dynamics of jams.
In this paper, we are interested in the numerical simulations of the system \eqref{PGCD}, and 
in the asymptotic regime $\ve\to 0$ in \eqref{RMAR}.
The difficulty is two--fold. On the one hand, in the free flow case, it is well-known that the pressureless gas dynamics system \eqref{PGD} can lead to delta-shocks formation, which makes it difficult to treat numerically \cite{BJL}. 
On the other hand, with the formation of a congestion, there is no direct access to the limit velocity offset $\pi$ which is defined in a quite abstract way. 
Therefore, in order to go beyond the simple particulate approach in \cite{BeDe}, we wish to develop numerical simulations of the RMAR model \eqref{RMAR} with the velocity offset \eqref{pr1} for small values of $\ve$. 
We are still facing several numerical challenges. 
Firstly, the model prohibits that the density exceeds the threshold $\rho_\star$.
Secondly, (one of) the characteristic speeds of the system become very large in congested region, which makes the time step shrink: the smaller $\ve$, the more severe the stability constraint.
Therefore, we need to design a scheme which can preserve the natural estimates of the problem, in particular 
the density limitation. 
As already observed in \cite{ChG, CG08} standard schemes may fail this objective due to the very specific structure of the PDE system.
We also refer the reader to \cite{Ju} for further examples related to fluid mixtures.
Moreover, we would like to relax the stability constraints on the time step.
To this end, a first attempt would be to adapt the explicit-implicit method proposed in \cite{DHN} for the 
Euler system with congestion constraint. However, this method applied to \eqref{RMAR} 
is not  satisfactory: it produces excessively 
smoothed density profiles and it overestimates the velocity 
 in congested regions. 
Thus we design a different splitting strategy, partly inspired from \cite{DeTa}.
Beyond the conception of a numerical scheme able to handle the stiffness of the problem,
 we will also discuss different asymptotic approaches of the constrained problem \eqref{PGCD}, based on different definitions of the scaled offset velocity in \eqref{RMAR}
 which all lead asymptotically to \eqref{PGCD}. 
 It is interesting to study how the shape of the pseudo--pressure affects the intermediate states (for not so extreme values of the scaling parameter), and the numerical costs.

The outline of this article is the following. In Section~\ref{Prop}, we go back to some properties of the Aw-Rascle-Zhang system and we detail the numerical difficulties we face.
Additionally, we propose different velocity offsets and scaling that can be used to recover asymptotically the constrained system \eqref{PGCD}. 
Then, in Section~\ref{scheme}, we propose a new explicit-implicit scheme based on a splitting strategy. 
The splitting is constructed to reduce the characteristic speeds in the explicit part
so that we can expect to use larger time steps.
Finally, in Section~\ref{simu}, we display some numerical simulations in order to prove the efficiency of the scheme and to compare the behavior of the system when using different velocity offsets.

\section{Properties of the Aw-Rascle-Zhang model and numerical difficulties}
\label{Prop}

We will describe in this Section the main numerical difficulties we have to deal with, when computing solutions of system \eqref{RMAR}.

\subsection{Different velocity offsets}

With the velocity offset \eqref{pr1}, it is forbidden to produce numerical densities larger than the threshold $\rho_\star$: if, due to any numerical error, the code returns a density larger than $\rho_\star$, we cannot update further the system and the simulation breakdowns.
 To cope with this difficulty, we propose to slightly modify the law, replacing \eqref{pr1} by a function $\rho\mapsto \tilde p^{\ve}(\rho) $ which is defined for any positive entry, which behaves like $p^{\ve}$ for $\rho<\rho_\star$, and which blows up as $\varepsilon\to 0$ for $\rho\geq \rho_\star$.
For instance, we set
\begin{align}
\label{pr1bis}\tag{{\bf{VO2}}}
\tilde p^{\ve}(\rho)=
\left\{
\begin{array}{ll}
\ve \left(\ds\frac{\rho_\star\rho}{\rho_\star-\rho}\right)^{\gamma}, & \text{ if } \rho \leq \rho_{tr}^{\ve},
\\
c^0_\ve + c^1_\ve(\rho-\rho_{tr}^{\ve})+ c^2_\ve \ds\frac{(\rho-\rho_{tr}^{\ve})^2}{2}, & \text{ if } \rho>\rho_{tr}^{\ve}.
\end{array}
\right.
\end{align}
In this formula, $\rho_{tr}^{\ve}$ is a transition density, which has a modeling nature; it should satisfy $\rho_{tr}^{\ve}
\to \rho_\star$ as $\ve$ tends to 0.
Beyond the transition, $\tilde p^{\ve}$ is a second order polynomial, computed so that $\tilde p^{\ve}$ remains a $\mathcal C^2$ function.
We thus set
$$\begin{array}{l}
\rho_{tr}^{\ve}=\rho_\star-h(\ve),
\\
c^0_{\ve}=p^{\ve}(\rho_{tr}^{\ve})=\ve\left(\ds\frac{\rho_\star(\rho_\star-h(\ve))}{h(\ve)}\right)^\gamma,
\qquad
c^1_{\ve}=(p^{\ve})'(\rho_{tr}^{\ve}),\qquad
c^2_{\ve}=(p^{\ve})''(\rho_{tr}^{\ve}).
\end{array}$$ 
The expected behavior holds for instance with $h(\ve)=\ve$, since it satisfies the two following properties when $\ve \to 0$: $h(\ve) \to 0$ and $\ds c^0_{\ve} \to +\infty$.
We point out that $\ds \rho_{tr}^{\ve}$ is purely a modeling parameter
and the model \eqref{pr1bis} leads us to the same difficulties as \eqref{pr1} in terms of stability issues, as we will see in the next Section.
\\

Another option is to use the following velocity offset
\begin{align}\label{pr2}\tag{{\bf{VO3}}}
p^{\gamma}(\rho)=V_{\mathrm{ref}}
\left(\ds\frac{\rho}{\rho_\star}\right)^{\gamma},\qquad \gamma> 1,
\end{align}
for large values of the exponent $\gamma$.
In this formula $V_{\mathrm{ref}}>0$ is a reference velocity, to bear in mind the physical meaning of $p$ (in the numerical simulations below 
we will simply set $V_{\mathrm{ref}}=1$).
This approach is used in fluid mechanics, for modeling 
certain free boundary problems where bubbles are immersed in a gas \cite{LiMa}.
This function is defined on $[0,\infty)$, it behaves proportionally to 
the gas-law $\rho^\gamma$ for $\rho \to 0 $ and it blows up as $\gamma\to +\infty$ for $\rho\geq \rho_\star$. 
Using \eqref{pr2} and the regime of large $\gamma$'s in traffic flows modeling is quite new; we shall see that this simple law has certain advantages in the numerical simulations of congested situations.
\medskip

In what follows, $p$ refers to \eqref{pr1}, \eqref{pr1bis} or \eqref{pr2}.
We will see that similar behaviors, corresponding to what can be expected for \eqref{PGCD}, are captured asymptotically (namely as $\ve\to 0$ or $\gamma\to \infty$) by these 
velocity offsets. However, 
the intermediate 
behaviors can significantly differ and the definition of the velocity offset seriously impacts the numerical costs.


\subsection{Stability issues}

As long as the functions $\rho$ and $v$ are smooth enough, we can rewrite system~\eqref{AW1} in the fully non--conservative form
\[
\partial_t \begin{pmatrix}\rho\\ v \end{pmatrix}
+A(\rho,v)\partial_x \begin{pmatrix}\rho\\ v \end{pmatrix}=0,\qquad
A(\rho,v)= \begin{pmatrix}v \quad & \rho\\ 0\quad & v-\rho p'(\rho) \end{pmatrix}.
\]
The two eigenvalues related to the system are therefore equal to 
\begin{align}\label{eig}
\lambda_1=v-\rho p'(\rho) \leq \lambda_2=v
\end{align}
 with related eigenvectors
\begin{align*}
	r_1=\begin{pmatrix}1\\- p'(\rho)\end{pmatrix},	\qquad
	r_2=\begin{pmatrix}1\\0\end{pmatrix}.
\end{align*}
The system is strictly hyperbolic, away from the regions where $\rho=0$.
Let us just note that the largest eigenvalue $\lambda_{2}$ is always linearly degenerate, leading to contact discontinuities and that $\lambda_{1}$ is genuinely non-linear, except for certain forms of the velocity offset we are not considering here. Therefore, the first eigenvalue will admit shocks or rarefaction waves. One of the difficulties of the computations is that vacuum regions may appear.
Observe that the information does not travel faster than the actual cars speed $v$, and that the system preserves the natural properties $\rho\geq 0$, $v\geq 0$.

We are considering here some Finite Volume (FV) numerical schemes in order to compute the solutions of system~\eqref{AW1}. Let us denote by $\Delta t$ and $\Delta x$ the time step and the space step of the method, respectively. 
We consider the discrete times $\ds t_{n}=n \Delta t$, for $ n \in \N$ and the discretization cells $\ds\mathcal{C}_{j}=[x_{j-1/2},x_{j+1/2}]$, $ j \in \Z$ (neglecting for the time being the issue of the boundary conditions) where
 $\ds x_{j+1/2}=(j+1/2) \Delta x$. 
We go back to the conservative form \eqref{AW1b} and we denote
\[U(x,t)=\begin{pmatrix} \rho(x,t)\\
y(x,t)\end{pmatrix} \textrm{ with } y(x,t)=\rho (v+p(\rho))(x,t),\]
the conservative variables. 
In terms of the conservative variables $\rho$ and $y$, we simply have
\[
\left\{
\begin{array}{l}
\partial_t \rho+\partial_x(\rho v)=0,
\\
\partial_t y+\partial_x (y v)=0,
\end{array}
\right.
\]
which recasts as follows, using only the variables $\rho$ and $y$, 
\begin{align}\label{AW3}
\left\{
\begin{array}{l}
\partial_t \rho+\partial_x(y-\rho p(\rho))=0,
\\
\partial_t y+\partial_x \left(\ds\frac{y^2}{\rho}-y p(\rho)\right)=0.
\end{array}
\right.
\end{align}
The numerical unknown $U^n_j=(\rho_{j}^n, y_{j}^n)$ 
is thought of as an approximation of the mean value of $U(x,t)$ on the cell $\mathcal{C}_{j}$ at time $t_{n}$.
The FV scheme has the following general form
\begin{align*}
	U^{n+1}_j=U^n_j
	-\ds\frac{\Delta t}{\Delta x} \pare{F^n_{j+1/2}-F^n_{j-1/2}}
\end{align*}
which mimics what we obtain by integrating the continuous equation \eqref{AW3} over $\crochets{t^n,t^{n+1}}\times \mathcal C_j$. 
 For the simple schemes we wish to deal with, the numerical flux
at the interface $x_{j+1/2}$ is a function of the neighbouring cells $F_{j+1/2}^n=\mathbb F\pare{U^n_{j+1},U^n_j}$.
Without entering into the details of the schemes, the numerical stability of such a method relies (at least) on the following inequality, see for example \cite[Sect.~2.3.3]{Bouc} or \cite{toro},
\begin{align}\label{stab}
\Delta t \leq \frac12 \frac{\Delta x}{\max(|\lambda_{1}|, |\lambda_{2}|)},
\end{align}
where $\lambda_{1}$ and $\lambda_{2}$ are the two eigenvalues given by \eqref{eig}.
It imposes that the state $U\pare{x_{j+1/2},t^{n+1}}$ on the interface $x_{j+1/2}$ at time $t^{n+1}=t^n+\Delta t$ 
only depends on the states of the unknown at time $t^n$ on the neighbouring cells $\mathcal C_j$ and $\mathcal C_{j+1}$.
It can thus be obtained by solving the corresponding Riemann problem with data $U^n_j$ and $U^n_{j+1}$.

Let us first consider the case of the pressure \eqref{pr1}. 
In case of a congestion formation, $\ds \rho^{\ve} \to\rho_{\star}$ but we expect that 
 $\ds p^{\ve}(\rho^{\ve})$ remains bounded and admits the limit $ \pi$ as $\ve \to 0$; it leads to the ansatz
 $$\rho=\rho_{\star}- O\pare{\ve^{1/\gamma}}, \text{ when } \ve \to 0.$$ 
Accordingly the behavior of the characteristic speeds is given by 
$$\max\pare{\abs{\lambda_{1}^{(\ve)}}, \abs{\lambda_{2}^{(\ve)}}}=O\left(\ve^{-1/\gamma}\right), \text{ when }\ve \to 0, $$
since $v^{\ve}$ should remain bounded when $\ve \to 0$. Hence, as $\ve$ goes to zero, the time step $\Delta t$ shrinks 
due to the condition 
\eqref{stab}.
The same remarks apply to the velocity offset \eqref{pr1bis}, which essentially behaves like \eqref{pr1} when $\ve \to 0$.
For the velocity offset \eqref{pr2}, we find 
$$\max\pare{\abs{\lambda_{1}^{(\gamma)}}, \abs{\lambda_{2}^{(\gamma)}}}
=O\pare{\gamma}, \text{ when } \gamma \to +\infty,$$
which again imposes tiny time steps.
This observation motivates the design of a scheme based on splitting strategy 
so that the fast waves can be treated implicitly.


\subsection{Invariant regions}

Let us detail another difficulty which is very specific to the traffic flow system \eqref{AW1}.
The Riemann invariants for the system \eqref{AW1} are given by, see \cite{AwR1},
$$ z_{1}=v+p(\rho), \qquad z_{2}=v.$$
Therefore, the domain 
$$\big\{ (z_1,z_2) \in \r^2 \text{ with } z_1 \in [w_{m},w_{M}], \, z_2 \in [v_{m}, v_{M}]\big\} $$
is an invariant region for \eqref{AW1}:
if the initial datum lies in such a region, the solution will still be contained in the same region for all times. 
However, numerical difficulties arise due to the fact that such domains are non--convex for the conserved quantities $\rho, y$.
This point has been observed in \cite{ChG} for the traffic flows model, see also \cite{CG08, Ju} for similar problems.
At first sight, it would be tempting to define the numerical 
fluxes by using the Godunov scheme, which is a standard for systems of conservation laws. It works into two steps. 
Owing to the stability condition \eqref{stab}, we solve a set of uncoupled Riemann problems, centered at the interfaces $x_{j+1/2}$ with data $U_j^n,U^{n}_{j+1}$.
Then, we project the obtained piecewise constant solution to obtain the updated numerical unknown, constant on the cells $\mathcal C_j$.
This projection step does not preserve the invariant region, since the latter is non--convex (for the role of the convexity of the invariant domain we refer the reader to
 \cite{AYLR} for gas dynamics equations, and more generally to the textbook
 \cite[Prop.~2.11]{Bouc}).
A counter--example is detailed in \cite{ChG} to explain why the Godunov scheme fails to satisfy the maximum principle for the Riemann invariants of \eqref{AW1}, and especially for $v$.
To solve this difficulty, \cite{ChG} proposes a hybrid scheme 
which mixes the Glimm scheme to compute the contact discontinuities, and the Godunov procedure to compute 1-shock or 1-rarefaction wave. 
This hybrid method is 
well-adapted to handle the specific velocity offsets dealt with in \cite{ChG, CG08}, 
see Remark~\ref{Rm2} below,
which differ from the models we wish to consider here.

We bear in mind that, instead of using the mean of the solutions of the Riemann problems over the cells, the Glimm scheme 
uses a random sampling strategy in the reconstruction procedure. 
Hence, by construction, the Glimm scheme preserves the invariant regions, despite the defect of convexity. It is thus well adapted to the simulation of the system \eqref{AW1}.
Note however the final scheme is non conservative.

\begin{remark}
Note that, depending on the definition of the numerical fluxes, the stability condition can be even more 
constrained than \eqref{stab}, for instance in order to fulfill the bound from above on the density with the ``close--packing--like'' velocity offset \eqref{pr1}, see e.g.~\cite{BGM}. 
\end{remark}

\begin{remark}\label{Rm2}
In \cite{ChG,CG08}, the discussion focuses on the velocity offset 
\begin{align}\label{prCG}
p(\rho)=V_{\mathrm{ref}}\ \ln\pare{\ds\frac\rho {\rho_\star}},
\end{align}
which has some very specific features:
\begin{itemize}\item First of all, $\rho\mapsto p(\rho)$ is 
defined on $(0,\infty)$ and 
non decreasing. The model cannot
 treat vacuum regions since the velocity offset is not defined for $\rho=0$.
 With this model, the velocity offset 
is well--defined beyond $\rho_\star$, 
and we note that $p(\rho)<0$ for $0<\rho<\rho_\star$, which might be physically questionable.
\item Second of all, we have 
$ \rho p'(\rho)= V_{\mathrm{ref}}$.
In particular, we get $\lambda_1=v- V_{\mathrm{ref}}$
so that the CFL condition depends only on $v$ and it does not shrink
as the density approaches $\rho_\star$. 
\end{itemize}
\end{remark}

\section{Description of the scheme}\label{scheme}

Let us now explain in more details the construction of the scheme that we wish to use for the simulation of \eqref{AW1}, with a velocity offset $\rho\mapsto p(\rho)$
that introduces some stiffness in order to reproduce the expected behavior of the constrained system \eqref{PGCD}.
In what follows, $p$ thus refers to $p^{\ve}$ in \eqref{pr1}, to $\tilde p^{\ve}$ in \eqref{pr1bis} or to $p^{\gamma}$ in \eqref{pr2}.
\\

In order to get rid of the large characteristic speeds, 
the idea consists in splitting the velocity offset into two parts
\[p=p_{exp}+p_{imp},\]
so that the system with $p_{exp}$ 
has a stability CFL condition \eqref{stab} of order 1 with respect to the scaling parameters. 
Namely, the eigenvalue $\lambda_{exp,1}^{(\ve)}$ (resp. $\lambda_{exp,1}^{(\gamma)}$) is bounded with respect to $0<\ve\ll 1$ (resp. $\gamma\gg 1$).
The corresponding system can thus be treated explicitly by means of the Glimm scheme, 
which preserves the invariant domains.
Next, only the stiff part that involves $p_{imp}$ 
 is treated implicitly.
 We expect also that the implicit part has a simple structure that can be handled with a not too complicated scheme.
 Such splitting approach appeared in \cite{DH, DHN, DeTa} for more standard fluid mechanics systems, for instance as an efficient strategy to handle low Mach number regimes.
 However, as explained above, the structure of the PDE system 
 \eqref{AW1} significantly differs from the Euler equations, in particular with the lack of convexity of the invariant domains
 and these methods cannot be directly applied to \eqref{AW1}.

 \subsection{Definition of the explicit velocity offset}
 
In the first step of the splitting, we consider the system
\begin{align}\label{AW2p0}\left\{\begin{array}{l}
\partial_t \rho+\partial_x(\rho v)=0,
\\
\partial_t \big(\rho(v+p_{exp}(\rho))\big)+\partial_x \big(\rho v(v+p_{exp}(\rho))\big)=0.
\end{array}\right.
\end{align}
 It has the same structure as the Aw-Rascle-Zhang system \eqref{AW1b}, just replacing the full velocity offset $p$, that can be \eqref{pr1}, \eqref{pr1bis} or \eqref{pr2}, by $p_{exp}$.
The characteristic speeds are 
 \[\lambda_1=v-\rho p'_{exp}(\rho),\qquad \lambda_2=v.\]
 In order to relax the stability constraint, we define $p_{exp}$
so that the characteristic speed does not blow up as the scaling parameters $\ve$ goes to 0 or $\gamma$ goes to $\infty$.
 It leads to require that $p'_{exp}(\rho)$ is bounded uniformly with respect to $\varepsilon$ (resp. $\gamma$ for the law \eqref{pr2}), when $\rho$ lies in a compact set of $(0,\infty)$.
 The definition depends on a truncation parameter $0<\rho_{num}<\rho_\star$.
Following an idea of \cite{DH, DHN}, we set
 \begin{align}\label{pr0}
 p_{exp}(\rho)=
\left\{
\begin{array}{l}
p(\rho), 
\hspace*{7cm} \text{ if } 0 \leq \rho \leq \rho_{num},
\\
p(\rho_{num})+ p'(\rho_{num})(\rho-\rho_{num})+ 
p''(\rho_{num})\ds\frac{(\rho-\rho_{num})^2}{2}, \\
 \hspace*{8cm}\text{ if } \rho>\rho_{num}.
\end{array}
\right.
\end{align}
We use a second order polynomial for $ \rho>\rho_{num}$ to ensure that $p_{exp}$ is a $\ds \mathcal{C}^2$ function.
We will see below that it is important to reach such a regularity. 
The transition density $\rho_{num}$ is chosen so that $p_{exp}'(\rho_{num})$ is bounded when $\ve \to 0$ (resp. $\gamma \to +\infty$). If such a condition is satisfied, then the two eigenvalues are bounded with respect to $\ve$ (resp. $\gamma$) and therefore the stability condition \eqref{stab} is independent of $\ve$ (resp. $\gamma$). Let us now explain how to choose $\rho_{num}$ for the velocity offsets under consideration. 
 \begin{itemize}
\item[a)] For the laws \eqref{pr1} and \eqref{pr1bis}, when $\ve \to 0$ and $\rho^{\ve}\to \rho_{\star}$, we expect that $ \rho^{\ve}=\rho_{\star} - O(\ve^{1/\gamma})
$, since $\ds p^\ve (\rho^{\ve})$ remains finite.
We set $\ds\rho_{num}=\rho_{\star}(1-\delta \rho)$ with 
 $\delta \rho>0$. A simple computation shows that
$$p_{exp}'(\rho_{\star})= p'(\rho_{num})+ p''(\rho_{num})(\rho_{\star}-\rho_{num})=O\left(\ve (\delta \rho) ^{-(\gamma+1)}\right),$$
which leads us to set 
$$\delta \rho=\ve^{1/(\gamma+1)} \text{ and }\rho_{num}=\rho_{\star}\left(1-\ve^{1/(\gamma+1)}\right).$$
For the velocity offset \eqref{pr1bis}, we point out that $\rho_{num}$ is a truncation parameter 
of numerical nature, designed to ensure a certain stability property of the scheme, while $\rho_{trans}$ relies on modeling consideration and does not prevent the blow up of the velocity offset.
In any case, we have $0<\rho_{num}<\rho_{trans}<\rho_\star$.
\item [b)]
For the law \eqref{pr2}, we require that $p_{exp}'(\rho_{num})$ remains bounded when $\gamma \to +\infty$.
 Denoting again $\ds\rho_{num}=\rho_{\star}(1-\delta \rho)$, a simple computation leads to
 $$p'_{exp}(\rho_{num})=\frac{\gamma}{\rho_{\star}}\left( 1-\delta \rho\right)^{\gamma-1} +\frac{\gamma(\gamma-1)}{\rho_{\star}^2}\delta \rho\left( 1-\delta \rho\right)^{\gamma-2} .$$
 So, we are looking for $\delta \rho$ such that $\delta \rho \to 0$ and both terms in this sum are $O(1)$ when $\gamma \to +\infty$. A simple study of the two sequences of functions $\ds\left( \gamma \exp\left( -\gamma x\right) \right)_{\gamma \in \N}$ and $\ds \left( \gamma^2 x \exp\left( -\gamma x\right) \right)_{\gamma \in \N}$ shows that we should have $\ds \delta \rho=O(\gamma^{-\alpha})$ with $\alpha \in (0,1)$ in order to satisfy the required properties. 
\end{itemize}

In \cite{DeTa} the pressure term is split as 
$p=p_{exp}+p_{imp}$ 
with 
$p_{exp}=\alpha p$ 
and 
$p_{imp}=(1-\alpha)p.$
Here, due to the singularity of $p$ at $\rho_{\star}$, this approach does not permit to keep $p'_{exp}$ bounded and we have to define $p_{exp}$ and $p_{imp}$ as explained above.

 \subsection{A time-splitting scheme}

As said above, it is convenient to work on the conservative form \eqref{AW3} of the system \eqref{AW1}, 
dealing with the unknowns $U=\big(\rho,y\big)$ where $y=\rho(v+p(\rho))$.
With $\ds p=p_{exp}+p_{imp}$, we arrive at 
\begin{align*}
\begin{cases}
	\partial_t \rho+\partial_x \pare{y-\rho p_{exp}(\rho)-\rho p_{imp}(\rho)}=0,		\\
	\partial_t y+\partial_x\pare{ \ds\frac{y^2}{\rho} -y p_{exp}(\rho) -y p_{imp}(\rho)}=0.
\end{cases}
\end{align*}
We use now a time-splitting scheme. 
Knowing some approximate values $U^n=(\rho^n,y^n)$ at time $t_{n}$, we proceed as follows.
\begin{itemize}
\item Step 1: Solve with an explicit scheme the system of conservation laws
\begin{align*}
\begin{cases}
	\partial_t \rho+\partial_x \pare{y-\rho p_{exp}(\rho)}=0,	\\
	\partial_t y+\partial_x\pare{ \ds\frac{y^2}{\rho} -y p_{exp}(\rho)}=0.
\end{cases}
\end{align*}
As said above, this system has the same structure as the original problem \eqref{AW3}.
In particular the invariant domains are non--convex. 
It can be solved with the Glimm scheme adapted for the pressure $p_{exp}$.
More details will be given in Section \ref{exp}.
It defines some intermediate values $(\rho^{n+1/2}, y^{n+1/2})$. 
\item Step 2: Solve implicitly the system
\begin{align}
\label{imp_eq}
	\begin{cases}
		\partial_t \rho-\partial_x \pare{\rho \, p_{imp}(\rho)}=0,	\\
		\partial_t y-\partial_x\pare{ y \, p_{imp}(\rho)}=0.
	\end{cases}
\end{align}
Note that the system has a simple structure and the two equations decouple.
 The first equation is a non linear scalar conservation law for the density $\rho$, the second is a linear transport equation for $y$ where the velocity $-p_{imp}(\rho)$ can be considered as given. 
 The numerical method to solve the system \eqref{imp_eq} is thus not that complicated.
 More details will be given in Section \ref{imp}.
\end{itemize}


\subsection{A few words about the Glimm scheme for \eqref{AW1}}
\label{exp}

In order to use a Glimm scheme, we need to know the Riemann solutions of the problem. This computation has already been done in \cite{AwR1} 
and in \cite[Section~6]{BeDe} where all the details can be found. We 
refer the reader to some classical books \cite{smoller,toro} for general discussions about the role of Riemann problems in the theory of conservation laws
 and to \cite{AwR1, BeDe} for the specific case of the traffic flow system. We recap here only the Riemann solutions, omitting the details on the elementary waves and on the admissibility of solutions.

\subsubsection{A brief overview on the Riemann problem for \eqref{AW1}}

Let us just recall that the second eigenvalue $\lambda_{2}$ given by \eqref{eig} of system \eqref{AW1} is always linearly degenerate, leading to contact discontinuities and that $\lambda_{1}$ is genuinely non-linear, leading to shocks or rarefaction waves. One of the difficulties of the computations is that vacuum regions may appear. Therefore, the Riemann solution of system \eqref{AW1}, with an initial datum
\begin{align*}
(\rho,v)(x,0)=
\left\{
\begin{array}{ll}
\pare{\rho^L,v^L}, & \text{ for } x<0,
\\
\pare{\rho^R,v^R}, & \text{ for } x>0,
\end{array}
\right.
\end{align*}
can be computed according to the five following cases:
\begin{itemize}
\item if $\ds \rho^L>0,\;\rho^R>0$ and $\ds 0\leq v^R\leq v^L$, the solution consists of a 1-shock that connects $\pare{\rho^L,v^L}$ to the intermediate state $(\rho^*,v^*)$ and a contact discontinuity between $(\rho^*,v^*)$ and $\pare{\rho^R,v^R}$;
\item if $\ds \rho^L>0,\;\rho^R>0$ and $\ds v^L< v^R\leq v^L+ p(\rho^L)$, the solution consists of a 1-rarefaction wave that connects $\pare{\rho^L,v^L}$ to $(\rho^*,v^*)$ and a contact discontinuity between $(\rho^*,v^*)$ and $\pare{\rho^R,v^R}$;
\item if $\ds \rho^L>0,\;\rho^R>0$ and $\ds v^L+ p\pare{\rho^L}< v^R$, a vacuum region appears; the solution consists of a 1-rarefaction wave that connects $\pare{\rho^L,v^L}$ to $(0,v^*)$, then a vacuum region between $(0,v^*)$ and $(0,v^R)$ and a contact discontinuity between $\pare{0,v^R}$ and $\pare{\rho^R,v^R}$;
\item if $\ds \rho^L>0,\; \rho^R=0$, the solution is only a 1-rarefaction wave connecting $\pare{\rho^L,v^L}$ and $\pare{0,v^R}$;
\item if $\ds \rho^L=0,\; \rho^R>0$, the solution is only a 2-contact discontinuity connecting $\pare{0,v^L}$ and $\pare{\rho^R,v^R}$.
\end{itemize}

The intermediate state $(\rho^*,v^*)$ of the Riemann solution is computed with the Riemann invariants, that is to say
\begin{align*}
\begin{cases}
	v^*=v_{R},	\\
	\rho^*=p^{-1}\pare{ v^{L}-v^{R}+p\pare{\rho^{L}}}.
\end{cases}
\end{align*}
In the case of a shock, the speed of the shock between $(\rho^*,v^*)$ and $\pare{\rho^L,v^L}$ is given by $$ s=\frac{\rho^*v^*-\rho^L v^L}{\rho^*-\rho^L}.$$
In the case of a rarefaction wave, the self-similar solution $(\rho,v)(\xi)$ with $\xi=x/t$ is given by the following formulae
\begin{align}
\label{ssim}
\begin{cases}
	p(\rho(\xi))+\rho(\xi) p'(\rho(\xi))= p\pare{\rho^L}+v^L-\xi,	\\
	v(\xi)=v^L+ p\pare{\rho^L}-p(\rho(\xi)),
\end{cases}
\end{align}
which apply for $ \xi \in \crochets{\lambda_{1}\pare{\rho^L,v^L}, \lambda_{1}\pare{\rho^R,v^R}}$.

\begin{remark}
In practice, we compute the self-similar solutions of equation \eqref{ssim}, 
by using the Newton algorithm.
The method requires that $p$ has the $\ds \mathcal{C}^2$ regularity. 
This remark explains the construction of the explicit part of the velocity offset \eqref{pr0} (we have observed 
bad behaviors of the scheme 
when $p''$ has jumps).
\end{remark}

\subsubsection{Glimm scheme}

Hence, we have at hand formula to compute the solution of the Riemann problems, which are the elementary brick of the Glimm's scheme (like for Godunov's scheme).
This scheme has been introduced for theoretical purposes \cite{Gli}, 
and its implementation for hyperbolic systems is further discussed in \cite{ChG, colella, liu}. 
 Let $U^n_j=\ds \pare{\rho_{j}^n, y_{j}^n}$ be the approximated mean value on the cell $\mathcal{C}_{j}$ of $U=(\rho,y)$ at time $t_{n}$.
We proceed as follows to update the numerical unknown:
 \begin{itemize}
 \item 
We solve the associated Riemann problem at each interface $x_{j+1/2}$, namely all the Riemann problems with $U_L=U^n_j$, $U_R=U^n_{j+1}$.
\item 
Let $a_n$ be a number picked randomly in $ [0,1]$. We define the value $U^{n+1}_j$ of the numerical unknown in the cell $\mathcal C_j$ at time $t^{n+1}$ 
as to be the solution of the Riemann problem evaluated at the point $\ds x_{j-\frac{1}{2}}+a_n\Delta x \in \mathcal{C}_{j}.$
The scheme does not use any averaging or projection procedure and the obtained solution, by construction, remains in the invariant region of the PDE system. 
 In practice, we use the Van Der Corput quasi--random sequence $(a_n)_{n\in\N}$ (see \cite{colella}) defined by
$$
	a_n=\ds\sum_{k=0}^m i_k2^{-(k+1)}
$$
where $ n=\sum_{k=0}^m i_k2^k$, with $i_k\in\{0,1\}$, denotes the binary expansion of the integer $n$. 
\end{itemize}


\subsection{Treatment of the implicit part}\label{imp}

Let us now discuss how we handle the system \eqref{imp_eq} where we remind the reader that $p_{imp}$ contains the stiff part of the velocity offset. 
As said above, the system decouples and it has a very simple structure.
Let us set $\Phi:\rho \in [0,\infty)\mapsto - \rho p_{imp}(\rho)\in (-\infty,0]$. 
We solve the scalar conservation law for $\rho$ 
with the classical Engquist--Osher scheme. The Engquist--Osher flux is defined by
$$\mathbb F \pare{\rho_{j+1},\rho_j}=\ds\int_{\mathbb R} [\Phi']_+(\xi)\mathbf 1_{0\leq \xi\leq \rho_{j}}\ud \xi+
\ds\int_{\mathbb R} [\Phi']_-(\xi)\mathbf 1_{0\leq \xi\leq \rho_{j+1}}\ud \xi,$$
where $[X]_{+}=\max(X,0)$ and $[X]_{-}=\min(X,0)$.
The implicit scheme takes the following form
\begin{align}
\label{eo}
	\rho^{n+1}_j=\rho^{n+1/2}_j
	-\ds\frac{\Delta t}{\Delta x}\pare{\mathbb{F} \pare{\rho^{n+1}_{j+1},\rho^{n+1}_j}-\mathbb{F} \pare{\rho^{n+1}_{j},\rho_{j-1}^{n+1}}},
\end{align}
where $\rho^{n+1/2}_j$ is the result from the first (explicit) step of the scheme.

 Here, the velocity field $\Phi'(\rho)=-\frac{\ud }{\ud \rho}\big(\rho p_{imp}(\rho)\big) = -p_{imp}(\rho) - \rho p_{imp}'(\rho)$ of this scalar equation is always non positive.
Accordingly, 
the numerical flux reduces to a mere function of the right density
$$\mathbb F (\rho_{j+1}, \rho_j)=\ds
\ds\int_{\r} \Phi'(\xi)\mathbf 1_{0\leq \xi\leq \rho_{j+1}}\ud \xi=\Phi( \rho_{j+1}).$$
Consequently, the non-linear equation \eqref{eo} becomes
\begin{align*}
	\rho^{n+1}_j=\rho^{n+1/2}_j-
	\ds\frac{\Delta t}{\Delta x}\pare{\Phi\pare{ \rho_{j+1}^{n+1}}-\Phi\pare{ \rho_{j}^{n+1}}}.
\end{align*}
It forms a triangular system
 of non linear scalar equations.
 If $J$ stands for the number of grid points, we have to solve $J$ non linear scalar equations, which means $J$ executions of the scalar Newton algorithm to update the density.
 Then a mere linear system defines the updated $y$.
Indeed, once $\rho^{n+1}$ is determined, we solve the transport equation for updating $y$. To this end, we use the standard implicit upwind scheme
\begin{align*}
	y^{n+1}_j=y^{n+1/2}_j
	-\ds\frac{\Delta t}{\Delta x}\pare{G^{n+1}_{j+1/2}- G^{n+1}_{j-1/2}} 
\end{align*}
with
\[ G^{n+1}_{j+1/2}= y^{n+1}_j \crochets{-p_{imp}\pare{\rho^{n+1}_j}}_{+} + y^{n+1}_{j+1} \crochets{-p_{imp}\pare{\rho^{n+1}_{j+1}}}_{-} .\]
The specific case $p_{imp}\geq 0$, yields
\[y^{n+1}_j=y^{n+1/2}_j-\ds\frac{\Delta t}{\Delta x}\big(-p_{imp}\pare{\rho^{n+1}_{j+1}} y^{n+1}_{j+1}+p_{imp}\pare{\rho^{n+1}_{j}} y^{n+1}_{j}\big). \]
It forms a triangular linear system of equations that can be solved by backward substitution, leading to the straightforward formula
\begin{displaymath}
	y_j^{n+1}=\ds \frac{\ds y_{j}^{n+1/2}+\frac{\Delta t}{\Delta x}p_{imp}\pare{\rho^{n+1}_{j+1}}y^{n+1}_{j+1}}
	{\ds 1+\frac{\Delta t}{\Delta x}p_{imp}\pare{\rho^{n+1}_{j}}}.
\end{displaymath}

\section{Simulations}
\label{simu}


As indicated in the introduction,
it is far from clear how to design a  ``natural'' scheme for
 a direct simulation of the constrained model \eqref{PGCD}.
 We can only mention the recent approach 
 for crowd dynamics proposed in \cite{Maury}; here we are rather motivated by the asymptotic issues. 
Our aim is two--fold. On the one hand we wish to discuss the asymptotic behavior of the different models 
\eqref{pr1}, \eqref{pr1bis} and \eqref{pr2} for the velocity offset, which are all expected to capture asymptotically (for $\ve\to 0$ or $\gamma\to \infty$) the features of the limit system \eqref{PGCD}.
On the other hand, we shall discuss the numerical difficulties and the ability of the time--splitting strategy, which will be compared to the standard Glimm scheme with a small enough time step, in handling the asymptotic behaviour.

The simulations presented below are thought of as Riemann problems and we impose boundary conditions
that maintain constant the inflow conditions.
Of course, the method can be adapted to treat further boundary conditions. 
In particular imposing zero--influx produces vacuum regions, a numerical difficulty  that our method is able to handle, as shown with the decongestion case below.

\subsection{Case of simple transport}
\label{transport}

\begin{figure}
\subfloat[Density - Glimm scheme]{
	\includegraphics[trim=0 200 0 200,clip=true,scale=0.35]
		{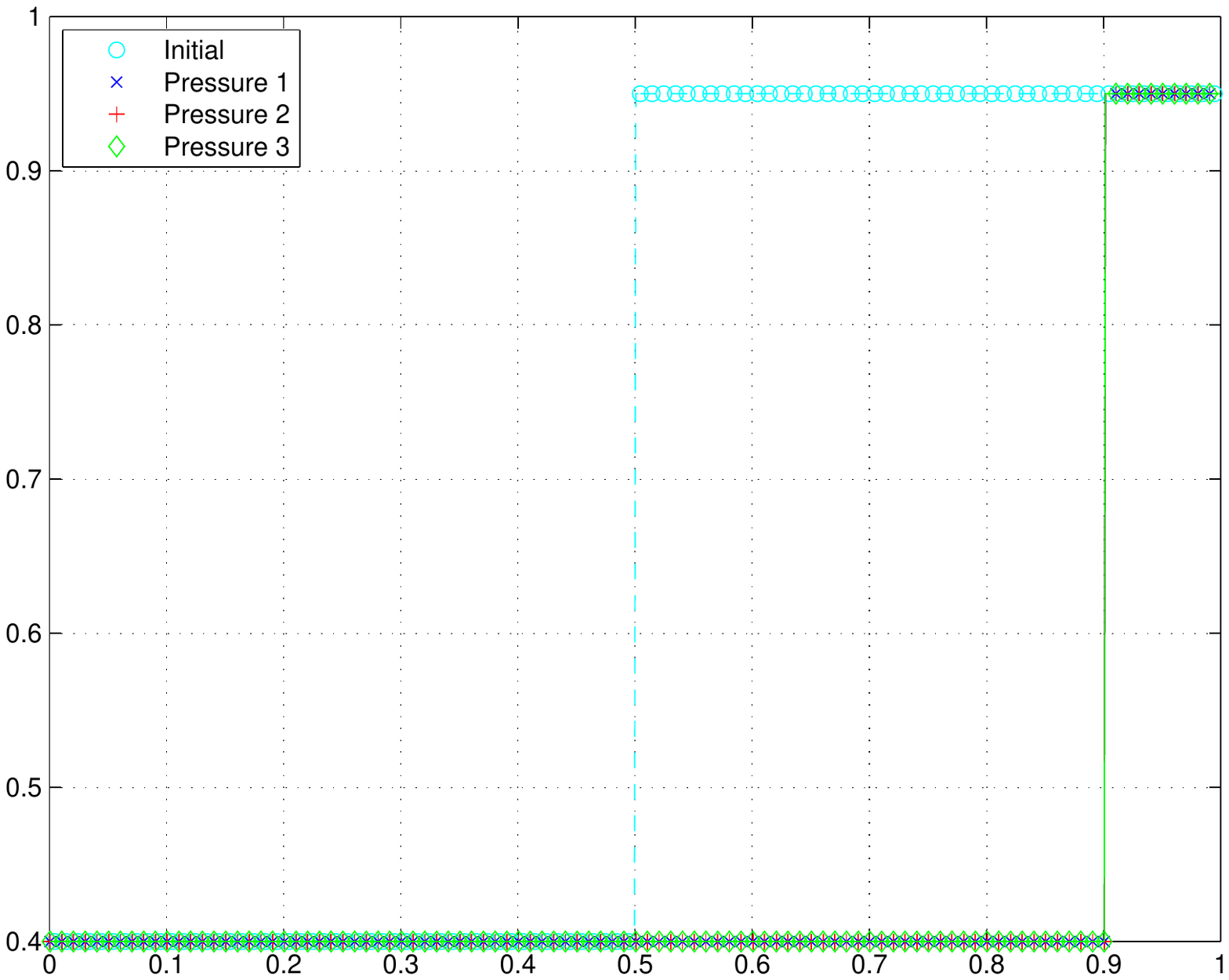}\label{Transp_11}
}
\subfloat[Velocity - Glimm scheme ]{
	\includegraphics[trim=0 200 0 200,clip=true,scale=0.35]
		{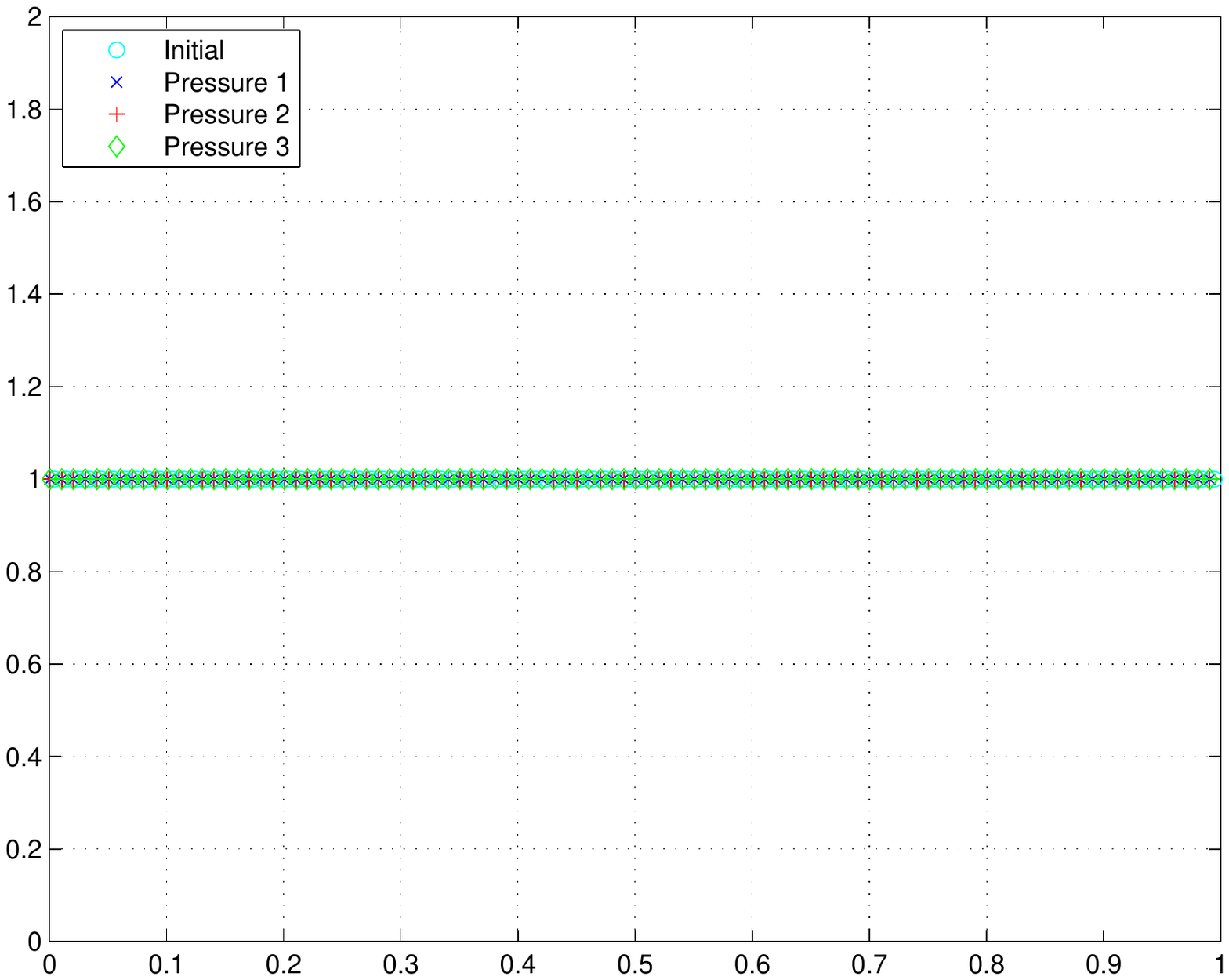}\label{Transp_12}
}
\\	
\subfloat[Density - implicit-explicit scheme ]{
	\includegraphics[trim=0 200 0 200,clip=true,scale=0.35]
		{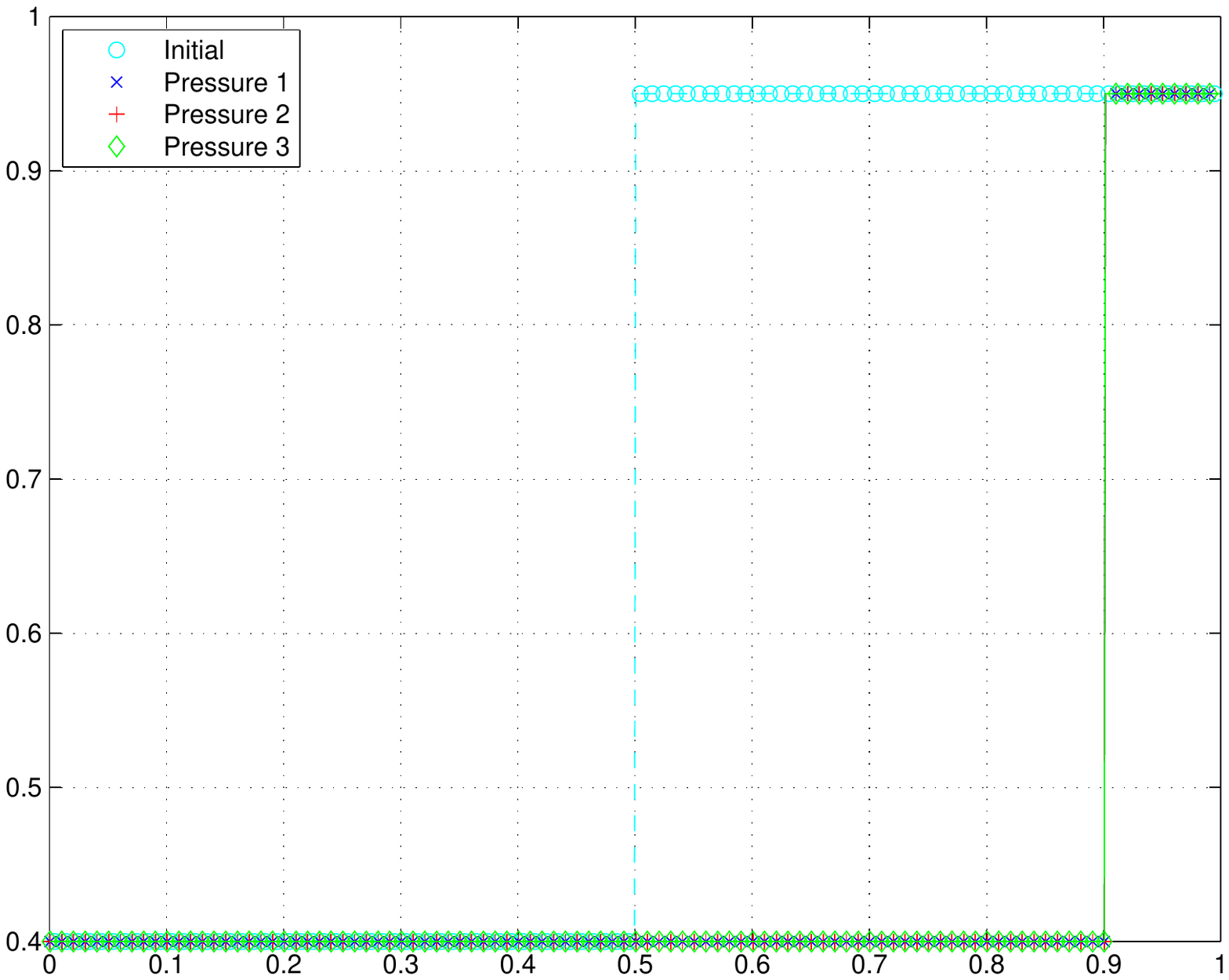}\label{Transp_21}
}
\subfloat[Velocity - implicit-explicit scheme]{
	\includegraphics[trim=0 200 0 200,clip=true,scale=0.35]
		{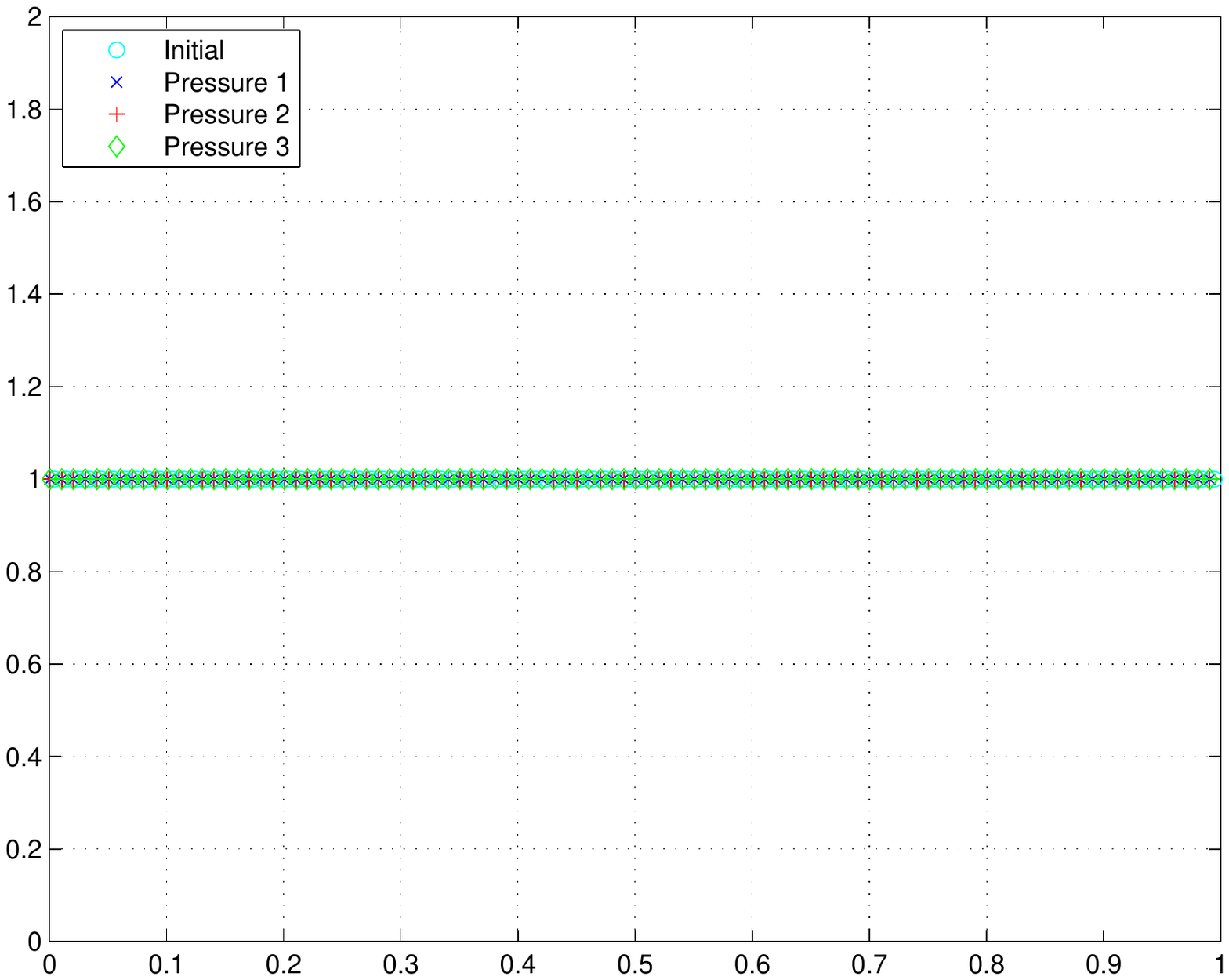}\label{Transp_22}
}
\caption{\textbf{Numerical results in the case of transport \eqref{Transp-init}--\eqref{Transp-exact}. }
Density (left) and velocity (right) with the Glimm scheme (top) and the explicit-implicit scheme (bottom). The results are given for the three different pressures under consideration: pressure \eqref{pr1} in blue, pressure \eqref{pr1bis} in red and pressure \eqref{pr2} in green. The initial conditions are plotted in cyan. \label{Transp} }
\end{figure}

To begin with, we test the case of a simple transport: 
the computational domain is the interval $[0,1]$ and for the initial data we set
\begin{align}
\label{Transp-init}
v^0(x)=1,\qquad 
\rho^0(x)=
\left\{\begin{array}{ll}
0.4, & \text{ if } x \in [0, 0.5[, \\
0.95, & \text{ if } x \in [0.5, 1],
\end{array}\right.
\end{align}
see Figure~\ref{Transp} (cyan curves).
We compare the six following situations:
\begin{itemize}
\item system \eqref{AW1} with pressure \eqref{pr1}, using the Glimm scheme,
\item system \eqref{AW1} with pressure \eqref{pr1}, using the scheme presented in Section~\ref{scheme},
\item system \eqref{AW1} with pressure \eqref{pr1bis}, using the Glimm scheme,
\item system \eqref{AW1} with pressure \eqref{pr1bis}, using the scheme presented in Section~\ref{scheme},
\item system \eqref{AW1} with pressure \eqref{pr2}, using the Glimm scheme,
\item system \eqref{AW1} with pressure \eqref{pr2}, using the scheme presented in Section~\ref{scheme}.
\end{itemize}
The solution at $T=0.4$ should be
\begin{align}\label{Transp-exact}
v(x,T)=1,\qquad 
\rho(x,T)=
\left\{\begin{array}{ll}
0.4, & \text{ if } x \in [0, 0.9[, \\
0.95, & \text{ if } x \in [0.9, 1].
\end{array}\right.
\end{align}
The space step is equal to $\Delta x=10^{-3}$ and the time step is computed in order to satisfy the stability condition (evaluated with the full $p$ for the Glimm scheme, and with $p_{\exp}$ for the implicit--explicit method).
 The parameters are taken as follows:
\begin{itemize}
\item Pressure \eqref{pr1} with $\gamma=2$, $\varepsilon=10^{-3}$, $\rho_{\star}=1$. For the explicit-implicit scheme, the numerical threshold is chosen as $ \rho_{num}=\rho_{\star}\left(1-\frac 1 5\ve^{1/(\gamma+1)}\right).$
\item Pressure \eqref{pr1bis} with $\gamma=2$, $\varepsilon=10^{-3}$, $\rho_{\star}=1$ and $\rho_{tr}^{\varepsilon}=\rho_{\star}-\varepsilon$. For the explicit-implicit scheme, the numerical threshold is chosen as $ \rho_{num}=\rho_{\star}\left(1-\frac 1 5 \ve^{1/(\gamma+1)}\right).$
\item Pressure \eqref{pr2} with $\gamma=4$. For the explicit-implicit scheme, the numerical threshold is chosen as $\ds \rho_{num}=\rho_{\star}\left(1-10^{-2}\right).$
\end{itemize}
The results of the numerical simulations for the three different pressures performed with the Glimm scheme are displayed at Figure~\ref{Transp_11}-\ref{Transp_12}, whereas the same simulations using the scheme constructed in Section~\ref{scheme} are exhibited at Figures~\ref{Transp_21}-\ref{Transp_22}.
All the results are equivalent and agree with the exact solution.

\subsection{Case of decongestion}

\begin{figure}
\subfloat[Density - Pressures \eqref{pr1} and \eqref{pr1bis}]{
\includegraphics[trim=0 200 0 200,clip=true,
 scale=0.35]
{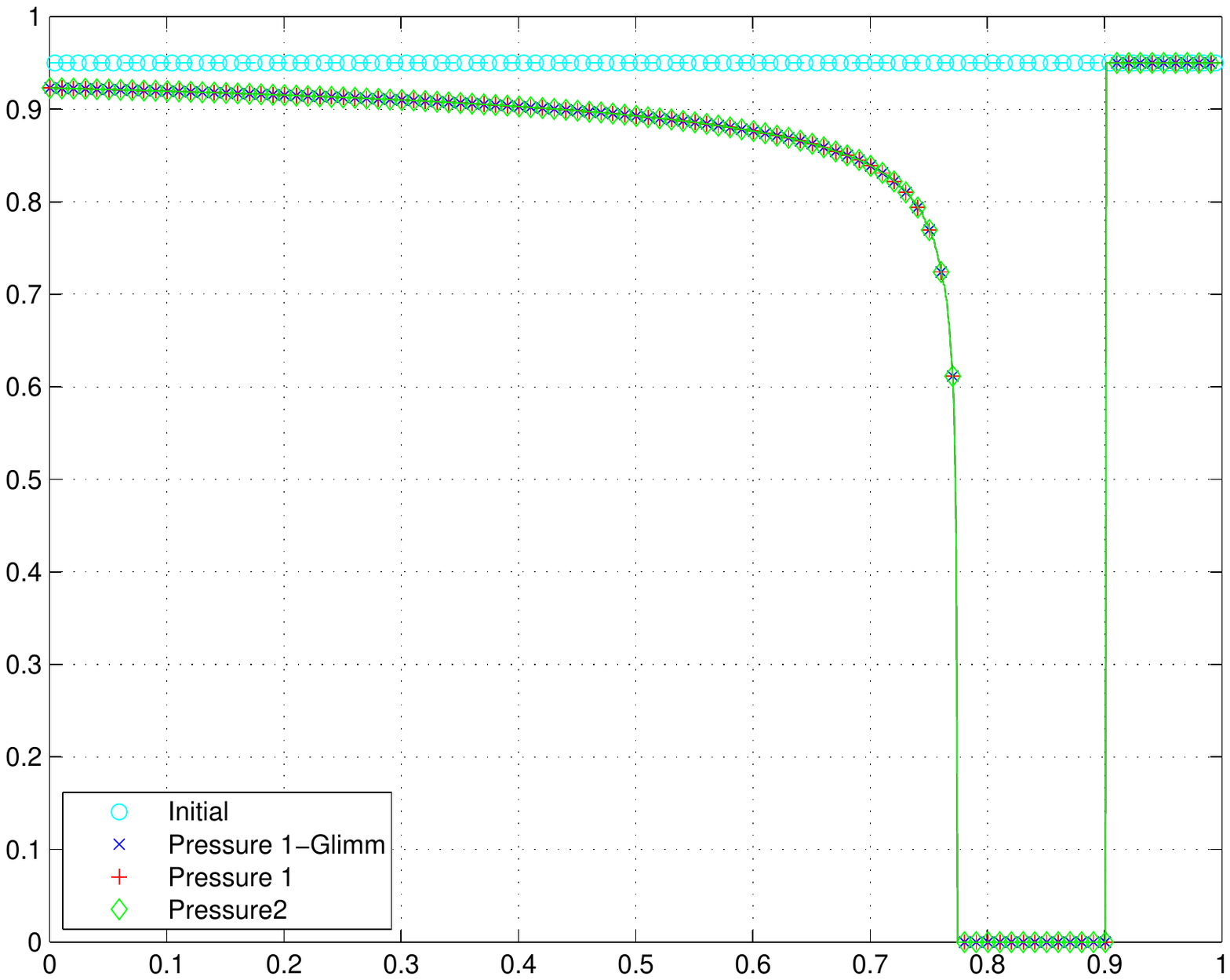}
}
\subfloat[Velocity - Pressures \eqref{pr1} and \eqref{pr1bis}]{
\includegraphics[trim=0 200 0 200,clip=true,
scale=0.35]
{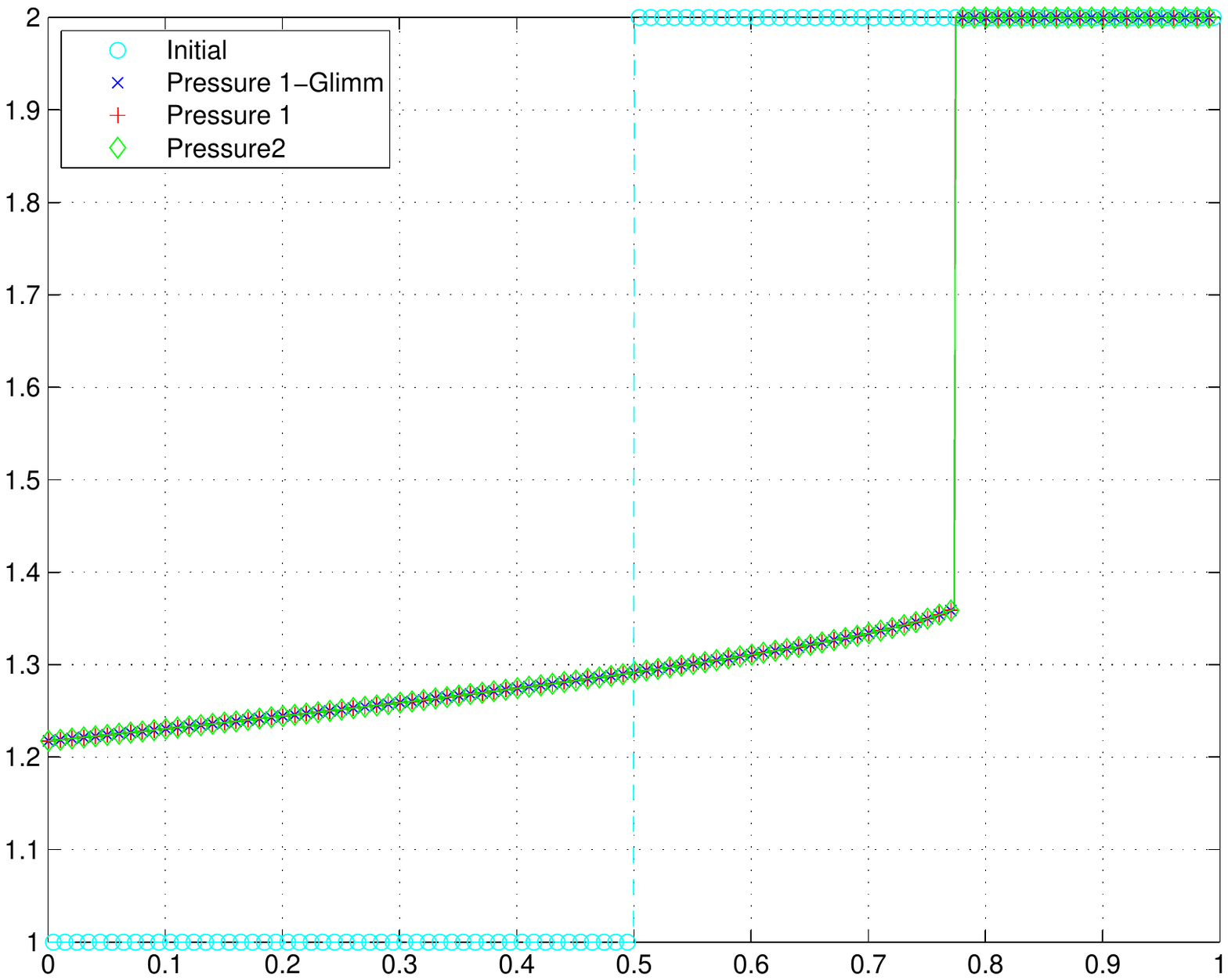}
}
\\
\subfloat[Density - Power-type pressure \eqref{pr2}.]{
\includegraphics[trim=0 200 0 200,clip=true,
scale=0.35]
{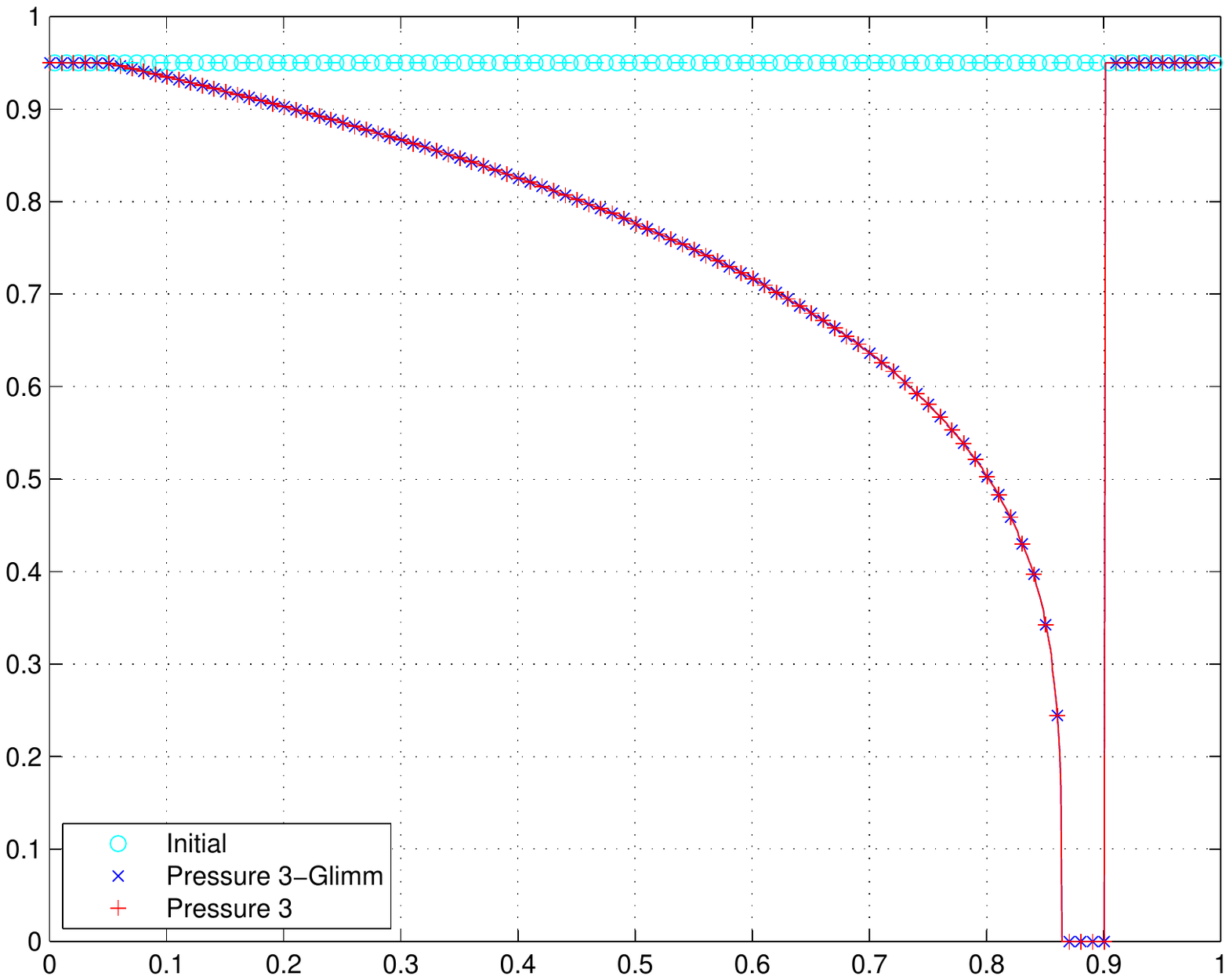}
}
\subfloat[Velocity - Power-type pressure \eqref{pr2}.]{
\includegraphics[trim=0 200 0 200,clip=true,
scale=0.35]
{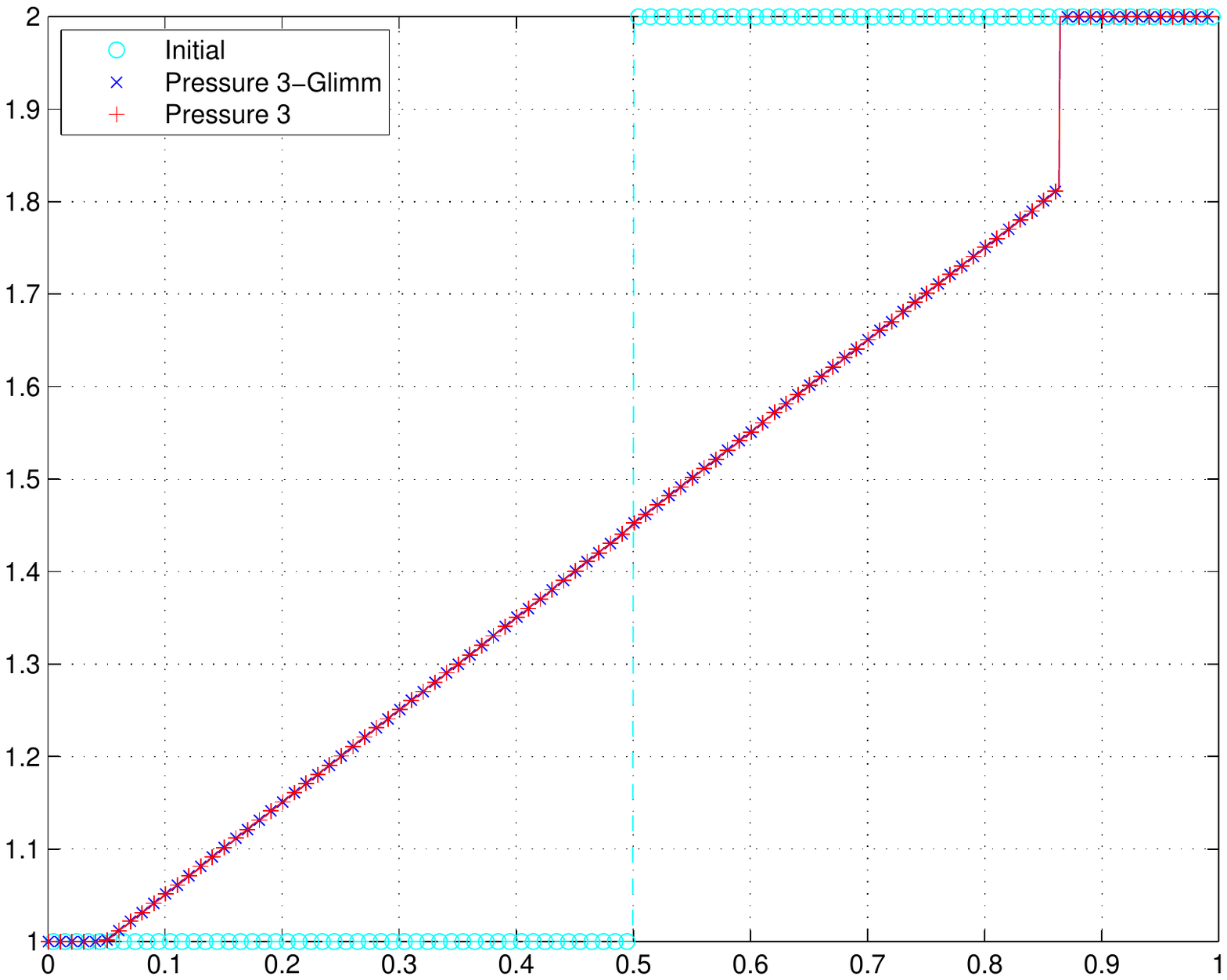}
}
\caption{\textbf{Numerical results in the case of decongestion \eqref{Decong-init}--\eqref{Decong-exact}.}
 Comparison of the two schemes. Density (left) and velocity (right) at final time $T=0.2$.
 Top: pressure \eqref{pr1} with Glimm scheme (blue),
 implicit-explicit scheme (red) and pressure \eqref{pr1bis} with implicit-explicit scheme (green). 
Bottom: pressure \eqref{pr2} with Glimm scheme (blue) and 
with the explicit-implicit scheme (red). The initial conditions are plotted in cyan. Parameters : $\gamma=2, \ve=10^{-3}$ for pressures  \eqref{pr1} and  \eqref{pr1bis}  and $\gamma=4$ for pressure \eqref{pr2}.
\label{Decong1}}
\end{figure}

\begin{figure}
\subfloat[Density - Pressure \eqref{pr2} for different $\gamma$]{
\includegraphics[trim=0 200 0 200,clip=true,
 scale=0.35]
{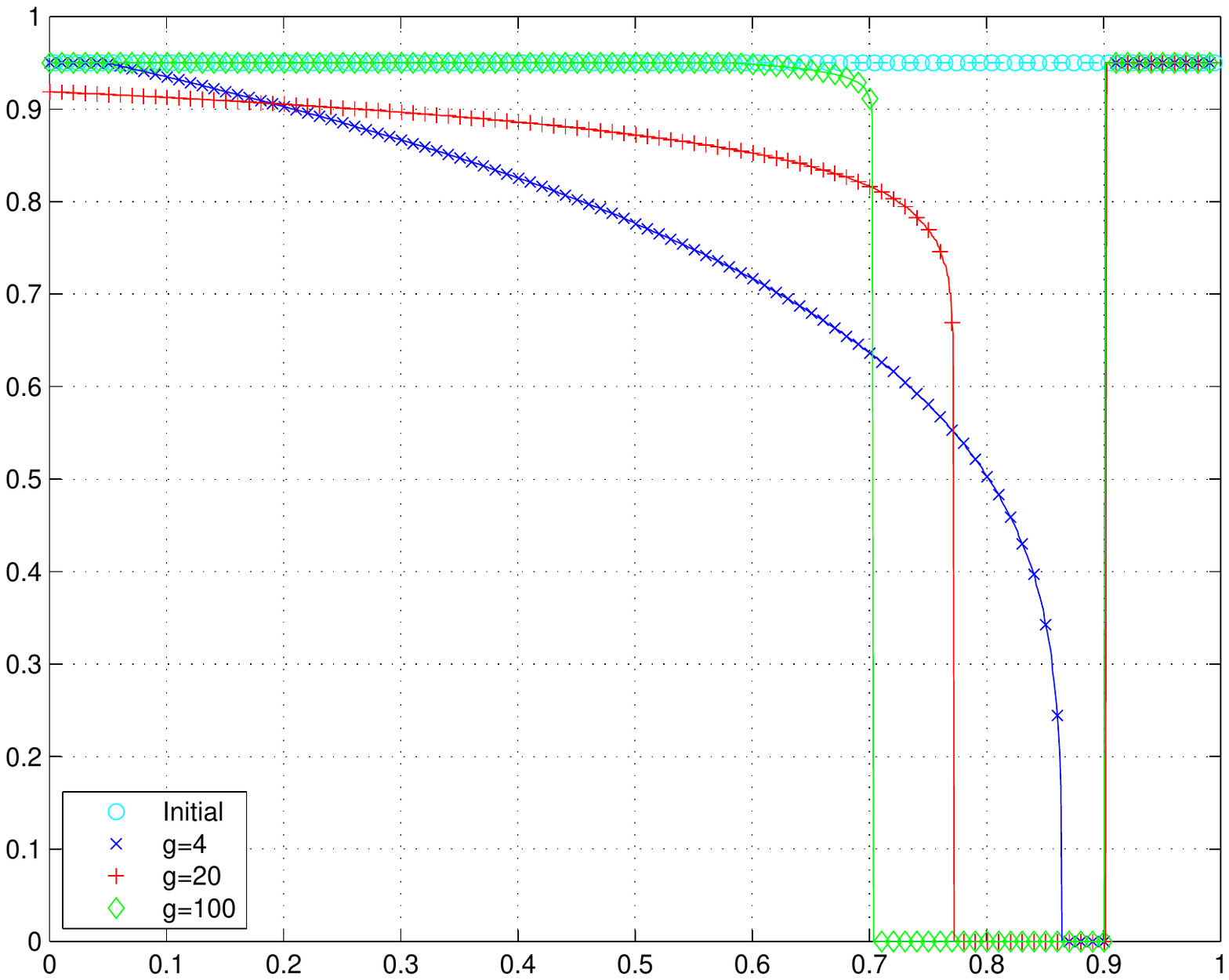}
}
\subfloat[Velocity - Pressure \eqref{pr2} for different $\gamma$]{
\includegraphics[trim=0 200 0 200,clip=true,
scale=0.35]
{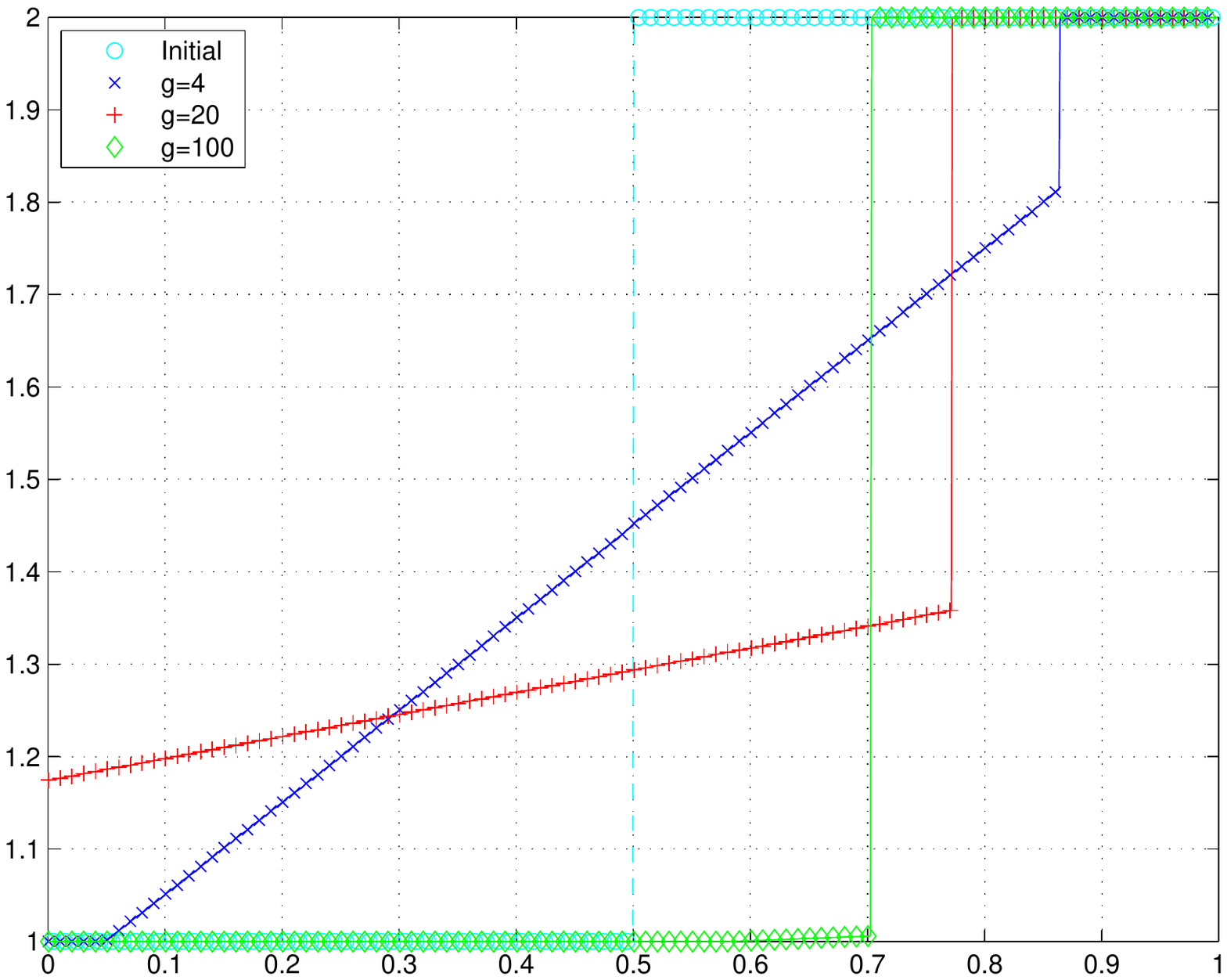}
}
\caption{\textbf{Numerical results in the case of decongestion \eqref{Decong-init}--\eqref{Decong-exact}. Pressure \eqref{pr2} for different values of $\gamma$.} 
Density (left) and velocity (right) at final time $T=0.2$: $\gamma=4$ (blue), $\gamma=20$ (red) and $\gamma=100$ (green). The simulations are performed with the implicit-explicit scheme and the initial condition is plotted in cyan.
\label{Decong3}}
\end{figure}

\begin{figure}
\subfloat[Density - Pressure \eqref{pr1bis} for different $\ve$ and $\gamma$]{
\includegraphics[trim=0 200 0 200,clip=true,
 scale=0.35]
{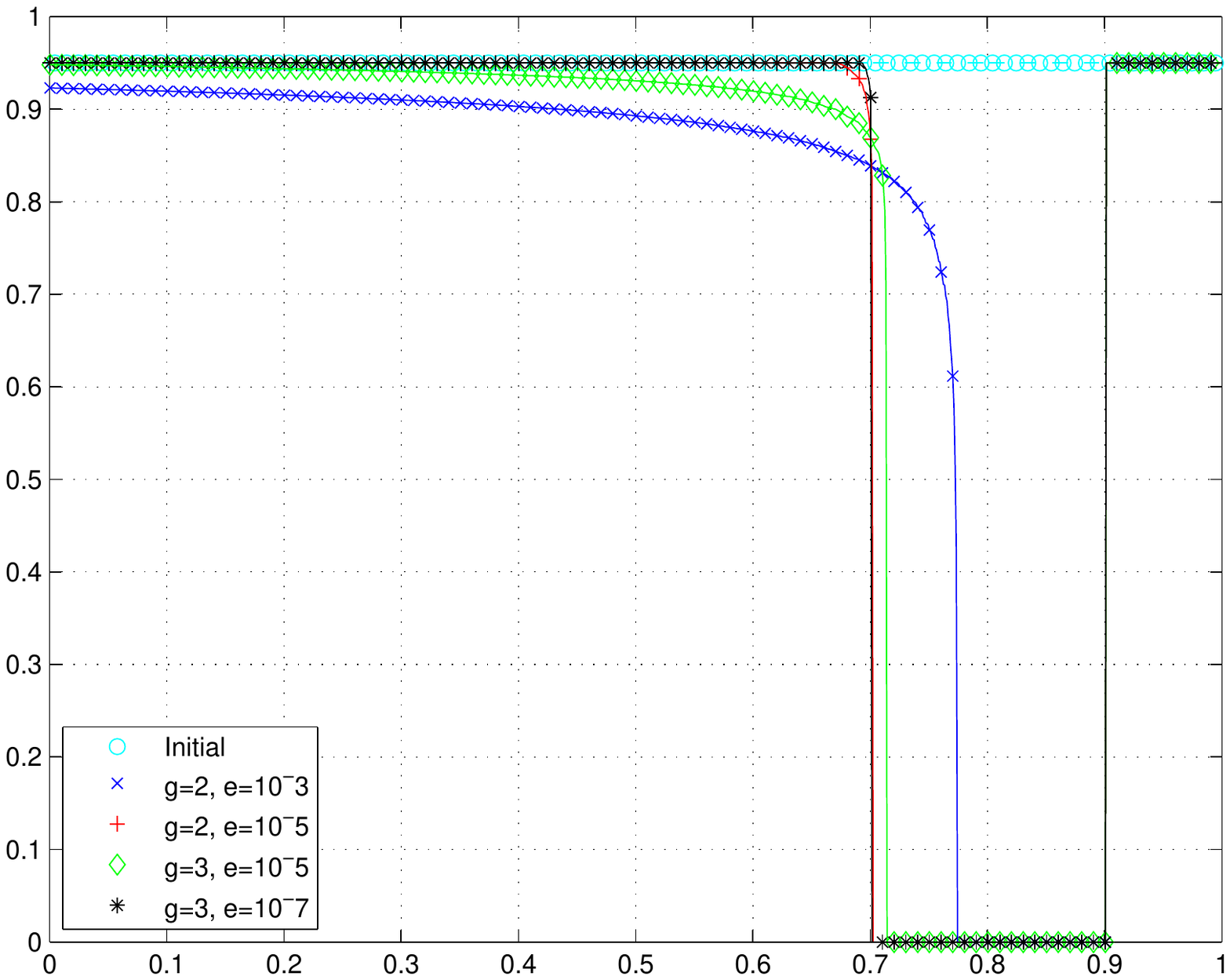}
}~~
\subfloat[Velocity - Pressure \eqref{pr1bis} for different $\ve$ and $\gamma$]{
\includegraphics[trim=0 200 0 200,clip=true,
scale=0.35]
{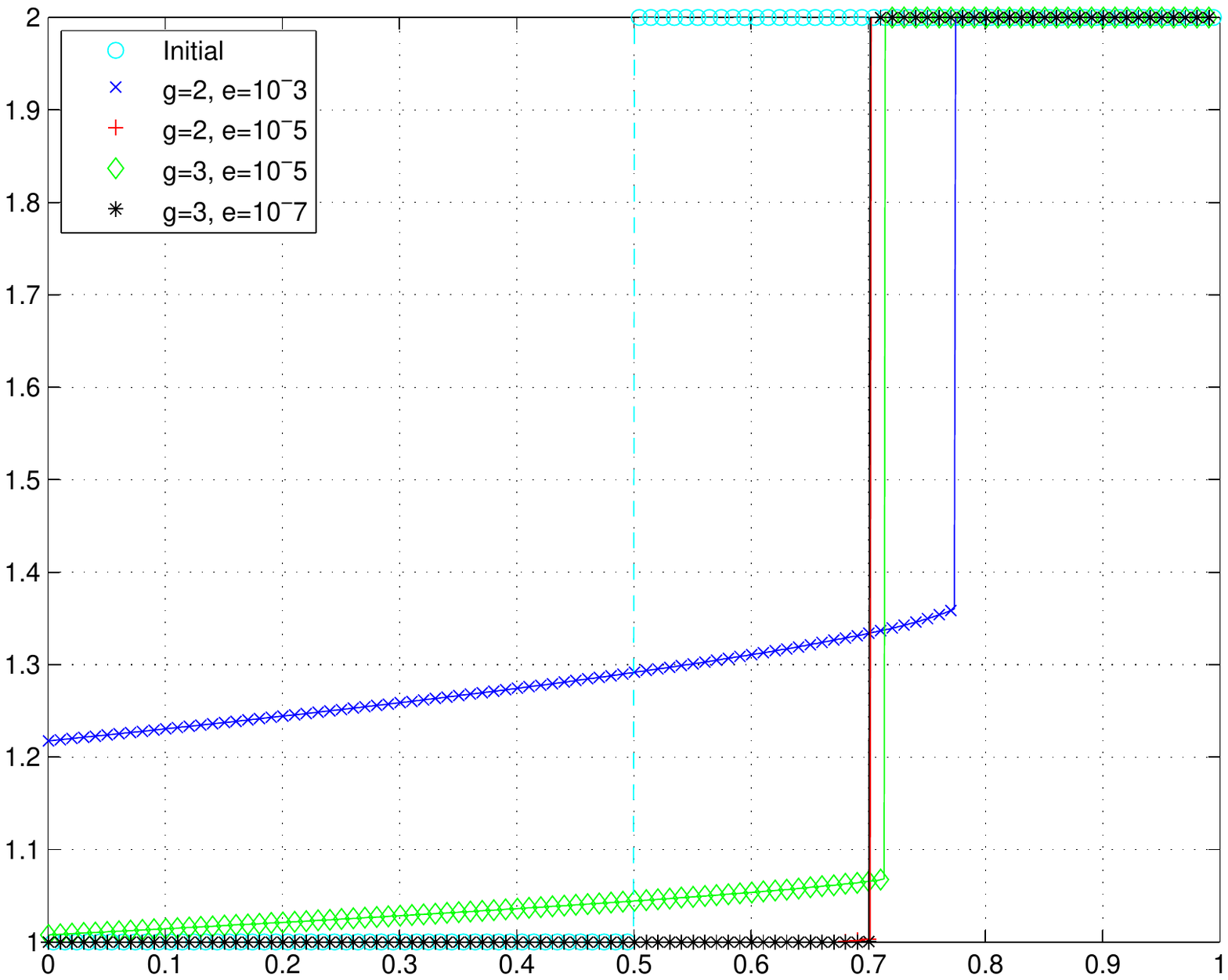}
}
\caption{\textbf{Numerical results in the case of decongestion \eqref{Decong-init}--\eqref{Decong-exact}. Pressure \eqref{pr1bis} for different values of $\ve$ and $\gamma$.} Density (left) and velocity (right) at final time $T=0.2$: \eqref{pr1bis} for $\gamma=2, \ve=10^{-3}$ (blue), $\gamma=2, \ve=10^{-5}$ (red), $\gamma=3, \ve=10^{-5}$ (green) and $\gamma=3, \ve=10^{-7}$ (black). The simulations are performed with the implicit-explicit scheme and the initial condition is plotted in cyan.
\label{Decong4}}
\end{figure}

Next, we study the case of a decongestion in the traffic. The data are defined by 
\begin{align}\label{Decong-init}
v^0(x)=
\left\{\begin{array}{ll}
1, & \text{ if } x \in [0, 0.5[, \\
2, & \text{ if } x \in [0.5, 1],
\end{array}\right.
\qquad \rho^0(x)=0.95.
\end{align}
The initial density is close to the threshold.
Since the vehicles ahead are going faster, a decongestion occurs.
The expected solution at $T=0.2$ should be
 equal to
\begin{align}\label{Decong-exact}
\rho(x,T)=
\left\{\begin{array}{ll}
0.95, & \text{ if } x \in [0, 0.7[, \\
0, & \text{ if } x \in [0.7, 0.9[, \\
0.95, & \text{ if } x \in [0.9, 1].
\end{array}\right.
\end{align}
Note the formation of a vacuum region, where  $v$ does not make sense.
It is indeed particularly interesting and relevant to check the ability of the models and of the numerical methods to handle the formation of vacuum regions, where the density vanishes.
The velocity offsets are defined as in the previous Section and we work with the same numerical parameters.
The results
can be found in Figure~\ref{Decong1}.
The Glimm scheme and the explicit-implicit scheme give the same results for all three pressures \eqref{pr1}, \eqref{pr1bis} 
and \eqref{pr2}.
Note that the density and velocity are the same for \eqref{pr1} and \eqref{pr1bis}, while, for the considered parameters,
\eqref{pr2} provides significantly different profiles.
We observe that   pressures \eqref{pr1} and \eqref{pr1bis} on the one hand and pressure  \eqref{pr2} on the other hand  give different results, especially in the vacuum region, none of them being totally in agreement 
with the ``expected'' result. 
 This can be explained by the moderate value of the parameters $\ve$ or $\gamma$.
Indeed,
the constrained behavior \eqref{Decong-exact} can be obtained by changing the parameters, see Figure~\ref{Decong3}, where 
we test \eqref{pr2} for different values of $\gamma \to +\infty$ and 
 Figure~\ref{Decong4}, where we make $\ve \to 0$ vary for \eqref{pr1bis}.
 In the former case, the limit behavior is captured with $\gamma=100$
 and for the latter, we get a satisfactory result with $ \ve=10^{-5}$ if $\gamma=2$, and $ \ve=10^{-7}$ if $\gamma=3$.
We notice also at Figure~\ref{Decong4} that if we take $\gamma=3$ for pressure \eqref{pr1bis} instead of $\gamma=2$, we need to take a smaller value of $\ve$, namely $\ve=10^{-7}$ instead of $\ve=10^{-5}$.
Note that the results for the original model \eqref{pr1} and for the modified model \eqref{pr1bis} are totally equivalent.

%
%

\subsection{Case of congestion}

\subsubsection{A jump of velocity creating a congestion}

Finally, we turn to the simulation of a congestion in the traffic. The initial conditions are given by 
\begin{align}\label{Cong-init}
\rho^0(x)=0.95,\qquad 
v^0(x)=
\left\{\begin{array}{ll}
2, \quad & \text{ if } x \in [0, 0.5[, \\
1, \quad & \text{ if } x \in [0.5, 1].
\end{array}\right.
\end{align}
The density is initially close to the threshold; since the cars ahead are slower, a congestion might occur
 and the Lagrange multiplier becomes active to prevent an excess of vehicles density.

Indeed, 
discontinuous solutions are characterized by the Rankine--Hugoniot conditions: with $t\mapsto s(t)$ 
the speed of the discontinuity curve, we have
\[ \dot s\ \big[\! \! \big[ \rho \big] \! \! \big]= \big[\! \! \big[ \rho v \big] \! \! \big],
\]
and
\[ \dot s\ \big[\! \! \big[ \rho (v+\pi) \big] \! \! \big]= \big[\! \! \big[ \rho v(v+\pi) \big] \! \! \big].\]
We can check that 
\[\rho_1(t,x)=
\left\{\begin{array}{ll}
0.95, \quad & \text{ if } x \in [0, 0.5-18t[, \\
1, \quad & \text{ if } x \in [0.5-18t, 0.5+t],
\\
0.95, \quad &\text{ if } x\in [0.5+t,1]
\end{array}\right.
\]
 with 
 \[v_1(t,x)=
\left\{\begin{array}{ll}
2, \quad & \text{ if } x \in [0, 0.5-18t[, \\
1, \quad & \text{ if } x \in [0.5-18t, 0.5+t],
\\
1, \quad &\text{ if } x\in [0.5+t,1]
\end{array}\right.
\]
and a Lagrange multiplier active only in the congestion domain
 \[\pi_1(t,x)=\mathbf 1_{0.5-18t \leq x\leq 0.5+t},
 \]
 is solution of \eqref{PGCD}.
%
The presence of slow vehicles 
ahead of the fast ones instantaneously creates a congestion behind the velocity jump:
the slow vehicles ahead
make the faster ones behind brake.
This is typical of
the Follow--the--Leader approach, which has led to a derivation of the Aw-Rascle-Zhang system \cite{AwR2, gazis}. 
However, it is likely that solutions of the constrained model \eqref{PGCD}
are not uniquely defined for such data; we refer the reader to \cite{BeBr} for such considerations. 

\begin{figure}
\subfloat[Density - Glimm scheme]{
\includegraphics[trim=0 200 0 200,clip=true,
scale=0.35]
{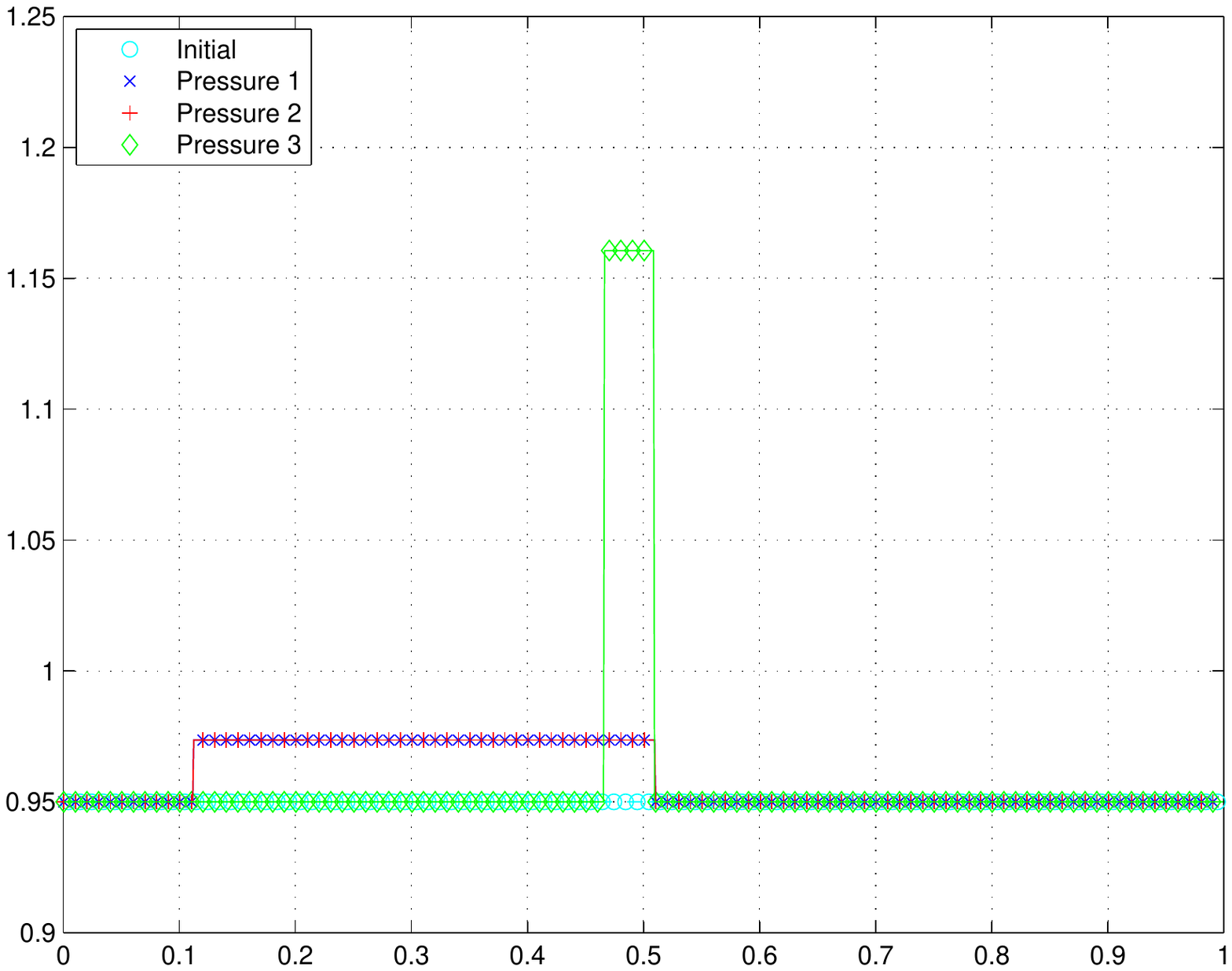}
}
\subfloat[Velocity - Glimm scheme]{
\includegraphics[trim=0 200 0 200,clip=true,
scale=0.35]
{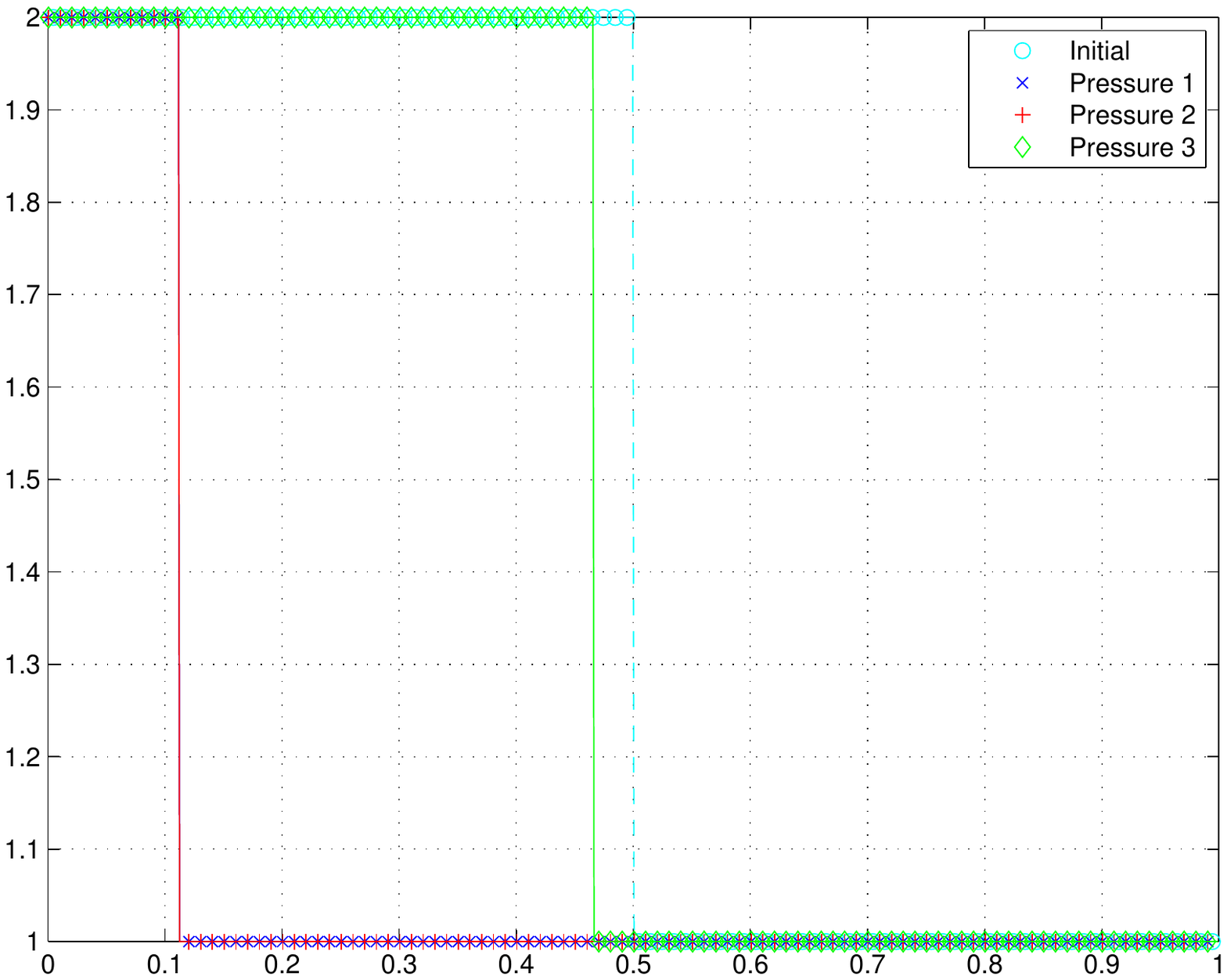}
}
\\
\subfloat[Density - implicit-explicit scheme]{
\includegraphics[trim=0 200 0 200,clip=true,
 scale=0.35]
{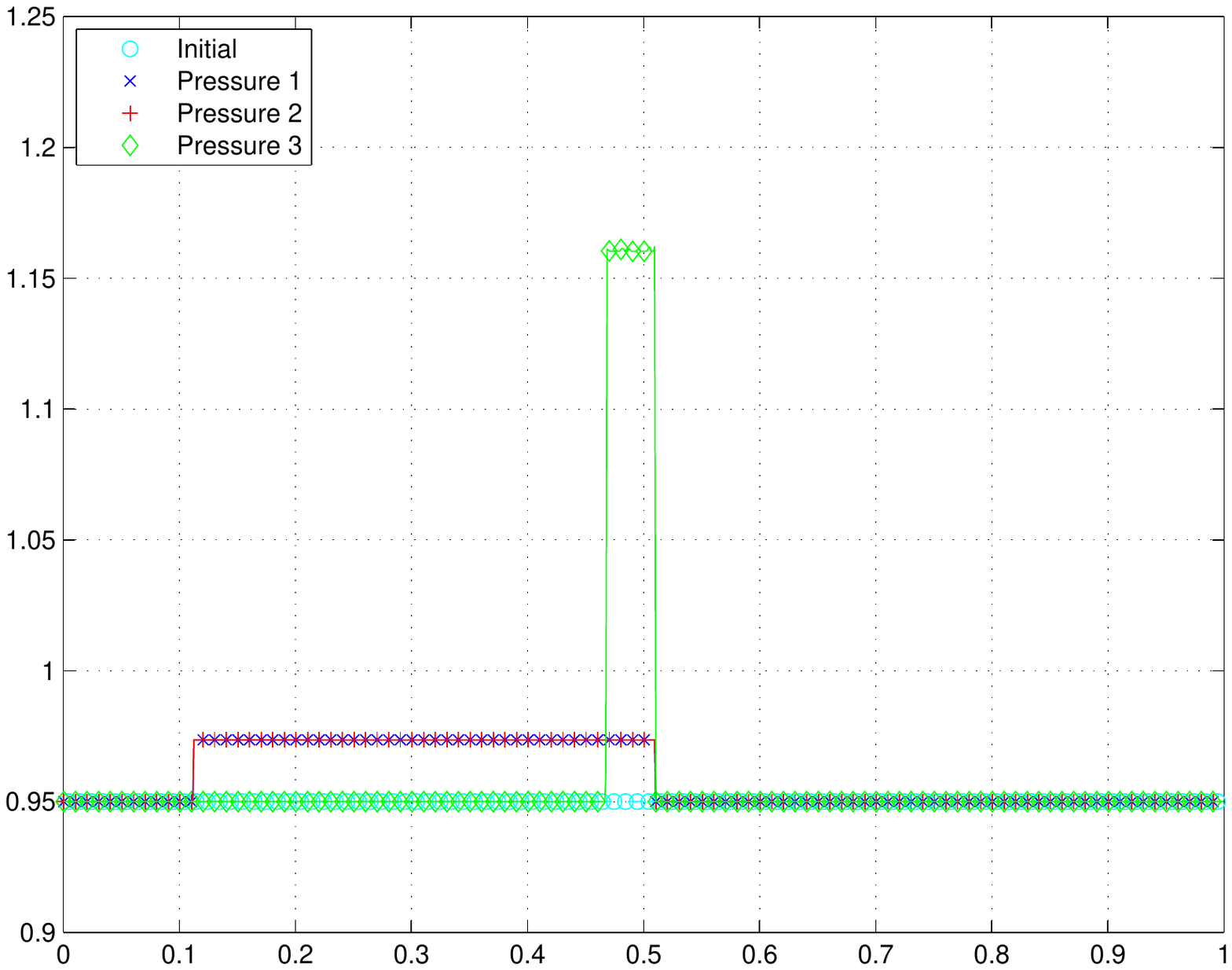}
}
\subfloat[Velocity - implicit-explicit scheme]{
\includegraphics[trim=0 200 0 200,clip=true,
scale=0.35]
{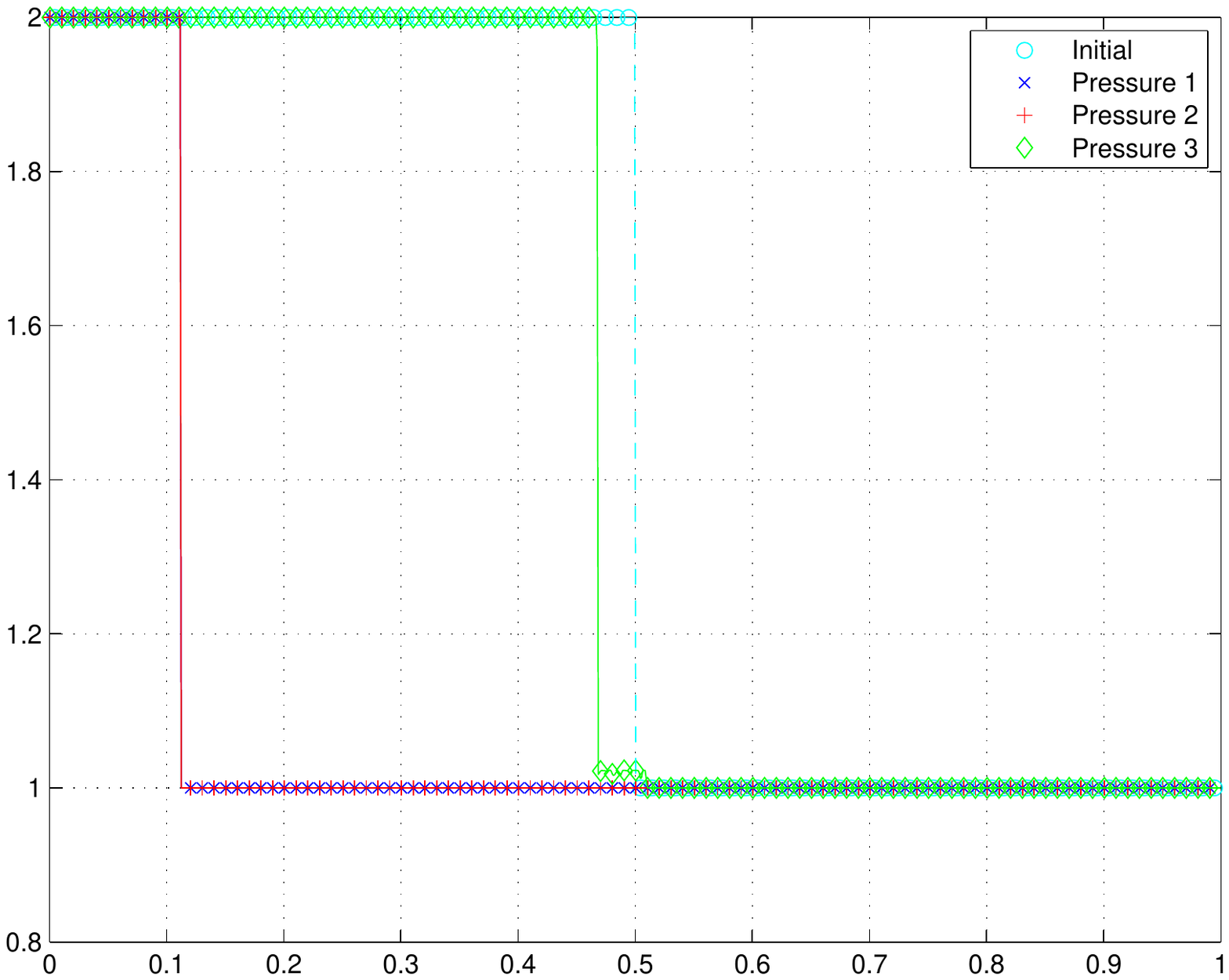}
}
\caption{\textbf{Numerical results in the case of congestion - Comparison of the two schemes.} Density (left) and velocity (right) at final time $T=0.01$. 
Glimm scheme (top)
and implicit-explicit scheme (bottom), with 
pressure \eqref{pr1} (blue), \eqref{pr1bis} (red) and \eqref{pr2} (green).
The initial conditions are plotted in cyan.
Parameters : $\gamma=2, \ve=10^{-3}$ for pressures  \eqref{pr1} and  \eqref{pr1bis}  and $\gamma=4$ for pressure \eqref{pr2}.
\label{Cong1} }
\end{figure}

The parameters are defined as in Section \ref{transport} and we show the solutions obtained at the final time $T=0.01$ in Figure~\ref{Cong1}.
We observe exactly the same behaviors between the Glimm scheme and the implicit-explicit scheme with \eqref{pr1} or \eqref{pr1bis}. 
We observe that with these parameters, the three models do not find the solution \eqref{Decong-exact}.
The time steps for the Glimm scheme are smaller than with the explicit-implicit scheme, but, quite surprisingly, by a factor 3 or 4 only.

\begin{figure}
\subfloat[Density - Glimm scheme]{
\includegraphics[trim=0 200 0 200,clip=true,
scale=0.35]
{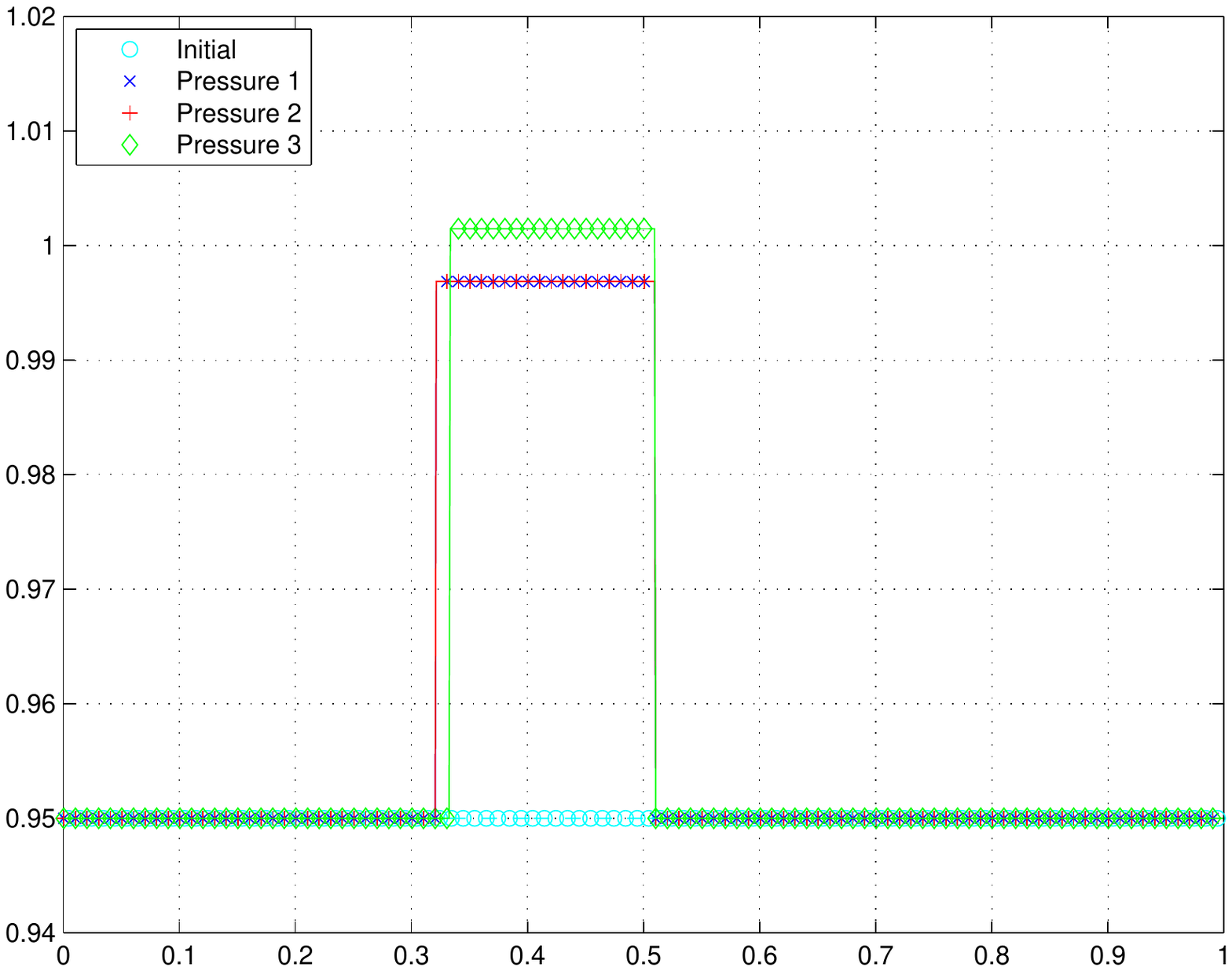}
}
\subfloat[Velocity - Glimm scheme]{
\includegraphics[trim=0 200 0 200,clip=true,
scale=0.35]
{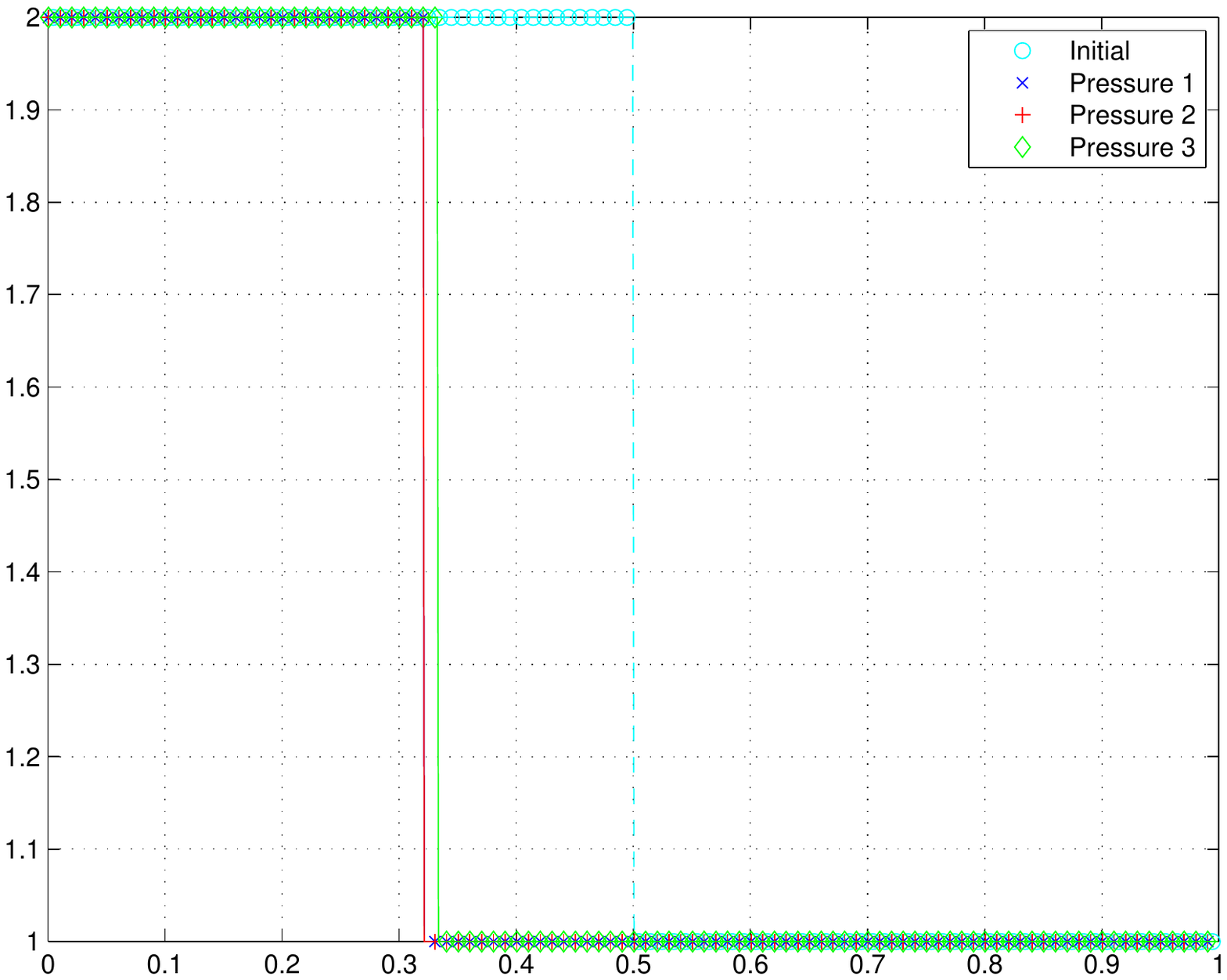}
}
\\
\subfloat[Density - implicit-explicit scheme]{
\includegraphics[trim=0 200 0 200,clip=true,
 scale=0.35]
{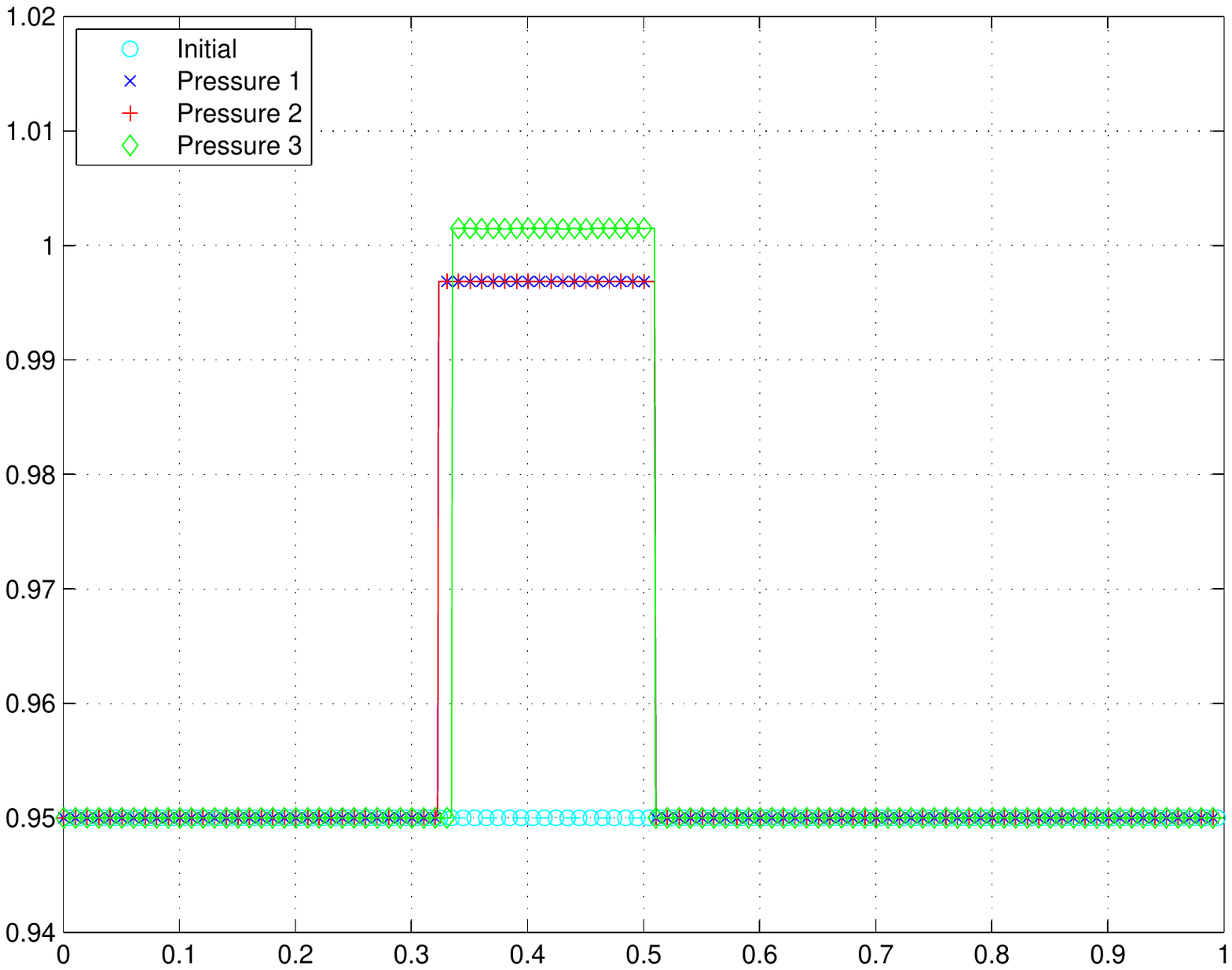}
}
\subfloat[Velocity - implicit-explicit scheme]{
\includegraphics[trim=0 200 0 200,clip=true,
scale=0.35]
{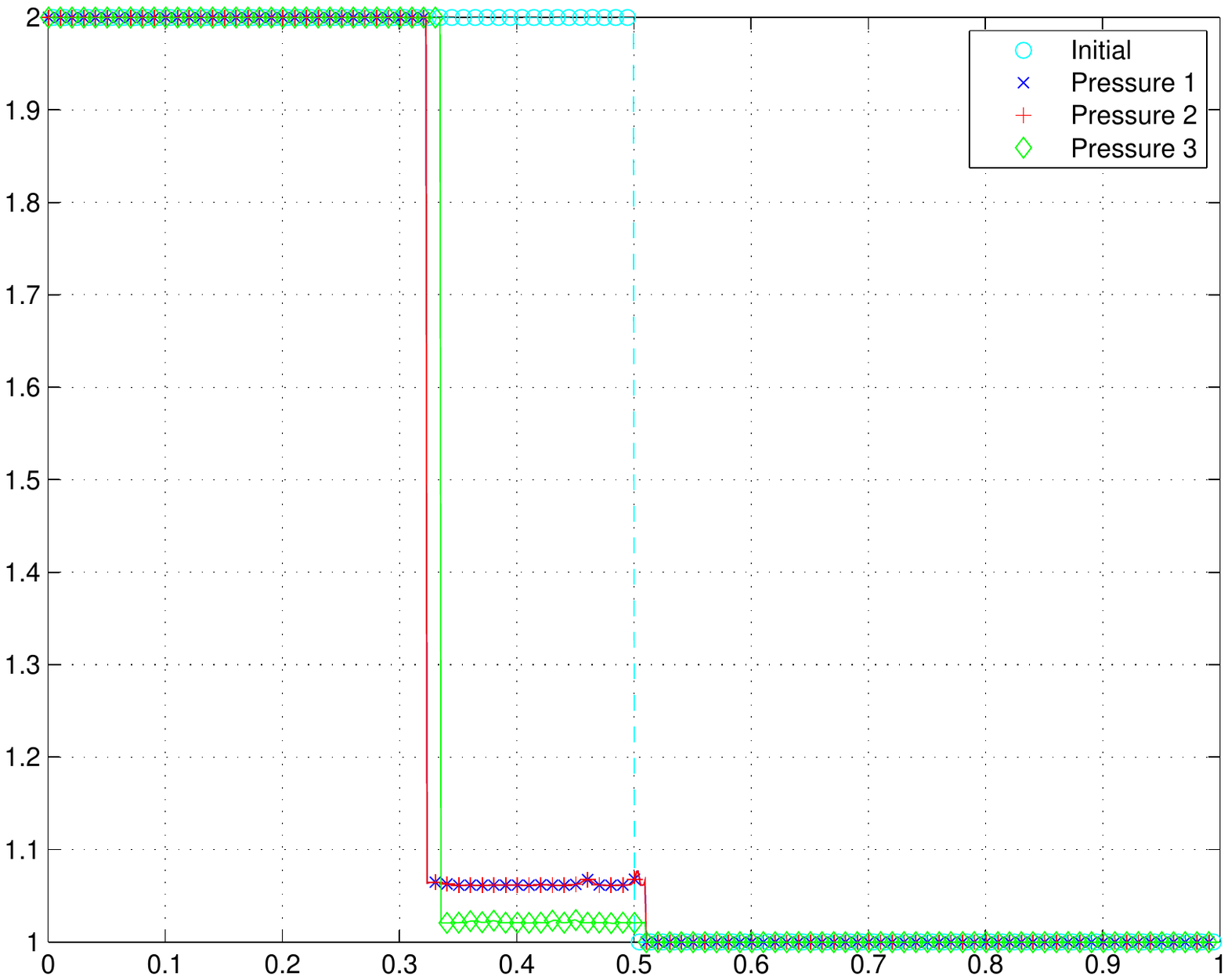}
}
\caption{\textbf{Numerical results in the case of congestion - Comparison of the two schemes - Different parameters.} For \eqref{pr1} and \eqref{pr1bis}, we use $\gamma=2$ and $\ve=10^{-5}$; for \eqref{pr2}, we take $\gamma=50$. 
Pressure \eqref{pr1} in blue,  pressure  \eqref{pr1bis} in red and pressure  \eqref{pr2} in green.
The initial conditions are plotted in cyan.
\label{Cong2}} 
\end{figure}

\begin{figure}
 \subfloat[Density - Glimm scheme]{
\includegraphics[trim=0 200 0 200,clip=true,
scale=0.35]
{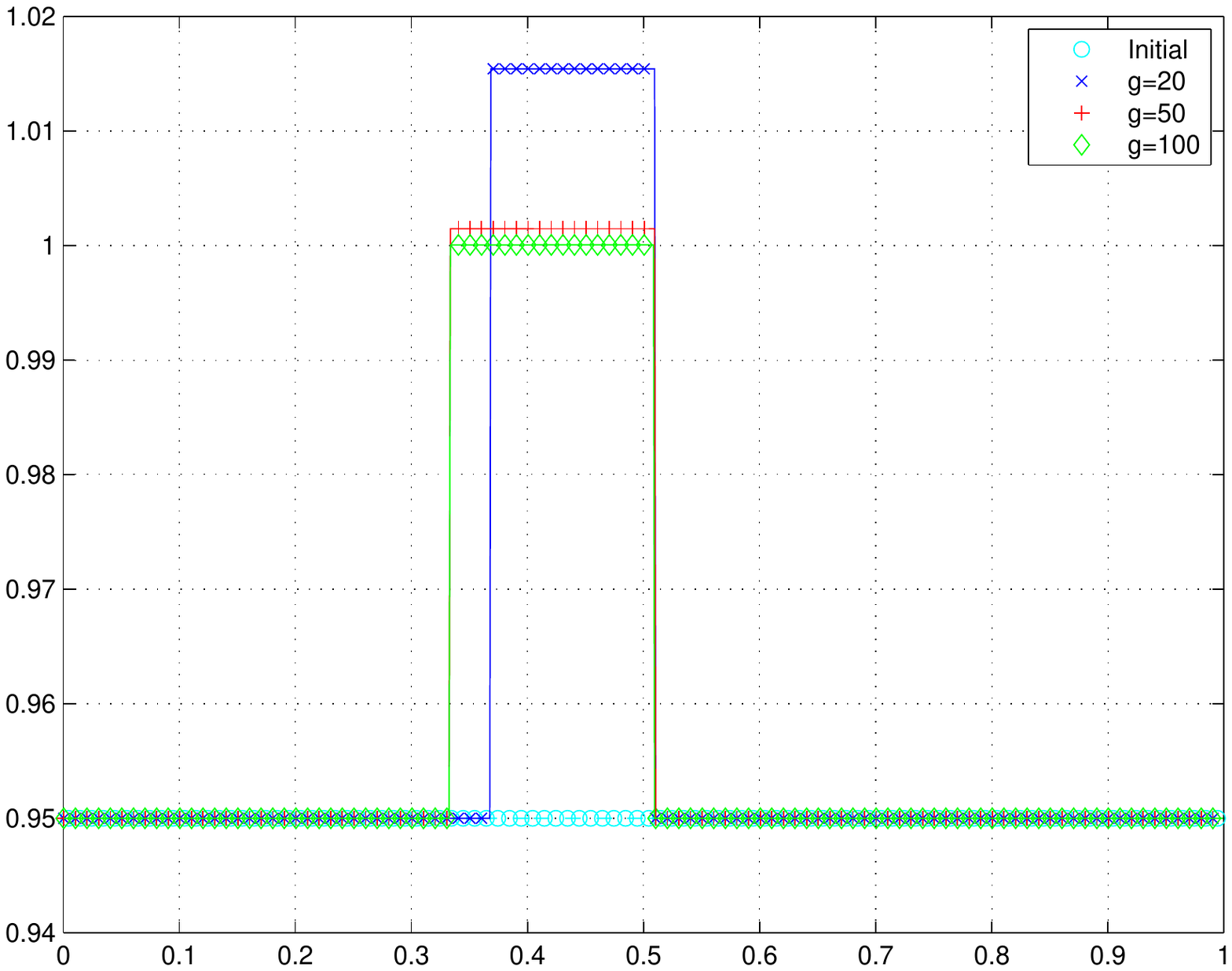}
}
\subfloat[Velocity - Glimm scheme]{
\includegraphics[trim=0 200 0 200,clip=true,
scale=0.35]
{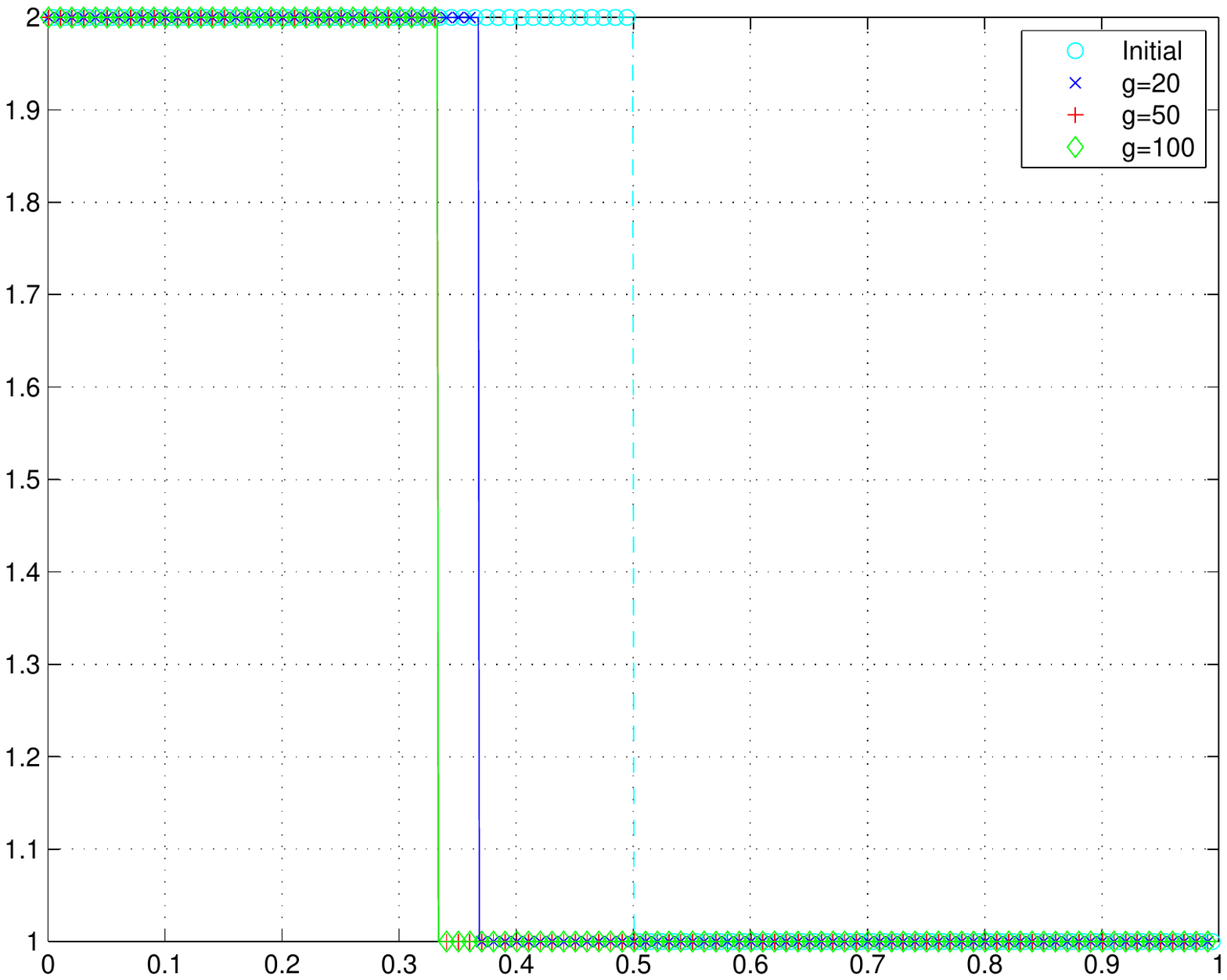}
}
\\
\subfloat[Density - implicit-explicit scheme]{
\includegraphics[trim=0 200 0 200,clip=true,
 scale=0.35]
{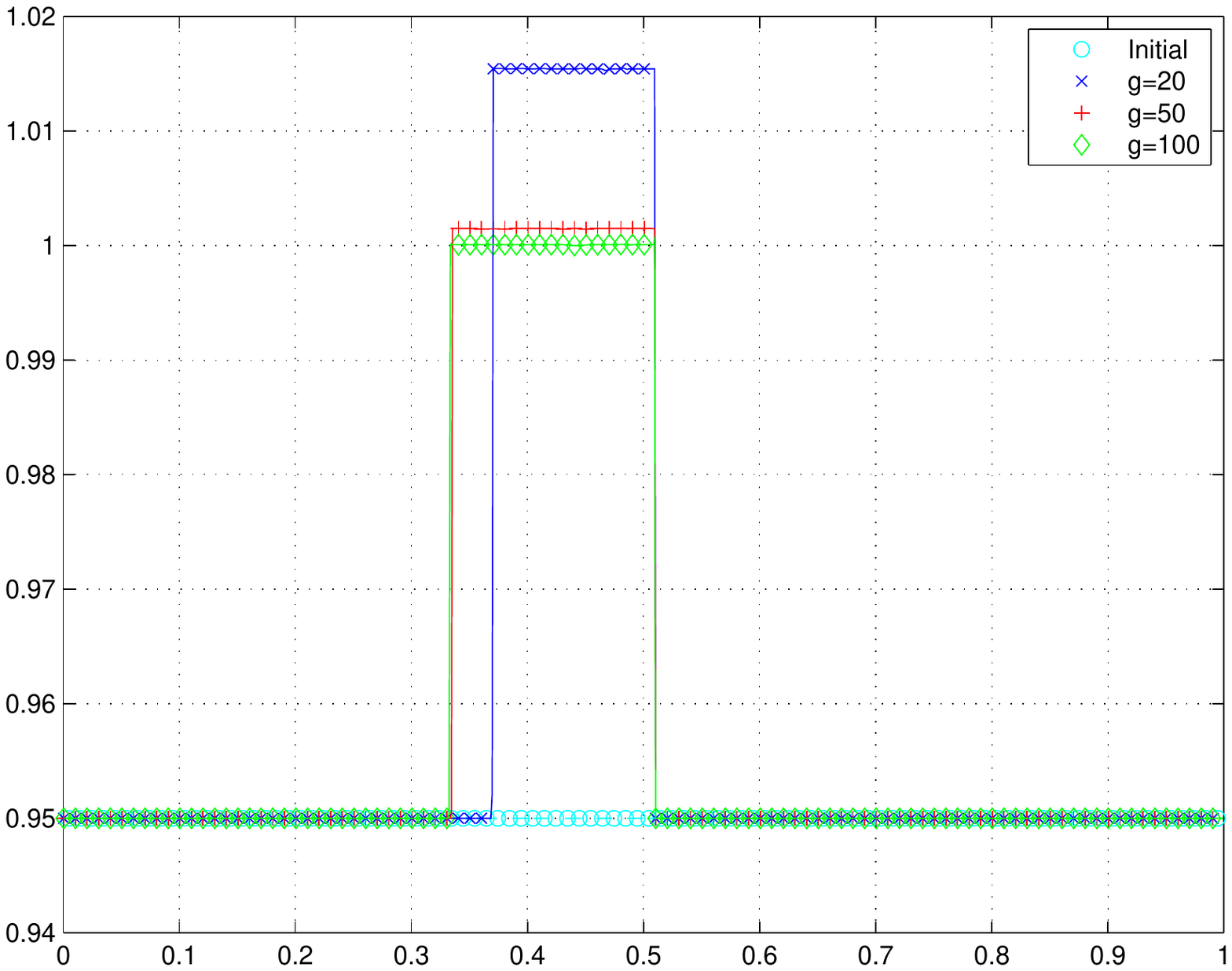}
}
\subfloat[Velocity - implicit-explicit scheme]{
\includegraphics[trim=0 200 0 200,clip=true,
scale=0.35]
{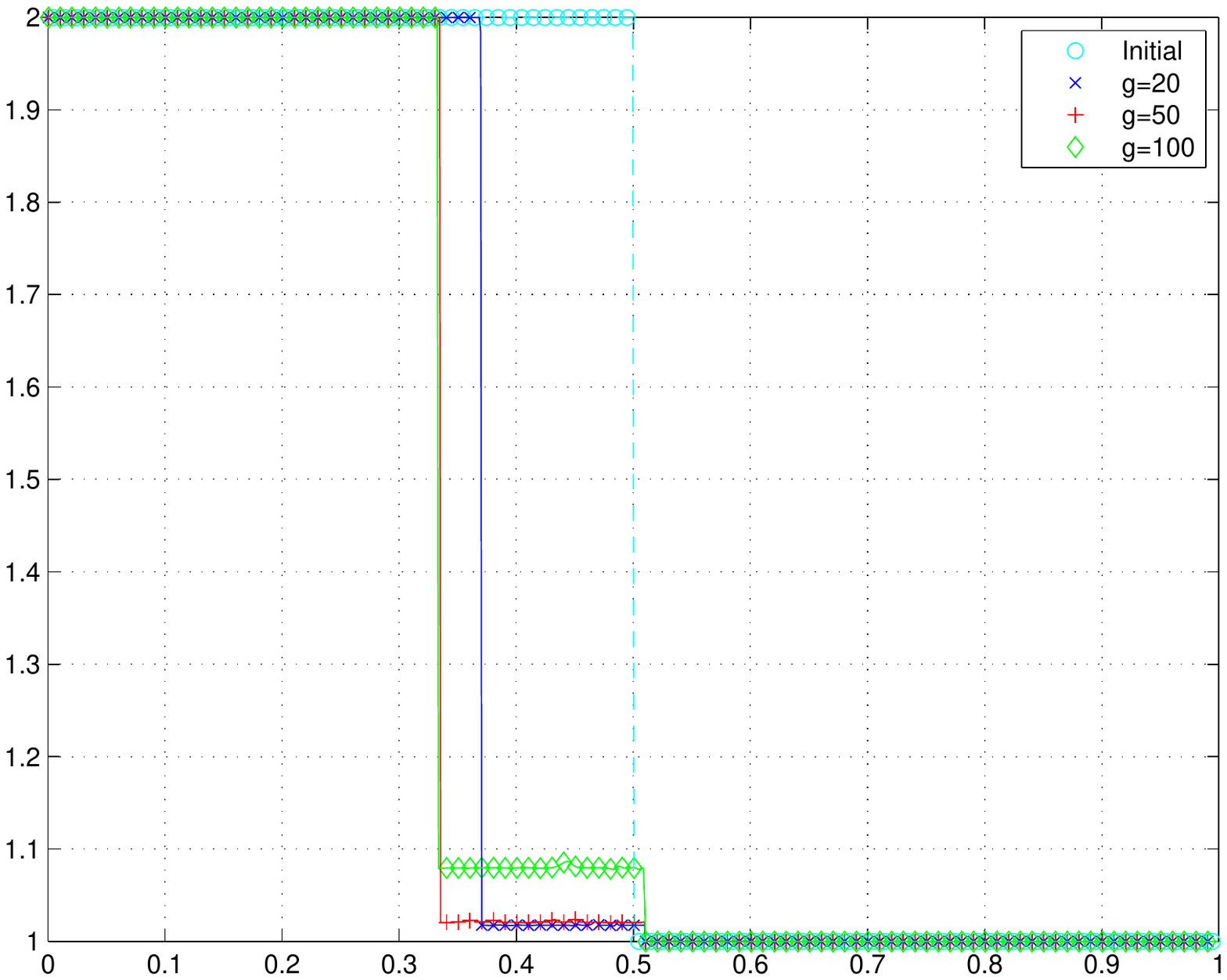}
}
\caption{\textbf{Numerical results in the case of congestion - Pressure \eqref{pr2} for different values of $\gamma$. } Density (left) and velocity (right) at final time $T=0.01$. 
Pressure \eqref{pr2} for $\gamma=20$ (red), $\gamma=50$ (blue) and $\gamma=100$ (green). The simulations are performed with the Glimm scheme (top) and the 
implicit-explicit scheme (bottom). The initial condition is plotted in cyan.
\label{Cong3}} 

\end{figure}

\begin{figure}
 \subfloat[Density - Glimm scheme]{
\includegraphics[trim=0 200 0 200,clip=true,
scale=0.35]
{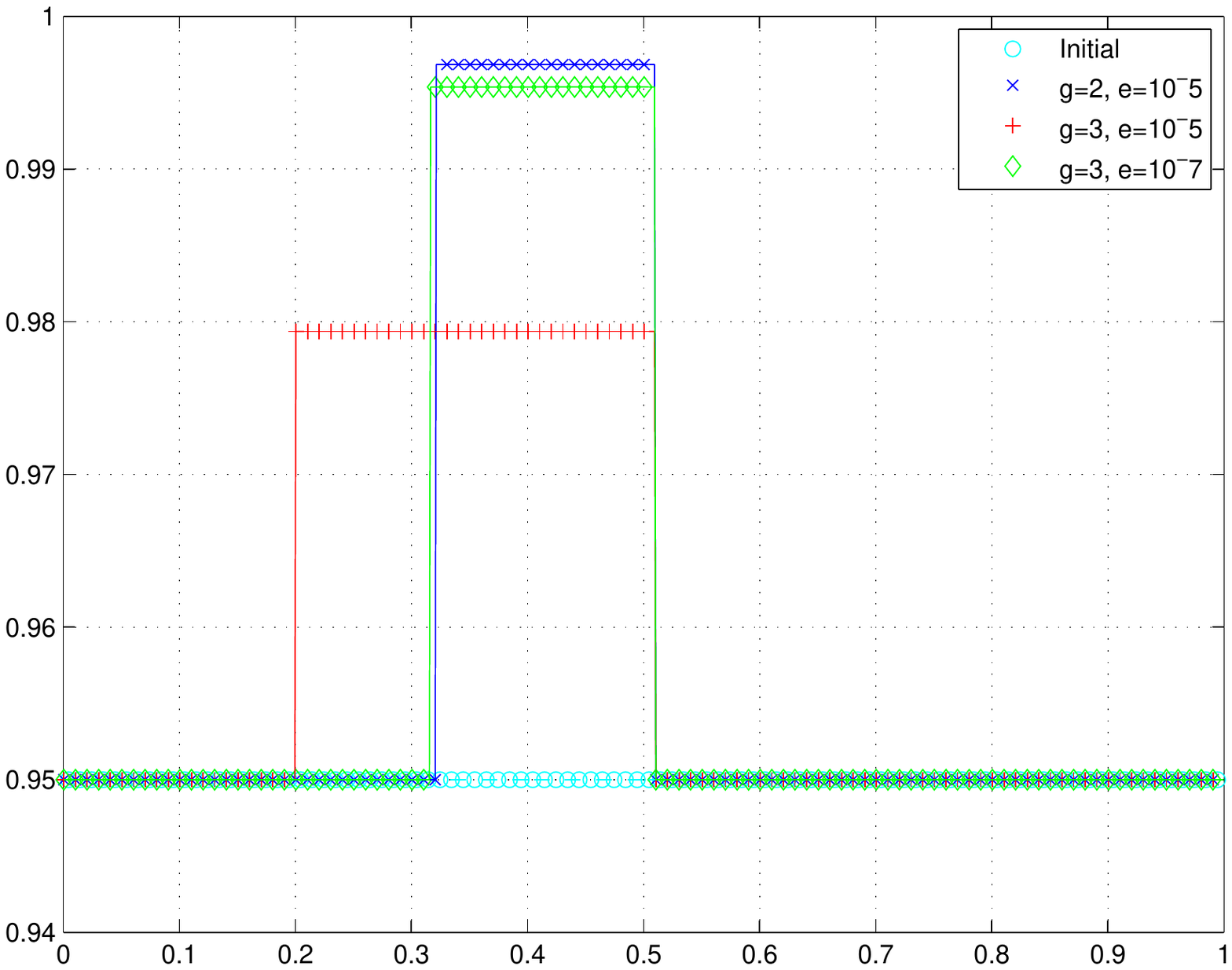}
}
\subfloat[Velocity - Glimm scheme]{
\includegraphics[trim=0 200 0 200,clip=true,
scale=0.35]
{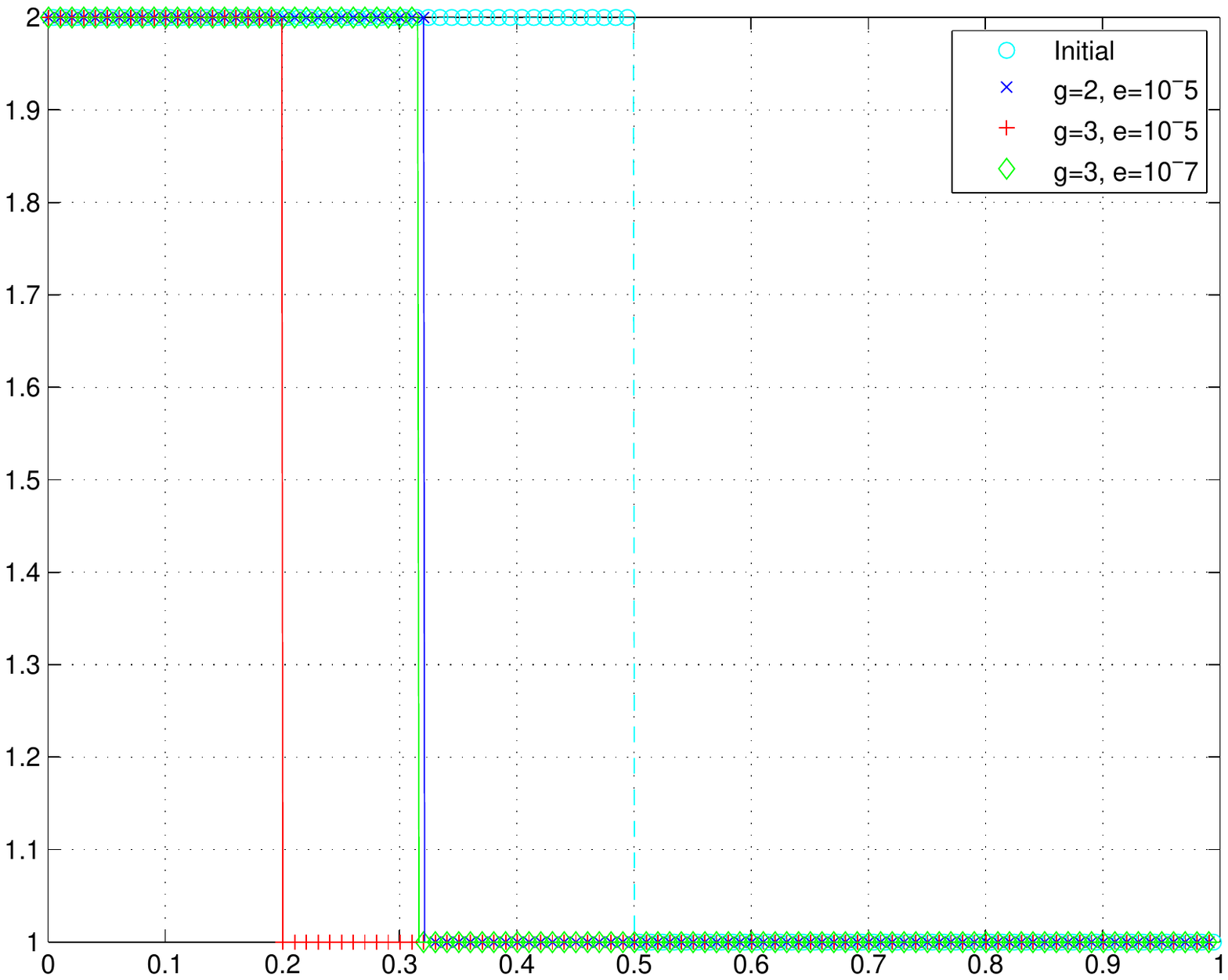}
}
\\

\subfloat[Density - implicit-explicit scheme]{
\includegraphics[trim=0 200 0 200,clip=true,
 scale=0.35]
{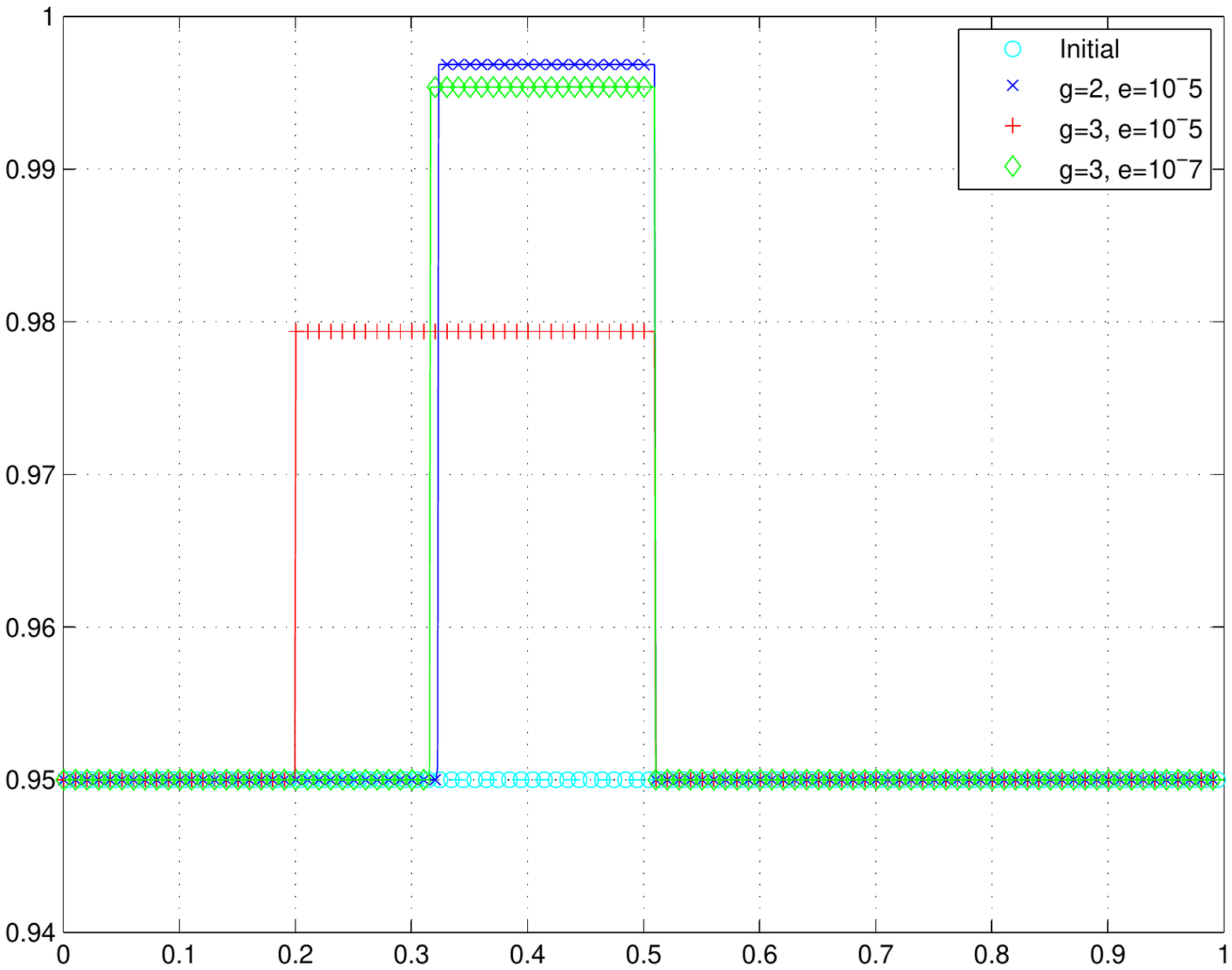}
}
\subfloat[Velocity - implicit-explicit scheme]{
\includegraphics[trim=0 200 0 200,clip=true,
scale=0.35]
{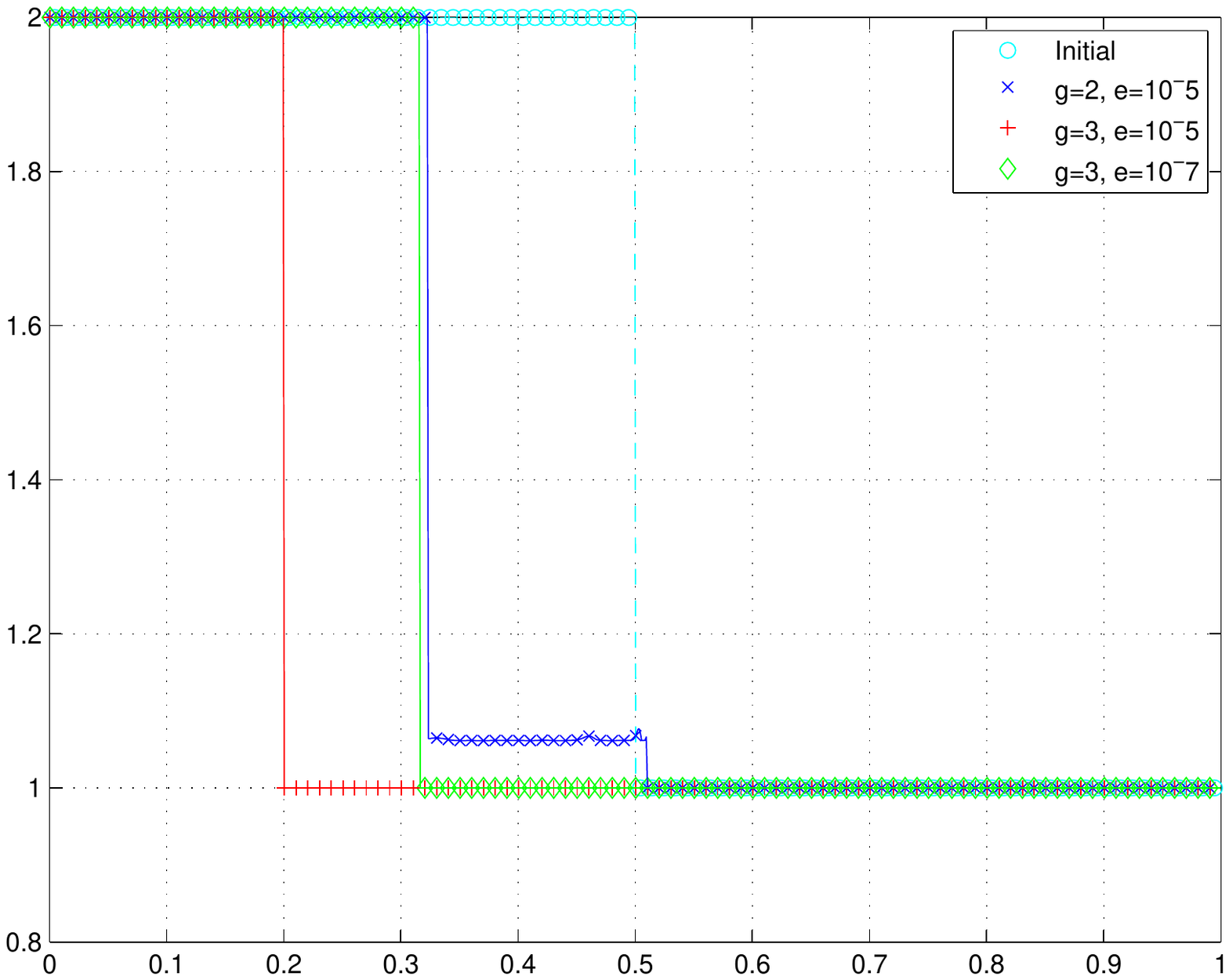}
}
\caption{\textbf{Numerical results in the case of congestion - Pressure \eqref{pr1bis} for different values of $\gamma$ and $\ve$. } Density (left) and velocity (right) at final time $T=0.01$. 
We compare the pressure \eqref{pr1bis} for $\gamma=2, \ve=10^{-5}$ (blue), $\gamma=3, \ve=10^{-5}$ (red) and $\gamma=3, \ve=10^{-7}$ (green). The simulations are performed with the Glimm scheme (top) and the 
implicit-explicit scheme (bottom). The initial condition is plotted in cyan.
\label{Cong4}} 
\end{figure}

Regarding the velocity offsets, \eqref{pr2} overshoots the maximal value of the density, equal to $1$, whereas the
 two other pressures \eqref{pr1} and \eqref{pr1bis} underestimate it. This is 
not surprising since \eqref{pr2} allows values larger than the threshold, but it contrasts with the behavior of the model \eqref{pr1bis} which has the same feature.
In Figure~\ref{Cong2}, we make the parameters vary as follows:
\begin{itemize}
\item Pressure \eqref{pr1} with $\gamma=2$, $\varepsilon=10^{-5}$, 
\item Pressure \eqref{pr1bis} with $\gamma=2$, $\varepsilon=10^{-5}$,
\item Pressure \eqref{pr2} with $\gamma=50$.
\end{itemize}
We observe that these parameters provide a result closer to the explicit solution $(\rho_1,v_1)$.
The results for \eqref{pr2} with different values of $\gamma$ are given at Figure~\ref{Cong3} and for \eqref{pr1bis} with different values of $\gamma$ and $\ve$ at Figure~\ref{Cong4}. The two schemes behave equivalently in that case, with some advantage in terms of time step for the implicit--explicit method, see Table~\ref{Tab}
(the smaller $\ve$, resp. the larger $\gamma$, the more important the gain). Note that in the case when $\gamma=3$, the results are exactly the same because the density is below the numerical threshold; consequently, the implicit part of the scheme is not used and the final result is the same as the one given by the explicit step, that is to say the final result is the same as the one given by the Glimm scheme.
As expected, making $\ve \to 0$ for \eqref{pr1bis} or $\gamma \to +\infty$ for \eqref{pr2} allows us to obtain a result compatible with $(\rho_1,v_1)$.
In particular, we point out that the approach 
of \eqref{PGCD} as an asymptotic model 
 from \eqref{RMAR} provides the solution 
 where the fastest cars should brake behind the slow vehicles, independently of the density of slow vehicles ahead.
 This behavior corresponds to the derivation originally introduced in \cite{AwR1}.

\begin{table}[htbp!]
\centering
\label{table} 
\begin{tabular}{|c|c|c|c|}
\hline Pressure & Time step & Time step & Factor \\
 \& param & Glimm scheme & explicit-implicit scheme & \\
 \hline
Pressure \eqref{pr1bis}, $\ve=10^{-4}$ & $\Delta t=2 \cdot 10^{-6}$ & $\Delta t=2 \cdot 10^{-6}$ & 1 \\
\hline
Pressure \eqref{pr1bis}, $\ve=10^{-5}$ & $\Delta t=7 \cdot 10^{-7}$ & $\Delta t=10^{-6}$ & 1.39 \\
\hline
Pressure \eqref{pr1bis}, $\ve=10^{-6}$ & $\Delta t=2 \cdot 10^{-7}$ & $\Delta t=7.7 \cdot 10^{-7}$ & 3.22\\
\hline
Pressure \eqref{pr1bis}, $\ve=10^{-7}$ &$\Delta t=7.5 \cdot 10^{-8}$ & $\Delta t=6.2 \cdot 10^{-7}$ & 8.18 \\
\hline
\hline

Pressure \eqref{pr2}, $\gamma=50$ & $\Delta t=9 \cdot 10^{-6}$ & $\Delta t=10^{-5}$ & 1.12 \\
\hline
Pressure \eqref{pr2}, $\gamma=100$ & $\Delta t=4.8 \cdot 10^{-6}$ & $\Delta t=6.4 \cdot 10^{-6}$ & 1.36 \\
\hline
Pressure \eqref{pr2}, $\gamma=200$ & $\Delta t=2.4 \cdot 10^{-6}$ & $\Delta t=5.6 \cdot 10^{-6}$ & 2.33 \\
\hline
Pressure \eqref{pr2}, $\gamma=500$ & $\Delta t=9.5 \cdot 10^{-7}$ &$\Delta t=2.7 \cdot 10^{-5}$ & 27.95 \\

\hline
\end{tabular}
\caption{\textbf{Time steps - comparison between Glimm scheme and the explicit-implicit scheme for the congestion case }. Pressure \eqref{pr2} for different values of $\gamma$ and pressure \eqref{pr1bis} for $\gamma=2$ and different values of $\varepsilon$. The time step is the smallest time step used during the simulation and the factor is the ratio of the time step for the explicit-implicit scheme over the time step for the Glimm scheme. \label{Tab}}
\end{table}

\subsubsection{A fast cluster reaching a slower one}

\begin{figure}
\subfloat[Density - Glimm scheme \label{2BlocPO:D}]
	{
\includegraphics[trim=0 200 0 200,clip=true,
 scale=0.35]
{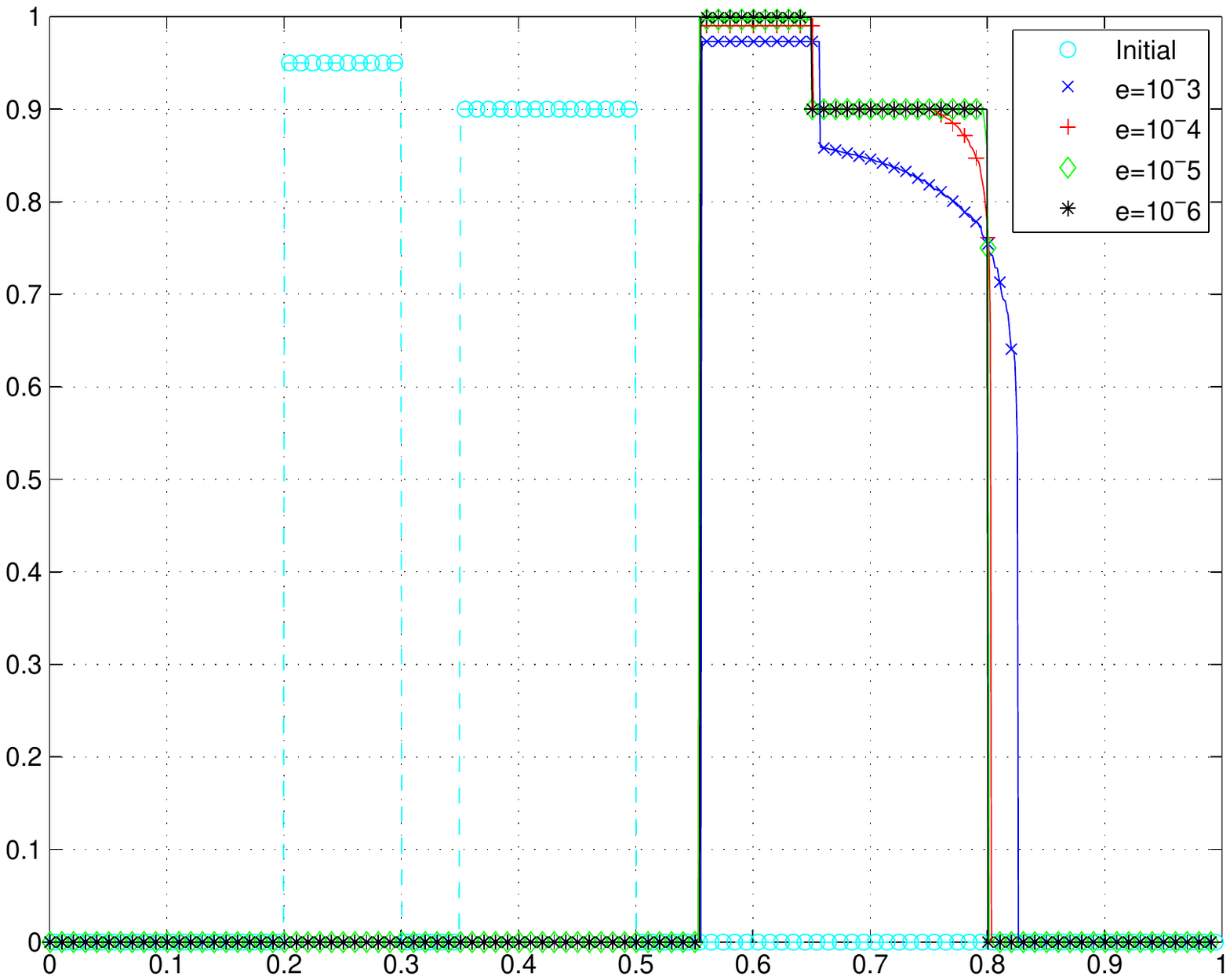}
	}
\quad
\subfloat[Velocity - Glimm scheme \label{2BlocPO:V}]
	{
\includegraphics[trim=0 200 0 200,clip=true,
 scale=0.35]
{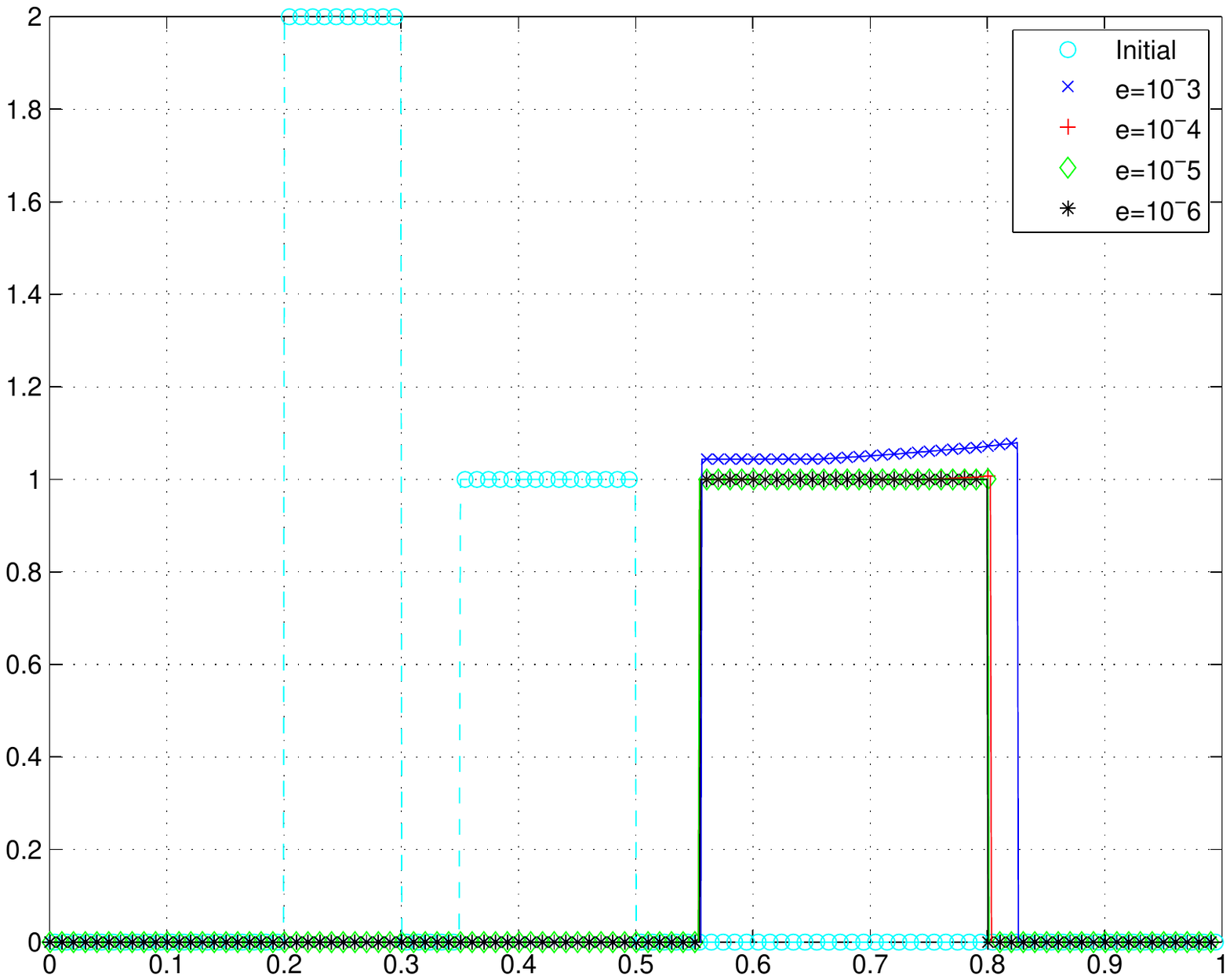}
	}
\caption{
\label{2BlocPO}
\textbf{Numerical results in the case of a shock between two blocks - Pressure \eqref{pr1} for $\gamma=2$ and different values of $\ve$}:
$\ve=10^{-3}$ (blue), $\ve=10^{-4}$ (red), $\ve=10^{-5}$ (green) and $\ve=10^{-6}$ (black). 
Figure~\ref{2BlocPO:D} (on the left) represents the densities whereas figure~\ref{2BlocPO:V}  (on the right) represents the velocities.
The initial condition is plotted in cyan.
}
\end{figure}

\begin{figure}
\subfloat[Density - Glimm scheme \label{2BlocNP:DT}]
{
\includegraphics[trim=0 200 0 200,clip=true,
 scale=0.35]
{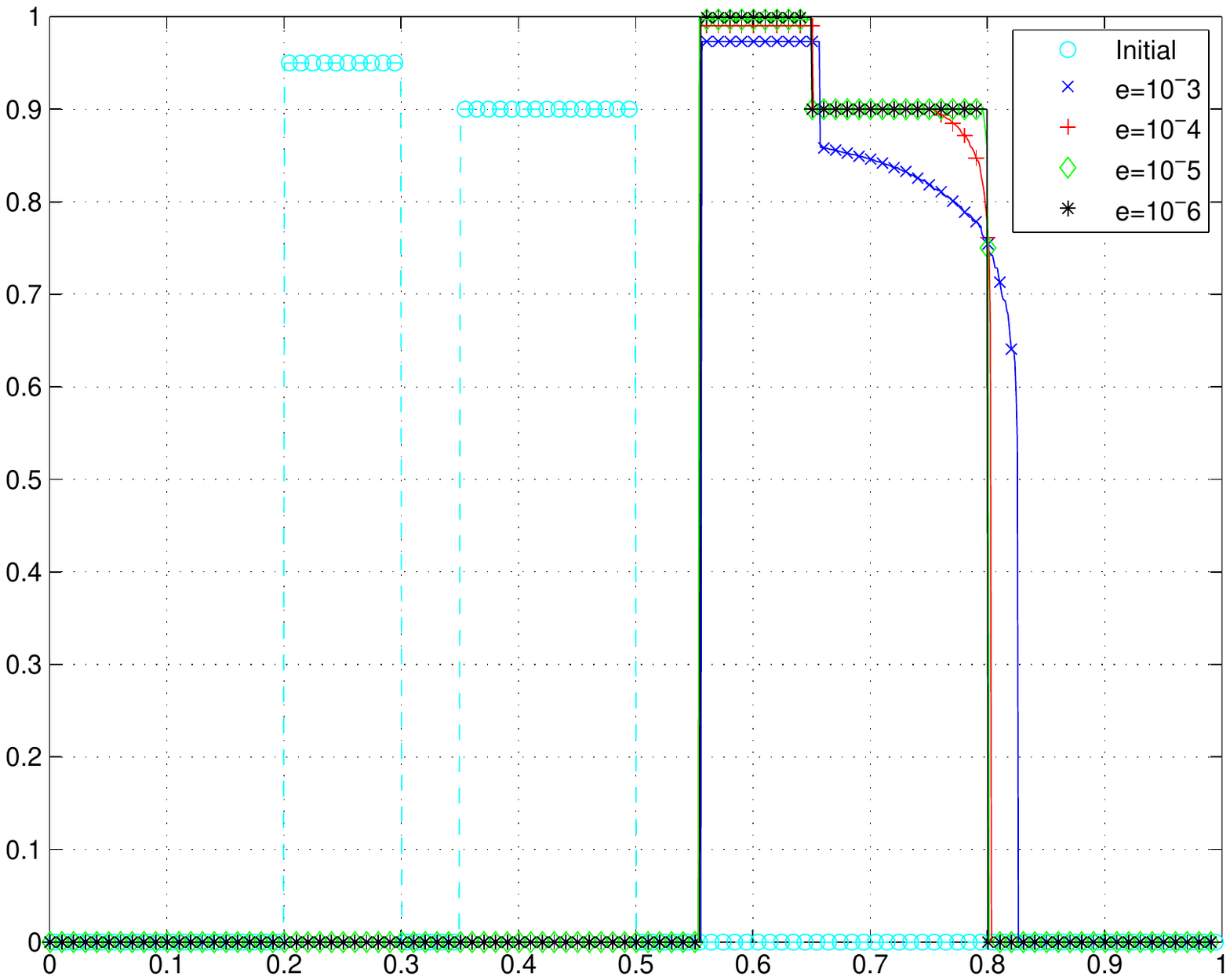}
	}
\quad
\subfloat[Velocity - Glimm scheme \label{2BlocNP:VT}]
	{
	\includegraphics[trim=0 200 0 200,clip=true,
 scale=0.35]
{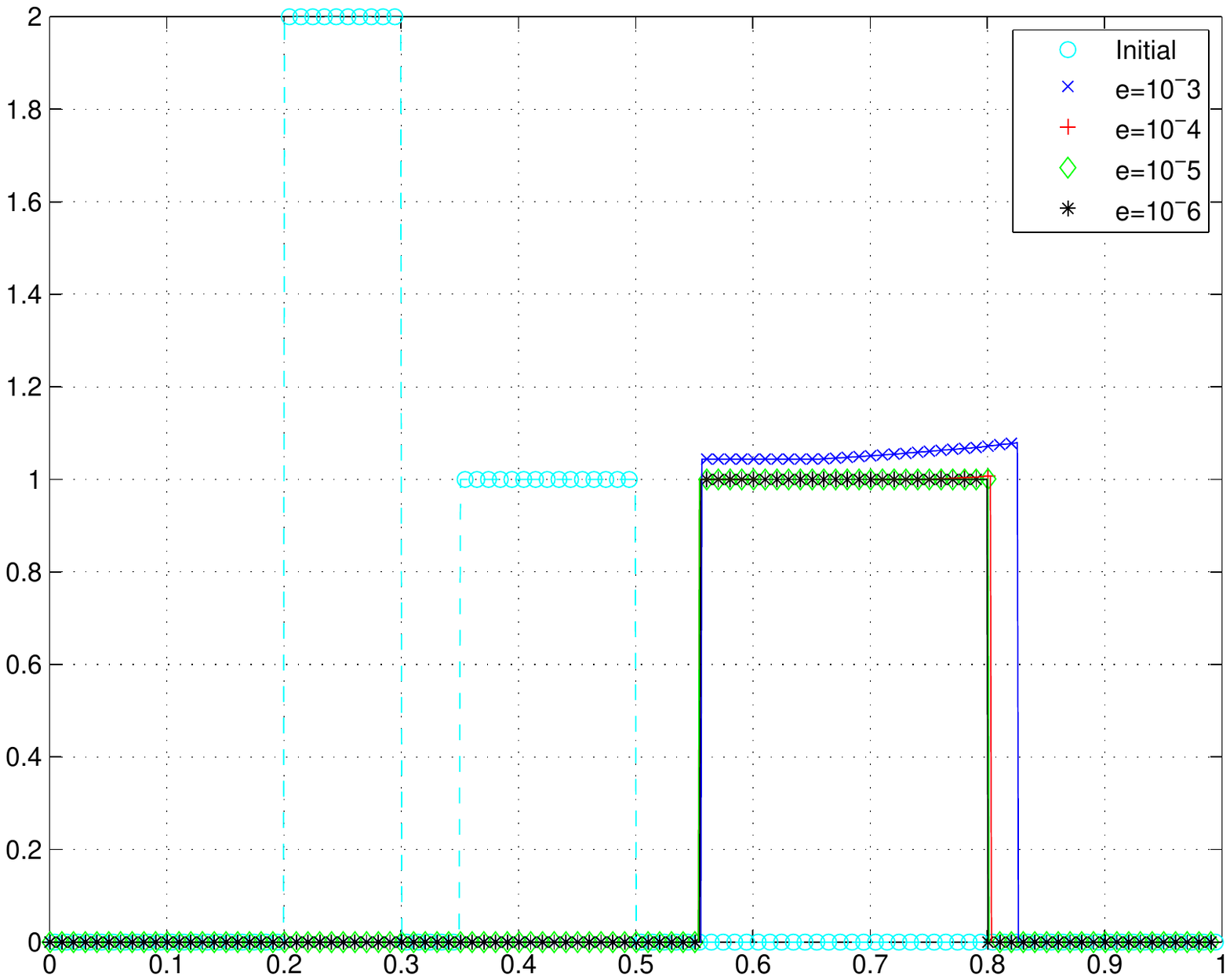}
	}
\\ 
\subfloat[Density - Implicit-explicit scheme \label{2BlocNP:DG}]
	{
\includegraphics[trim=0 200 0 200,clip=true,
 scale=0.35]
{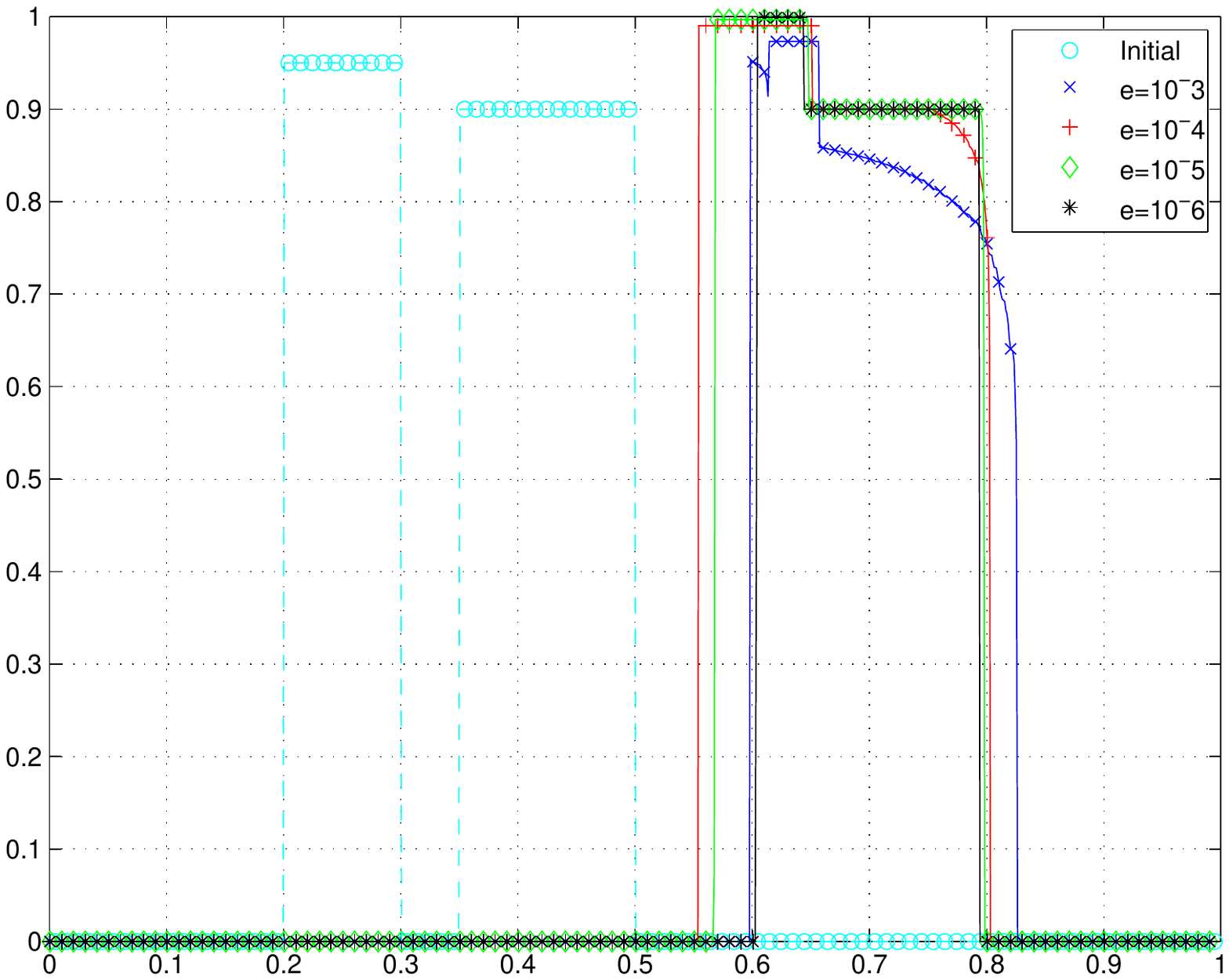}
	}
\quad
\subfloat[Velocity - Implicit-explicit scheme \label{2BlocNP:VG}]
	{
\includegraphics[trim=0 200 0 200,clip=true,
 scale=0.35]
{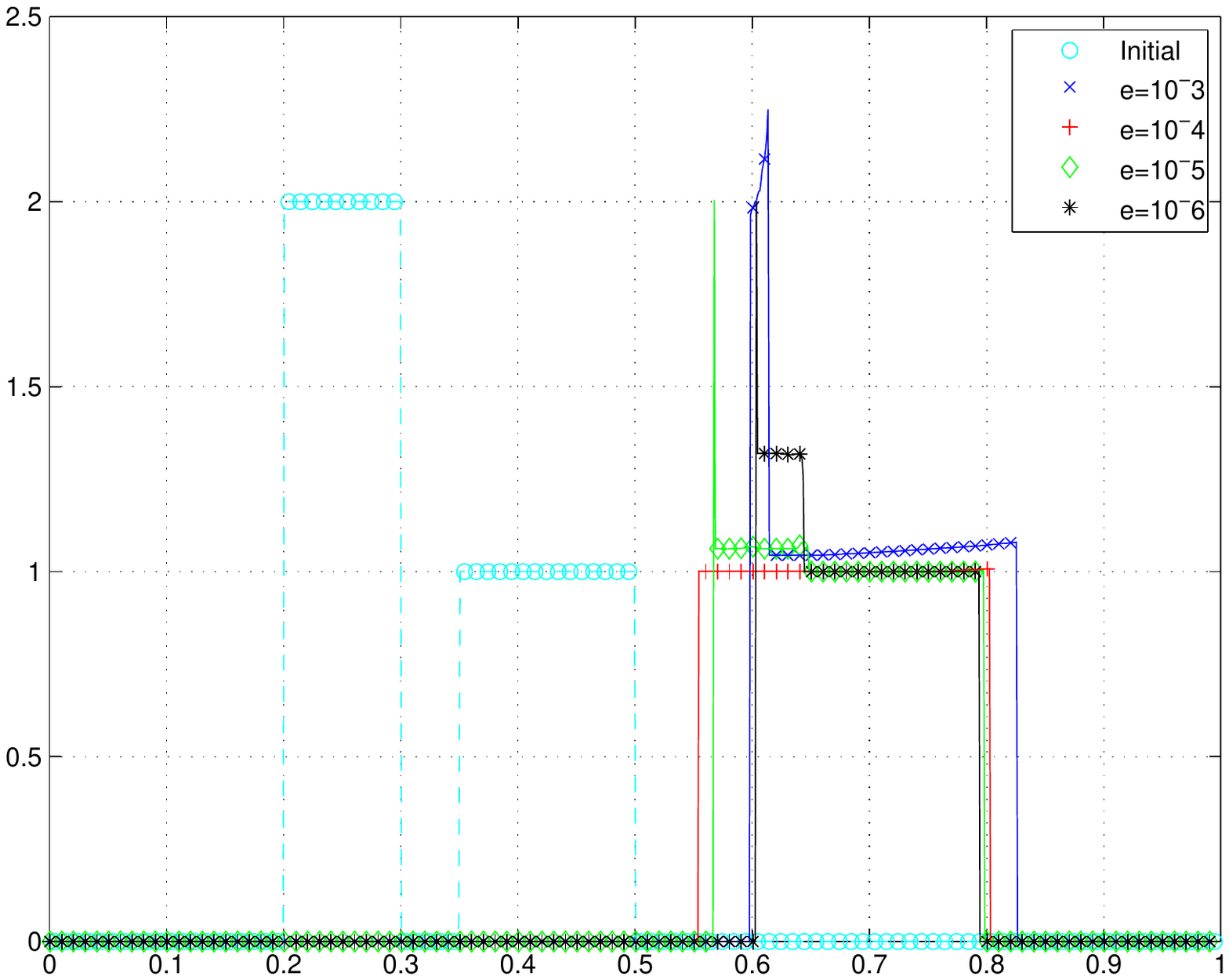}
	}
\caption{
\label{2BlocNP}
\textbf{Numerical results in the case of a shock between two blocks - Pressure \eqref{pr1bis} for $\gamma=2$ and different values of $\ve$}:
$\ve=10^{-3}$ (blue), $\ve=10^{-4}$ (red), $\ve=10^{-5}$ (green) and $\ve=10^{-6}$ (black). 
On top, simulations are performed with the Glimm scheme and on the bottom, with the implicit-explicit scheme. We display the densities on the left and the velocities on the right.
The initial condition is plotted in cyan.
}
\end{figure}

\begin{figure}
\subfloat[Density - Glimm scheme\label{2BlocPP:DT}]
	{
	\includegraphics[trim=0 200 0 200,clip=true,
 scale=0.35]
{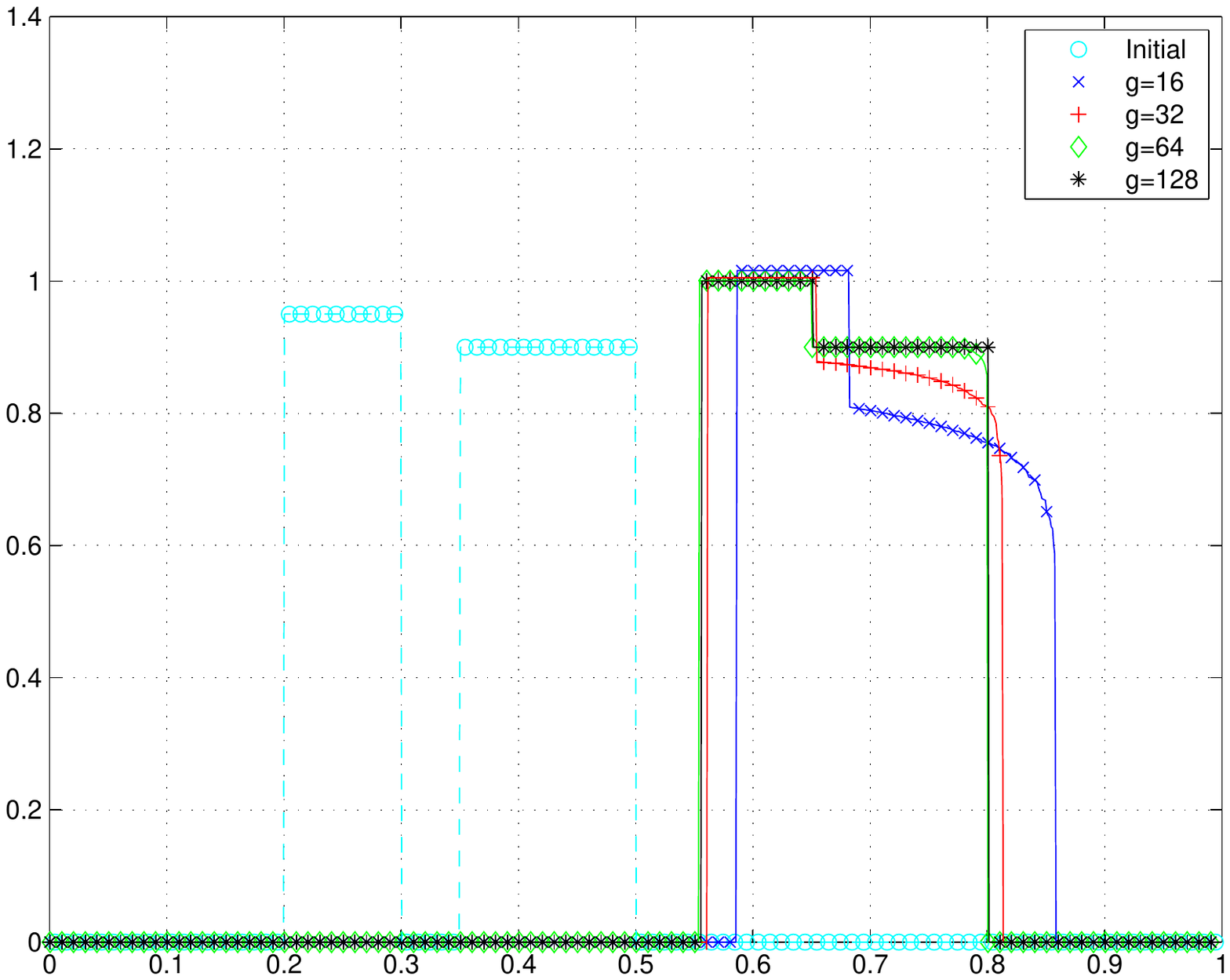}
	}
\quad
\subfloat[Velocity -  Glimm scheme\label{2BlocPP:VT}]
	{
	\includegraphics[trim=0 200 0 200,clip=true,
 scale=0.35]
{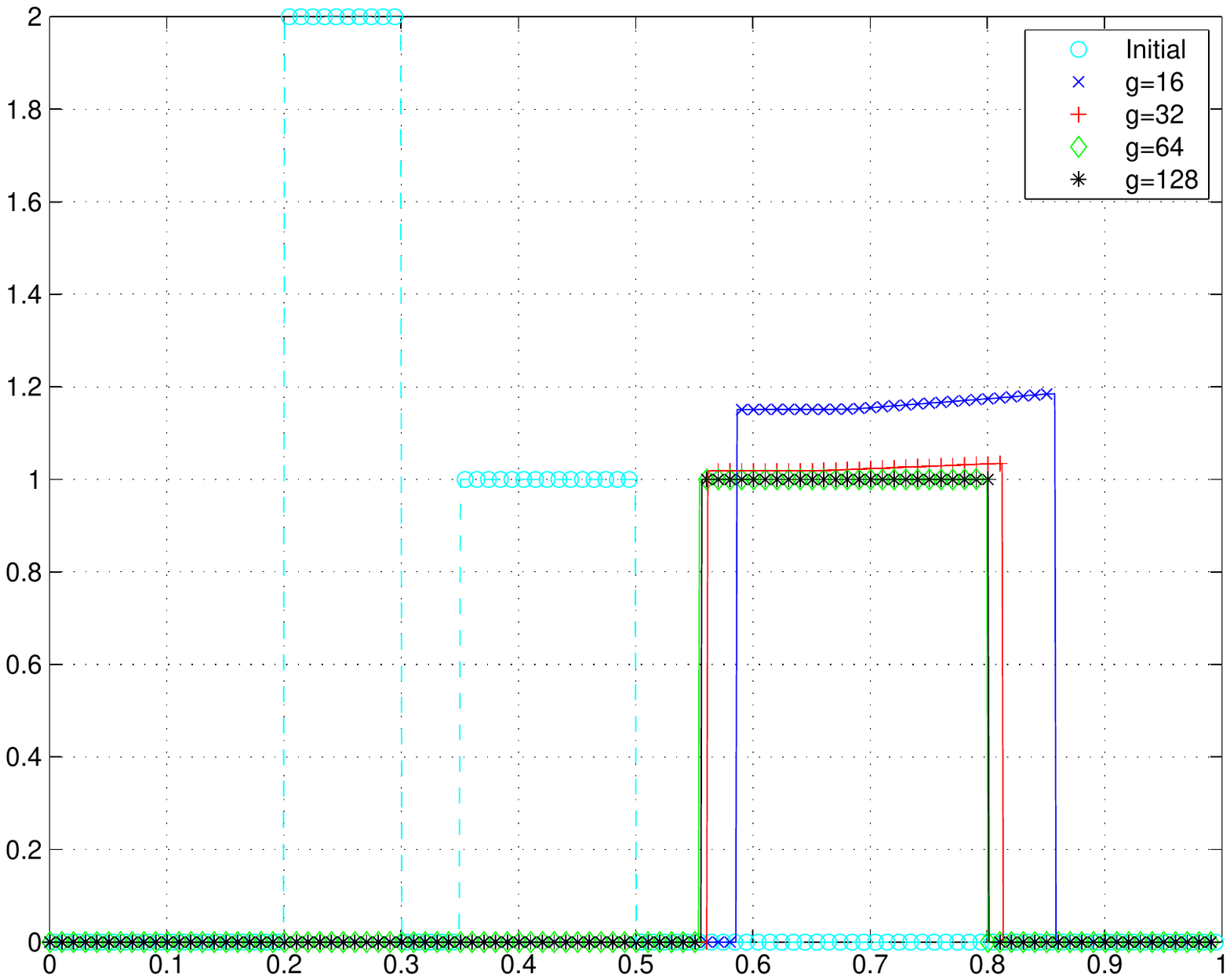}
	}
\\ 
\subfloat[Density - Implicit-explicit scheme\label{2BlocPP:DG}]
	{
	\includegraphics[trim=0 200 0 200,clip=true,
 scale=0.35]
{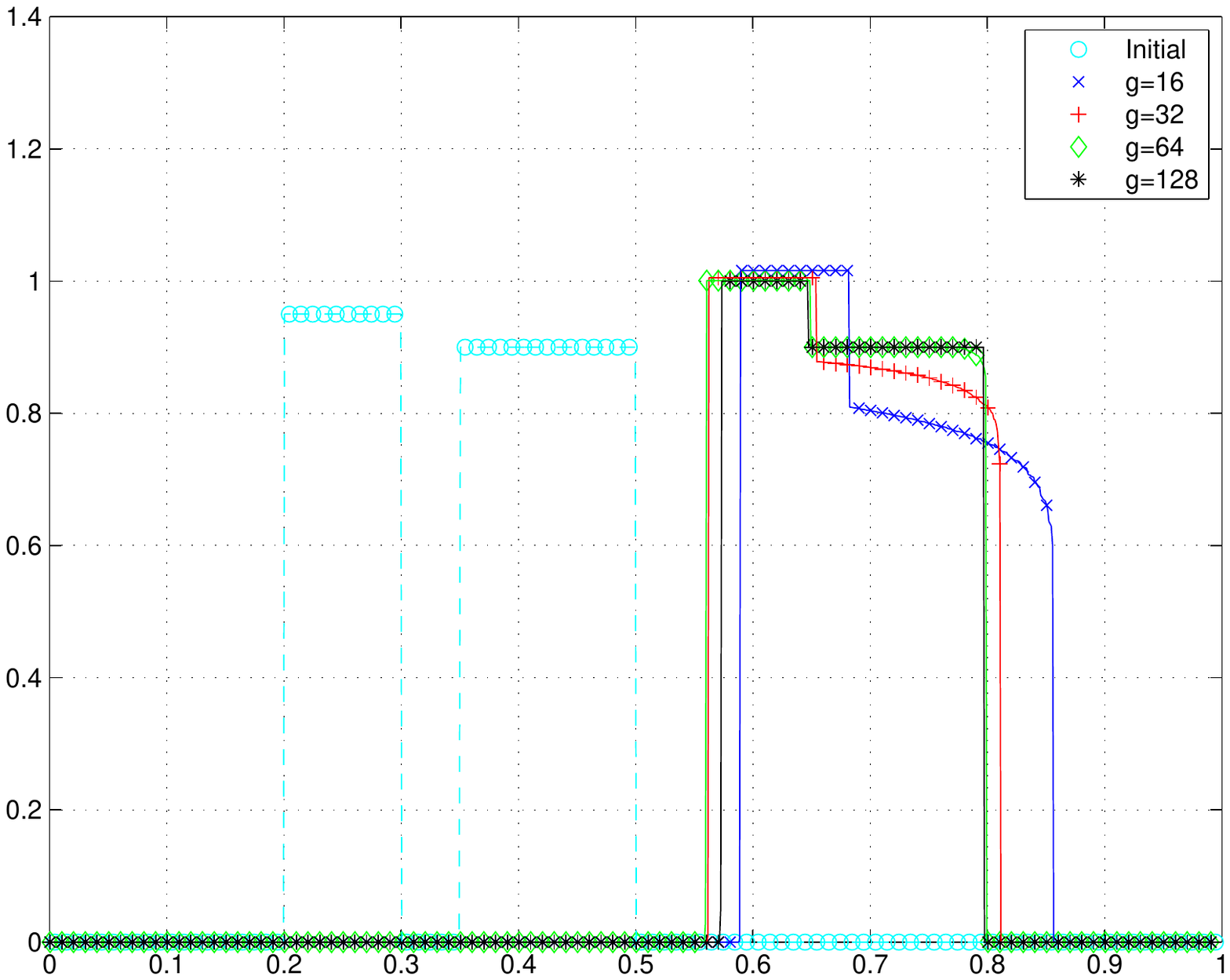}
	}
\quad
\subfloat[Velocity - Implicit-explicit scheme\label{2BlocPP:VG}]
	{
	\includegraphics[trim=0 200 0 200,clip=true,
 scale=0.35]
{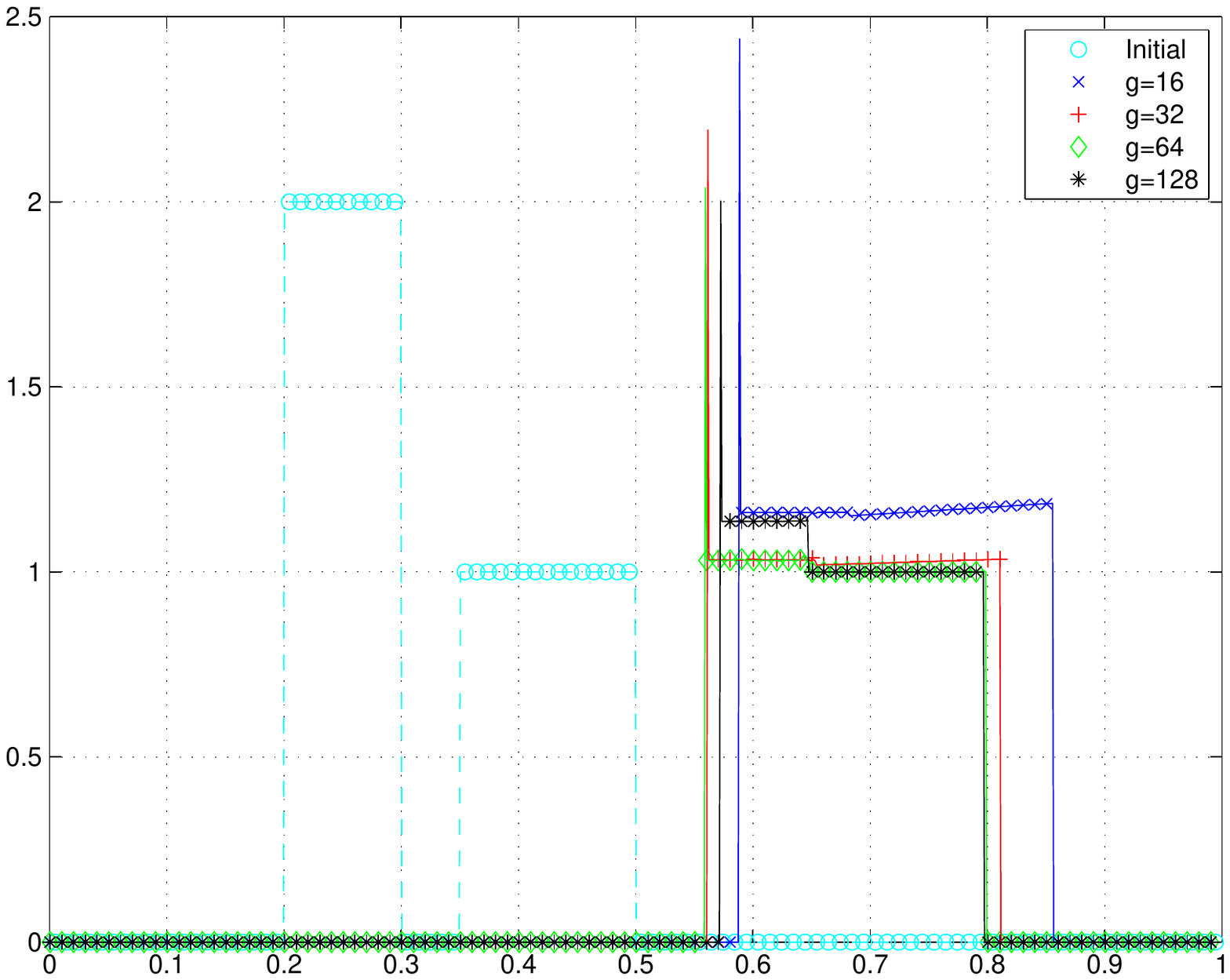}
	}
\caption{
\label{2BlocPP}
\textbf{Numerical results in the case of a shock between two blocks - Pressure \eqref{pr2} for different values of $\gamma$}:
$\gamma=16$ (blue), $\gamma=32$ (red), $\gamma=64$ (green) and $\gamma=128$ (black). 
On top, simulations are performed with the Glimm scheme and on the bottom, with the implicit-explicit scheme. We display the densities on the left and the velocities on the right.
The initial condition is plotted in cyan.
}
\end{figure}

We now consider the situation of a slow car cluster reached by a faster cluster of vehicles. The initial conditions are
\begin{align*}
	\rho^0(x)=
	\begin{cases}
      0.95 & \text{if } x\in\crochets{0.2,0.3}, \\
      0.9  & \text{if } x\in\crochets{0.35,0.5}, \\
      0    & \text{otherwise},
	\end{cases}
	&&
	v^0(x)=
	\begin{cases}
      2 & \text{if } x\in\crochets{0.2,0.3}, \\
      1 & \text{if } x\in\crochets{0.35,0.5}, \\
      0 & \text{otherwise}.
	\end{cases}
\end{align*}
and the expected solution at $T=0.3$ is
\begin{align}
\label{2blockSol}
	\rho(0.3,x)=
	\begin{cases}
      1    & \text{if } x\in\crochets{0.555,0.65}, \\
      0.9  & \text{if } x\in\crochets{0.65,0.8},   \\
      0    & \text{otherwise},
	\end{cases}
	&&
	v(0.3,x)=
	\begin{cases}
      1 & \text{if } x\in\crochets{0.555,0.8}, \\
      0 & \text{otherwise}.
	\end{cases}
\end{align}
Again, this solution presents a vacuum region.
The results for the different velocity offsets are represented in Figures~\ref{2BlocPO}, \ref{2BlocNP} and \ref{2BlocPP} for the scheme  presented in Section~\ref{scheme} and the Glimm scheme.

When using the Glimm scheme,
we obtain the expected limit solution as we make the parameters vary:
for  \eqref{pr1} as $\ve$ goes to 0, see 
Figure~\ref{2BlocPO}, for  \eqref{pr1bis} as $\ve$ goes to 0, see Figure~\ref{2BlocNP}, 
and for  \eqref{pr2} as $\gamma$ goes to infinity, see Figures~\ref{2BlocPP:DG} and \ref{2BlocPP:VG}.
The explicit-implicit  scheme  is able to compute precisely the car density
with the velocity offset \eqref{pr2}, see Figure~\ref{2BlocPP:DT}.
However, difficulties arise for the computation of the velocity, especially in the back of the 
jam, see  
Figure~\ref{2BlocPP:VT}: the velocity of the last cars in the congestion (coordinate $x=0.555$) is  over estimated.
Moreover, we observe in the case when $\gamma=128$ a loss of the quantity of cars between initial and final times of around $9 \%$.
This behavior is even worse with 
\eqref{pr1bis}, see Figures~\ref{2BlocNP:DG} and  \ref{2BlocNP:VG}. 
Here again, the quantity of mass is not conserved: the difference between the initial density and the final one is around $23\%$.
We point out that this simulation is quite tough; 
the results obtained with the Glimm scheme are neater
but  at the price of a significative numerical cost.

\subsubsection{Riemann problems: comparison with \cite{DeDe}}

In this section we present the numerical results obtained with the Glimm scheme and the implicit-explicit schemes of Section~\ref{scheme} for the sub-cases AI and AIII of \cite{BeDe,DeDe}. So we consider in the following a road  $\crochets{0,1}$ with the initial data
\begin{align}
\tag{{\bf{AI}}}\label{AI}
	\rho^0(x)=
	\begin{cases}
      0.7  & \text{if } x<0.5, \\
      0.5  & \text{if } x\geqslant0.5, \\
	\end{cases}
	&&
	v^0(x)=
	\begin{cases}
      0.5  & \text{if } x<0.5, \\
      0.1  & \text{if } x\geqslant0.5, \\
	\end{cases}
\end{align}
 and 
\begin{align}
\tag{{\bf{AIII}}}\label{AIII}
	\rho^0(x)=
	\begin{cases}
      0.7  & \text{if } x<0.5, \\
      0.5  & \text{if } x\geqslant0.5, \\
	\end{cases}
	&&
	v^0(x)=
	\begin{cases}
      0.1  & \text{if } x<0.5, \\
      0.5  & \text{if } x\geqslant0.5, \\
	\end{cases}
\end{align}
where we use the same labels as in \cite{BeDe,DeDe}. A solution for  \eqref{AI}  is 
\begin{align}
\label{caseAIsol}
	\rho(t,x)=
	\begin{cases}
      0.7  & \text{if } x<0.5-\frac{25}{30}t, \\
      1    & \text{if } x\in\crochets{0.5-\frac{25}{30}t,0.5+0.1t},	\\
      0.5  & \text{if } x\geqslant0.5+0.1t, \\
	\end{cases}
	&&
	v(t,x)=
	\begin{cases}
      0.5  & \text{if } x<0.5-\frac{25}{30}t, \\
      0.1  & \text{if } x\geqslant0.5-\frac{25}{30}t, \\
	\end{cases}
\end{align}
and
\begin{align}
\label{caseAIIIsol}
	\rho(t,x)=
	\begin{cases}
      0.7  & \text{if } x<0.5+0.1t, \\
      0    & \text{if } x\in\crochets{0.5+0.1t,0.5+0.5t},	\\
      0.5  & \text{if } x\geqslant0.5+0.5t, \\
	\end{cases}
	&&
	v(t,x)=
	\begin{cases}
      0.1  & \text{if } x<0.5+0.1t, \\
      0.5  & \text{if } x\geqslant0.5+0.5t, \\
	\end{cases}
\end{align}
is a  solution for \eqref{AIII} (note that $v$ does not make  sense in the vacuum region $x\in\crochets{0.5+0.1t,0.5+0.5t}$).
In order to compare our scheme with the results of \cite{DeDe} we take the same parameters, namely for the velocity offset \eqref{pr1} and \eqref{pr1bis} we use $\gamma=1$ and $\ve=10^{-3}$. 
For the velocity offset \eqref{pr2} we use $\gamma=64$.

\begin{figure}
\subfloat[Density at t=0.2\label{DeDe:CoD_t1}]
	{
\includegraphics[trim=0 200 0 200,clip=true,
 scale=0.3]
{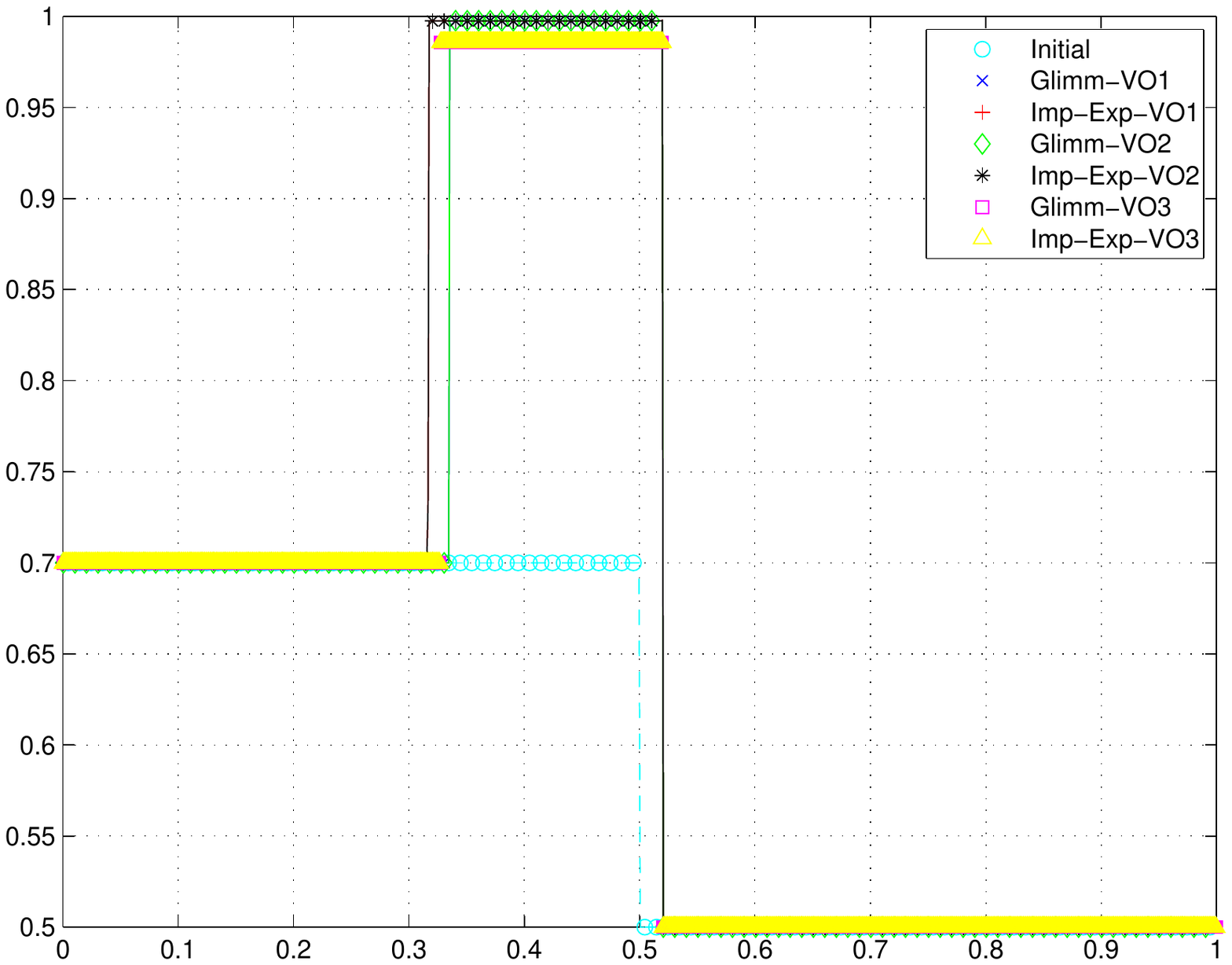}
	}
\quad
\subfloat[Velocity at t=0.2\label{DeDe:CoV_t1}]
	{
	\includegraphics[trim=0 200 0 200,clip=true,
 scale=0.3]
{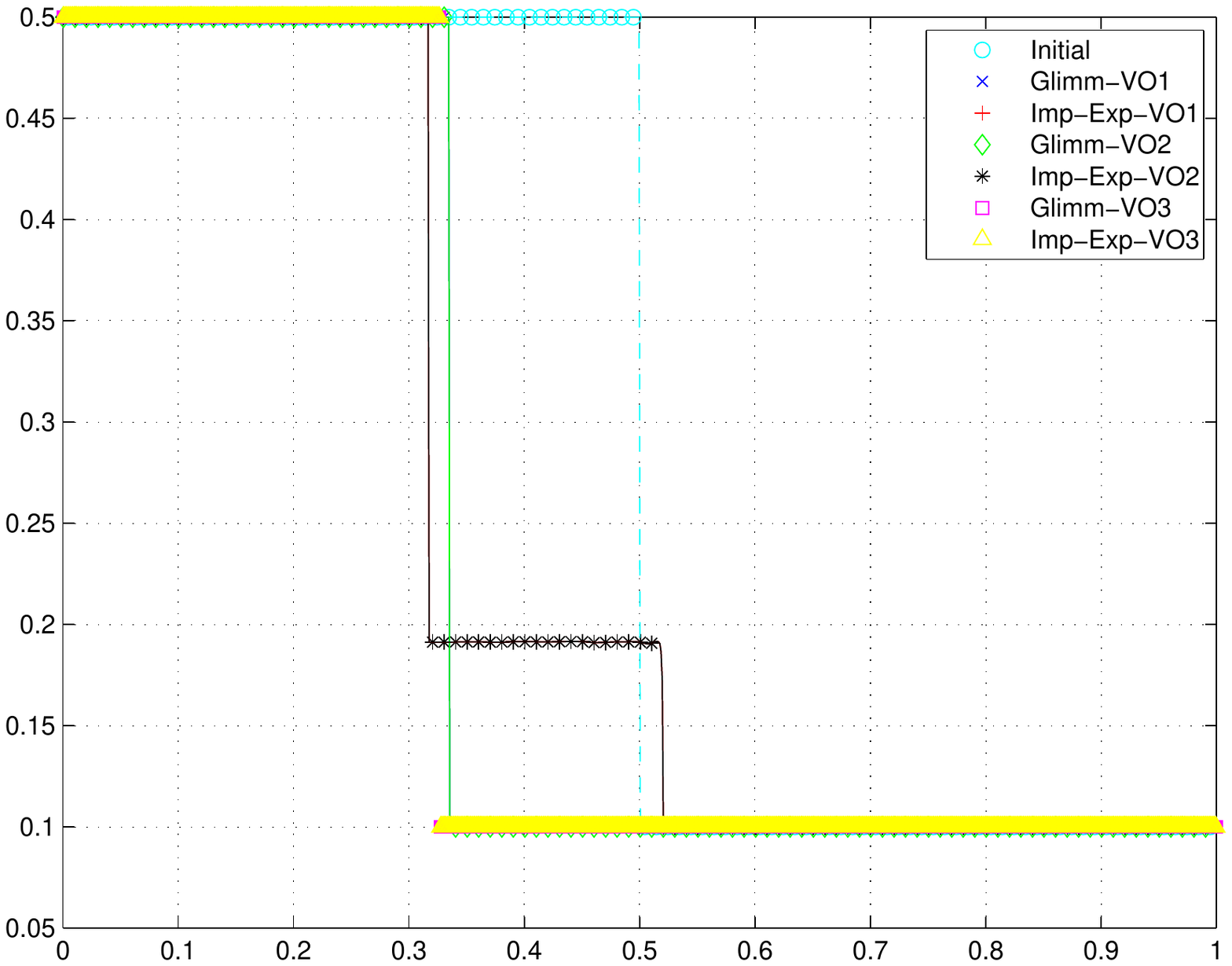}
	}
\\ 
\subfloat[Density at t=0.4\label{DeDe:CoD_t2}]
	{
	\includegraphics[trim=0 200 0 200,clip=true,
 scale=0.3]
{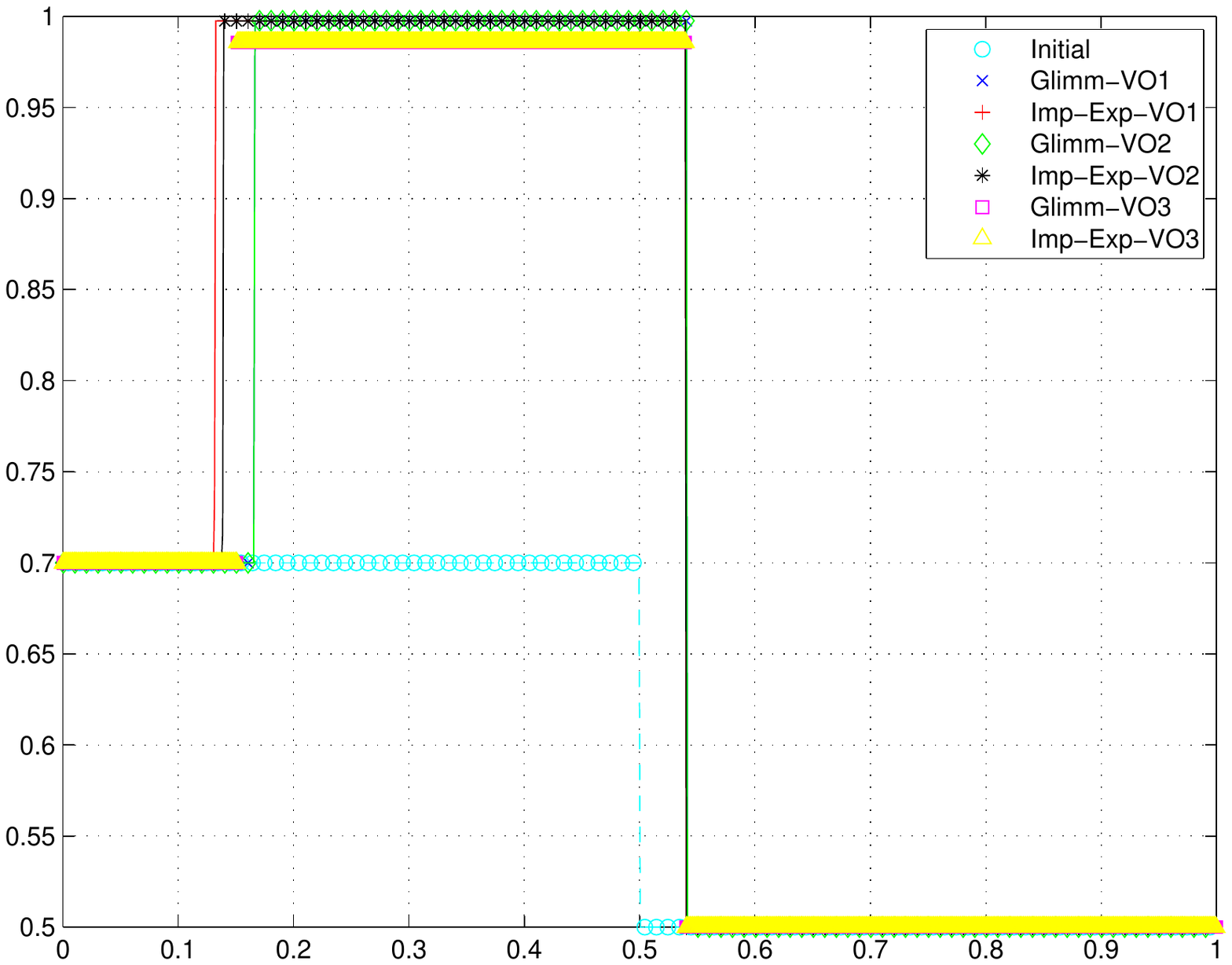}
	}
\quad
\subfloat[Velocity at t=0.4\label{DeDe:CoV_t2}]
	{
	\includegraphics[trim=0 200 0 200,clip=true,
 scale=0.3]
{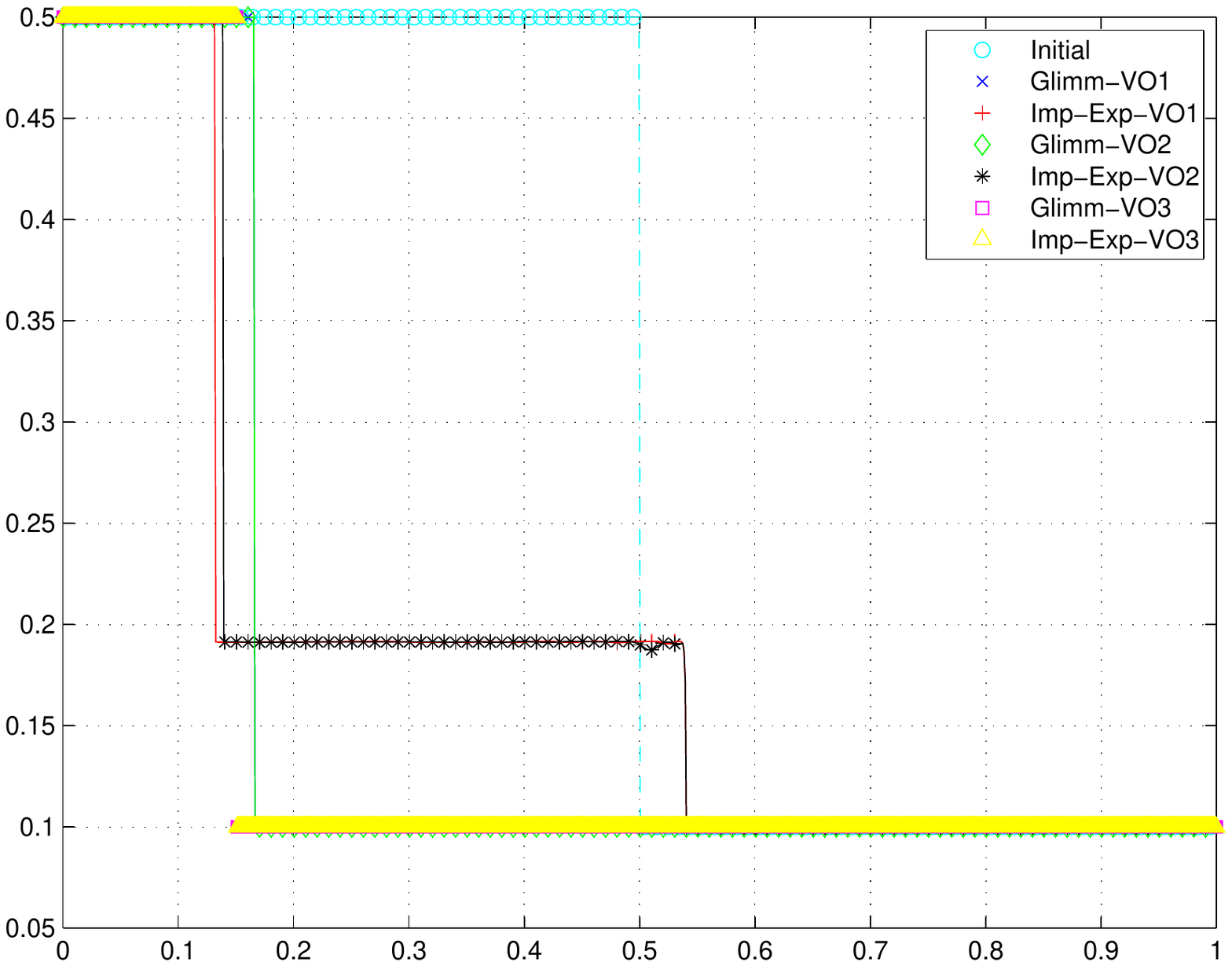}
	}
\\ 
\subfloat[Density at t=0.6\label{DeDe:CoD_t3}]
	{
	\includegraphics[trim=0 200 0 200,clip=true,
 scale=0.3]
{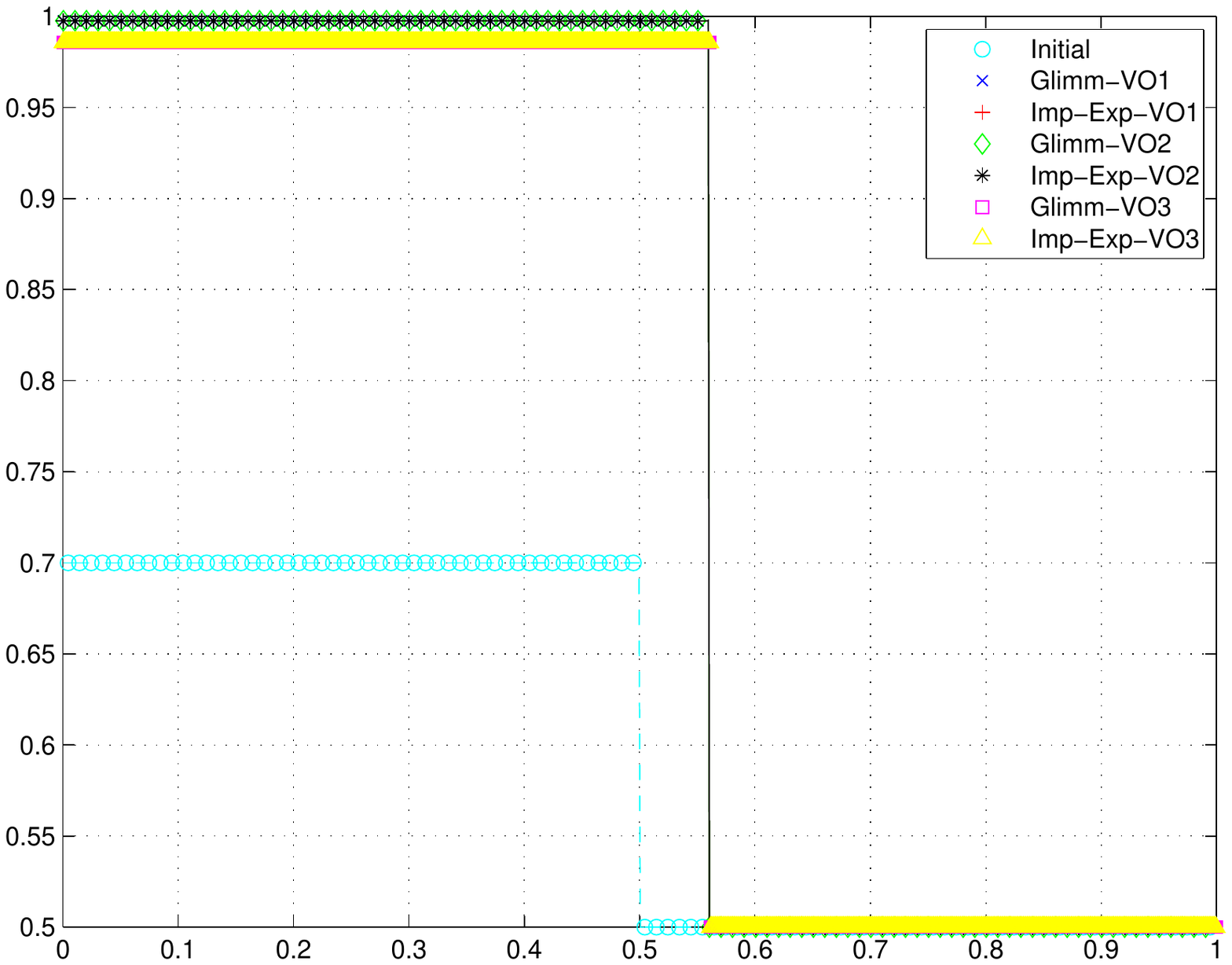}
	}
\quad
\subfloat[Velocity at t=0.6\label{DeDe:CoV_t3}]
	{
	\includegraphics[trim=0 200 0 200,clip=true,
 scale=0.3]
 {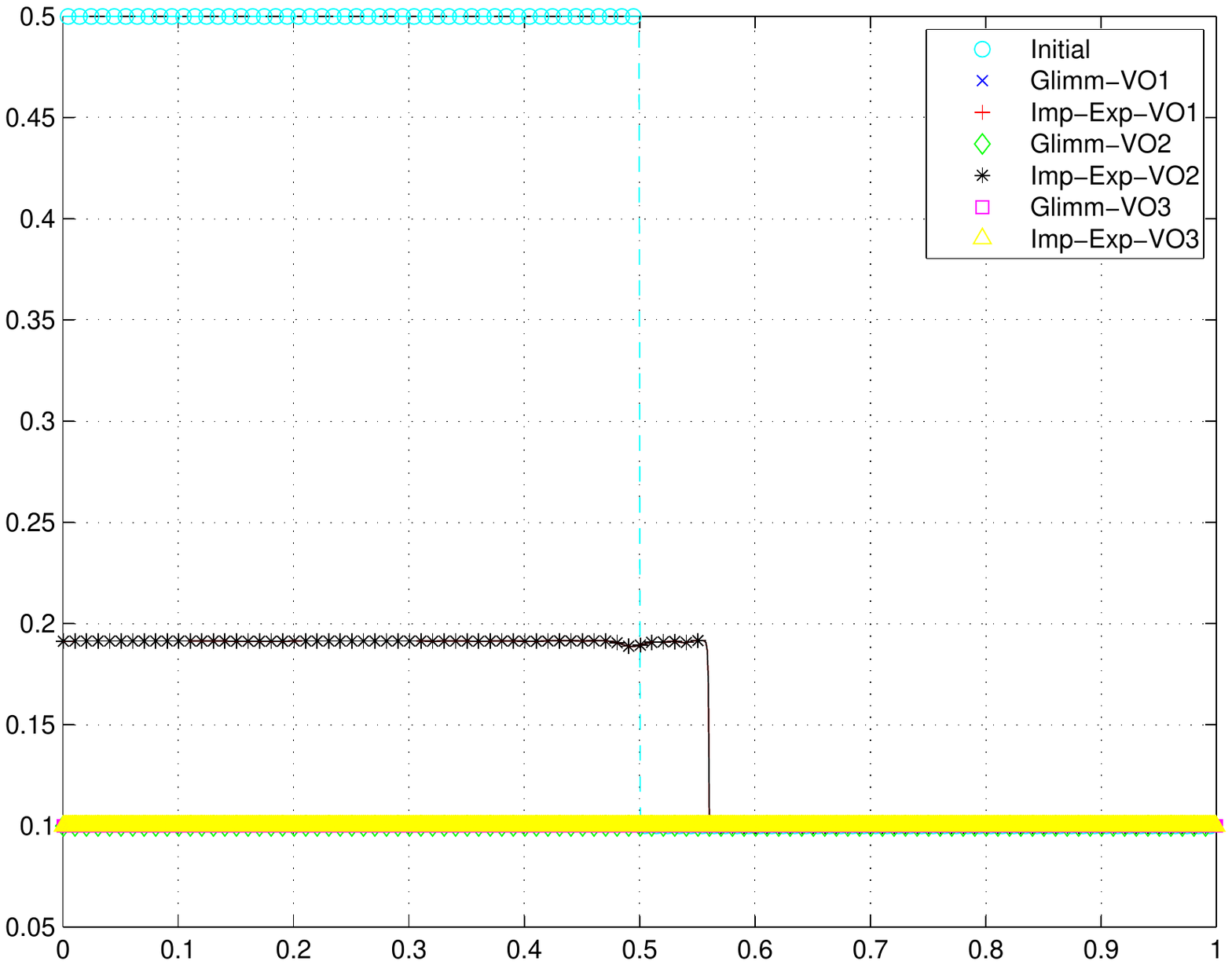}
	}
\caption{
\label{DeDe:Cofig}
\textbf{Comparaison of the numerical scheme and velocity offset for the initial data \eqref{AI}.}
We display the densities on the left and the velocities on the right, at $t=0.2$ (top), $t=0.4$ (middle) and $t=0.6$ (bottom).
Glimm scheme with pressures \eqref{pr1} (in  blue) , \eqref{pr1bis}  (in  green) and \eqref{pr2} (in pink); implicit-explicit  scheme with pressures \eqref{pr1} (in red) , \eqref{pr1bis}  (in black) and \eqref{pr2} (in yellow). 
}
\end{figure}

Figure~\ref{DeDe:Cofig} represents the results for the sub-case~\eqref{AI} computed with the three velocity offsets and with the Glimm  scheme or the explicit-implicit scheme, at time $t=0.2$, $t=0.4$ and $t=0.6$. According to Figures~\ref{DeDe:CoD_t1}, \ref{DeDe:CoD_t2} and \ref{DeDe:CoD_t3}, the explicit-implicit scheme for the velocity offset \eqref{pr1bis} overestimates the length of the congestion. For example at $t=0.4$ using \eqref{caseAIsol},  the tail of the congestion should be at $x=0.166$ instead of $x=0.13$. Moreover, the velocity in the congested area is over estimated, see Figures~\ref{DeDe:CoV_t1}, \ref{DeDe:CoV_t2} and \ref{DeDe:CoV_t3}. 

\begin{figure}
\subfloat[Density at t=0.27\label{DeDe:DeD_t1}]
	{
	\includegraphics[trim=0 200 0 200,clip=true,
 scale=0.3]
{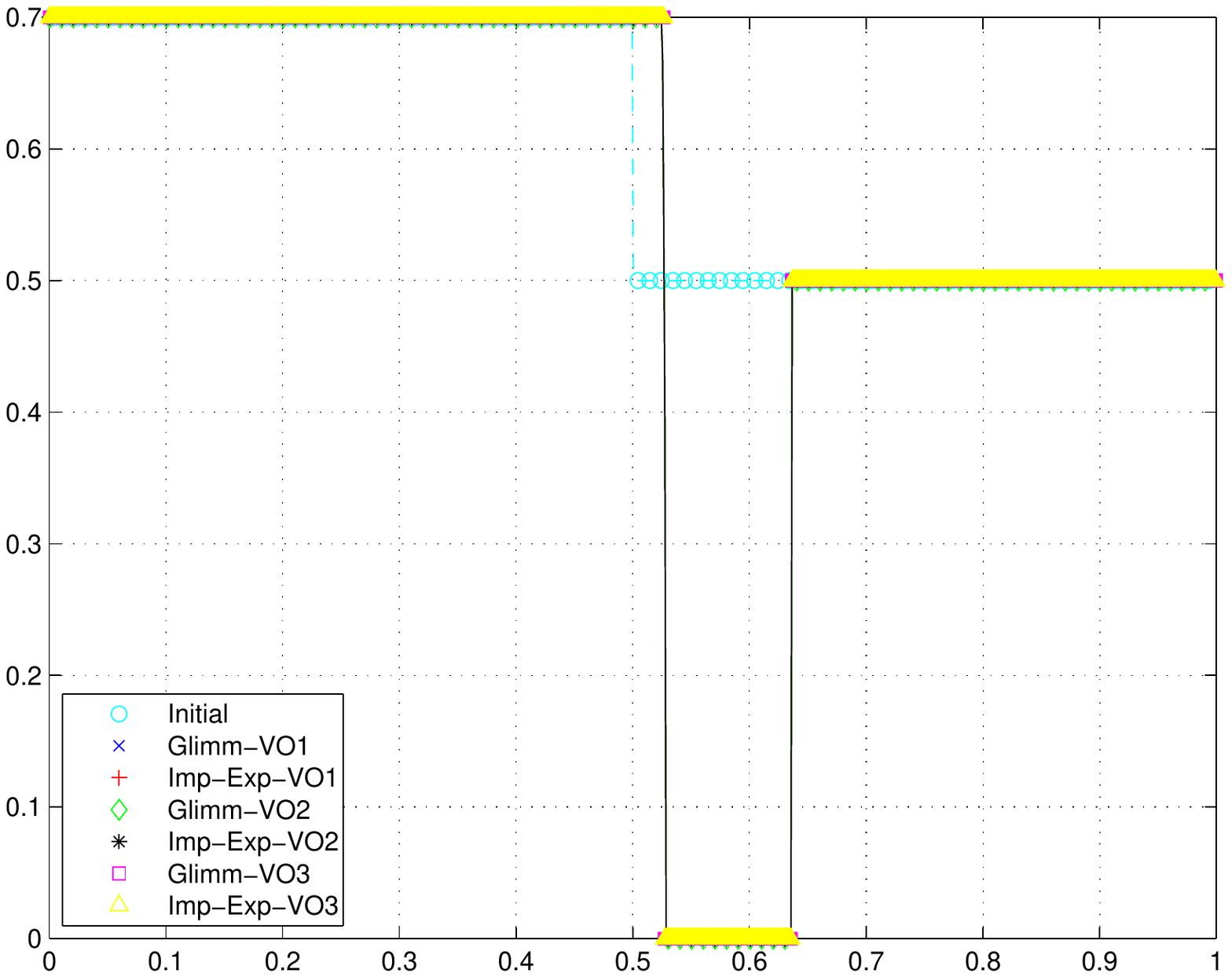}

	}
\quad
\subfloat[Velocity at t=0.27\label{DeDe:DeV_t1}]
	{
	\includegraphics[trim=0 200 0 200,clip=true,
 scale=0.3]
{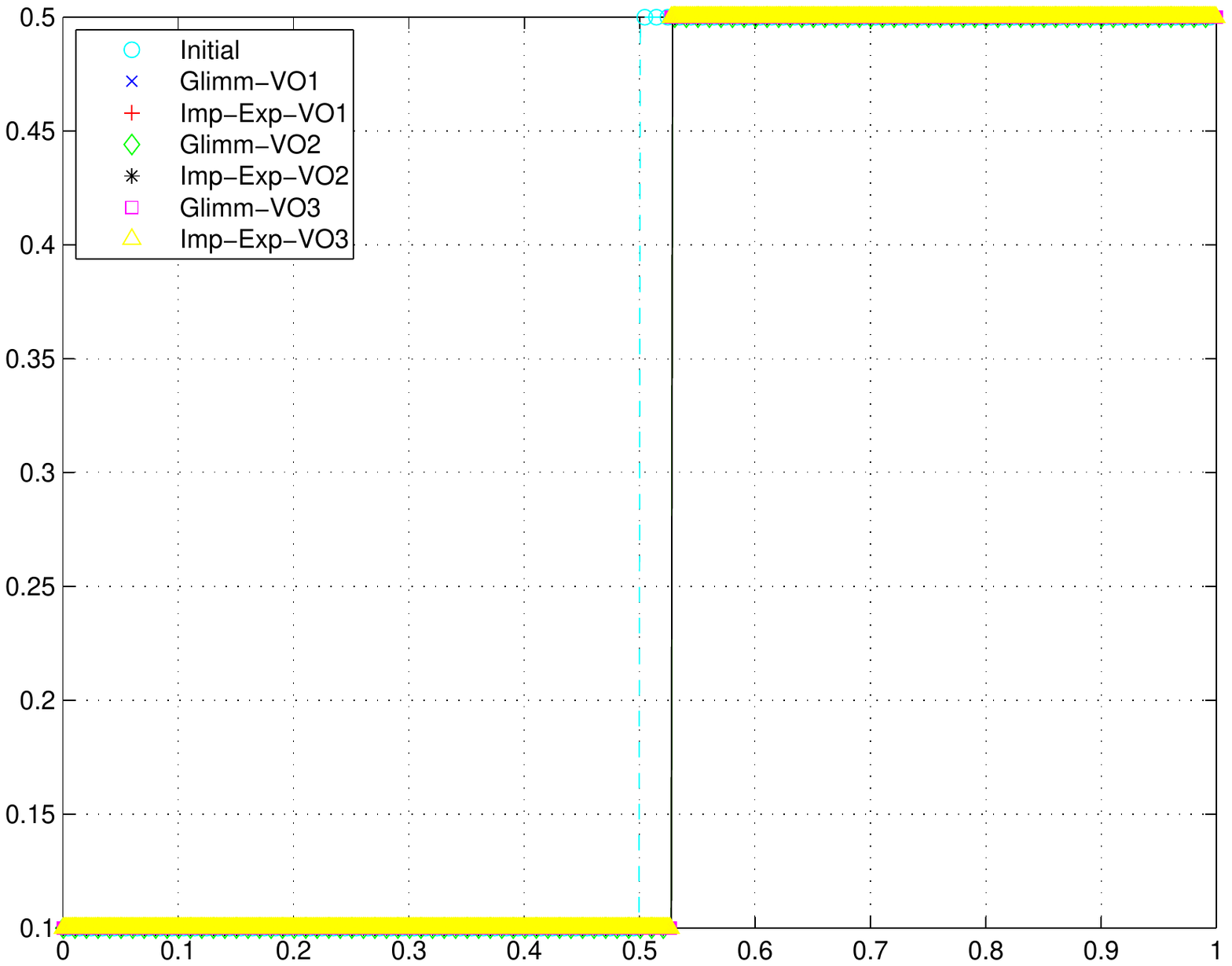}
	}
\\ 
\subfloat[Density at t=0.53\label{DeDe:DeD_t2}]
	{
	\includegraphics[trim=0 200 0 200,clip=true,
 scale=0.3]
{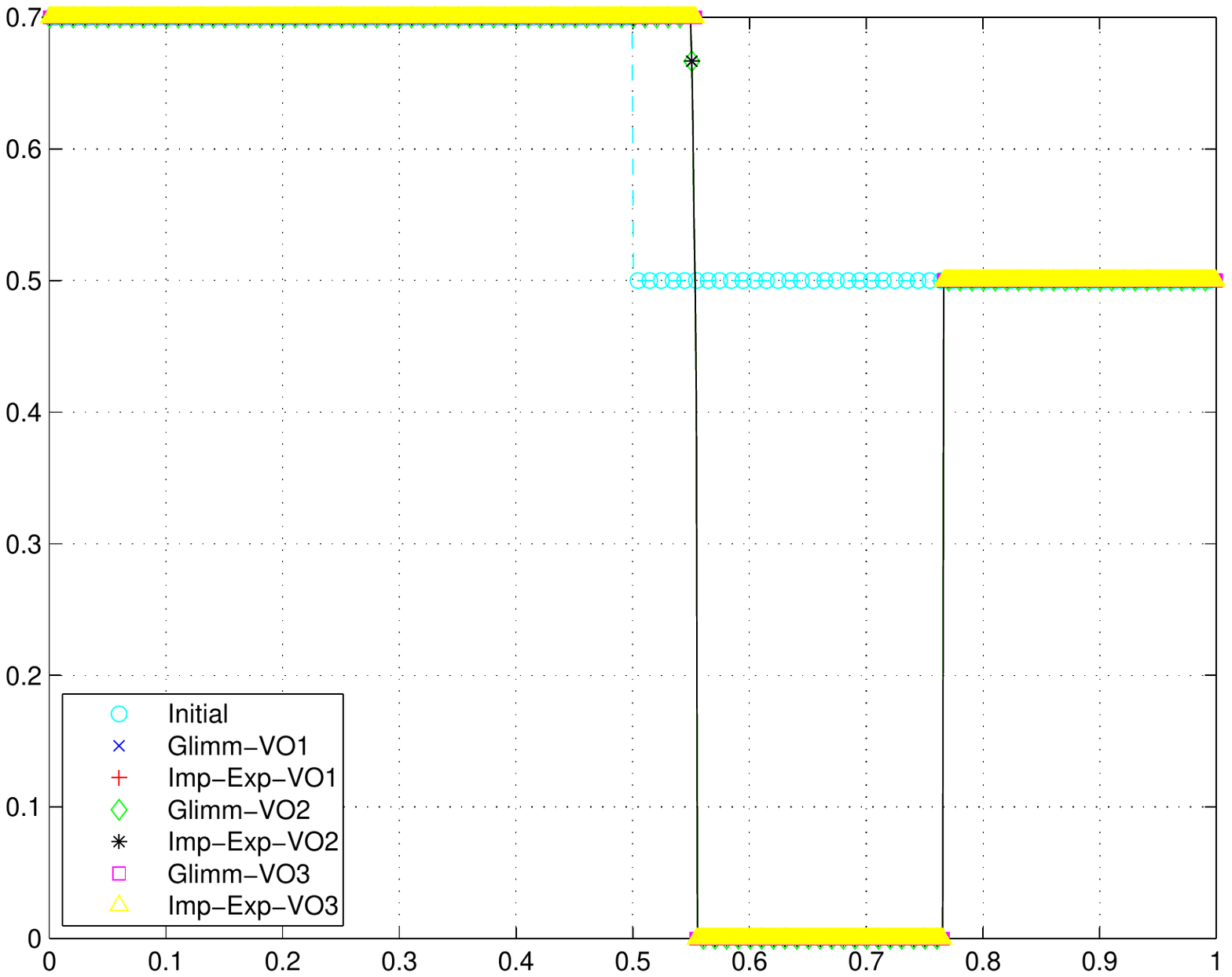}
	}
\quad
\subfloat[Velocity at t=0.53\label{DeDe:DeV_t2}]
	{
	\includegraphics[trim=0 200 0 200,clip=true,
 scale=0.3]
{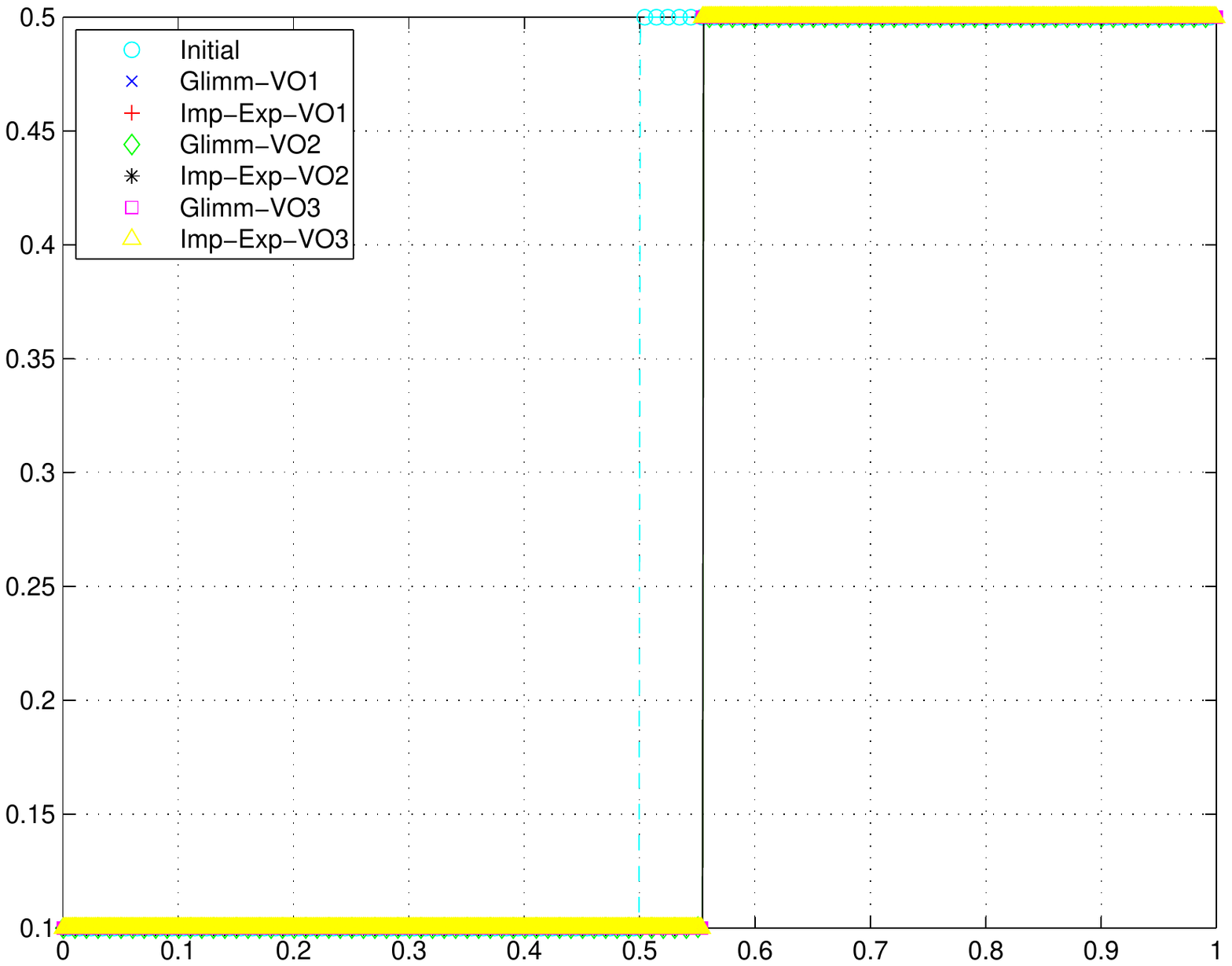}
	}
\\ 
\subfloat[Density at t=0.8\label{DeDe:DeD_t3}]
	{
	\includegraphics[trim=0 200 0 200,clip=true,
 scale=0.3]
{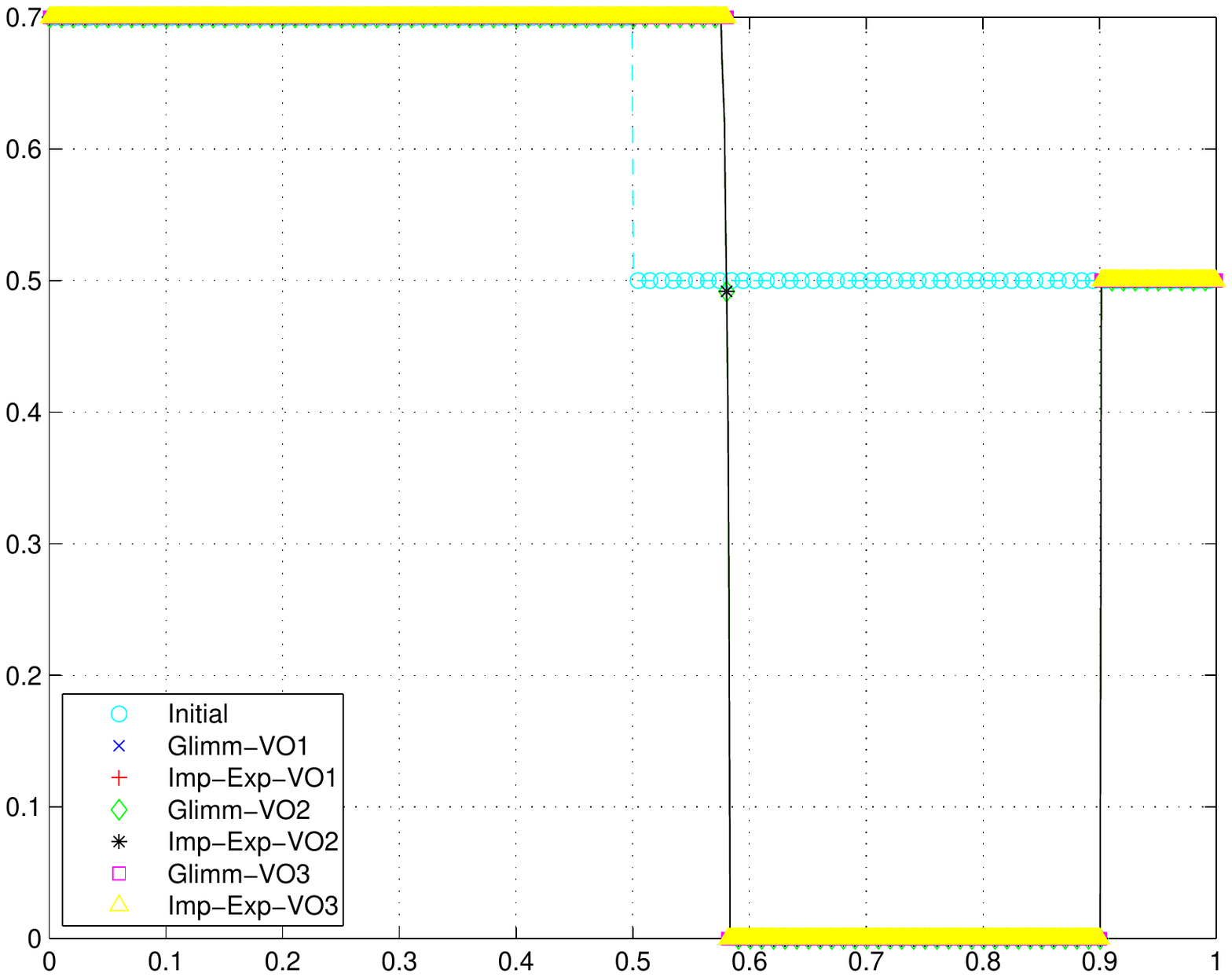}
	}
\quad
\subfloat[Velocity at t=0.8\label{DeDe:DeV_t3}]
	{
	\includegraphics[trim=0 200 0 200,clip=true,
 scale=0.3]
{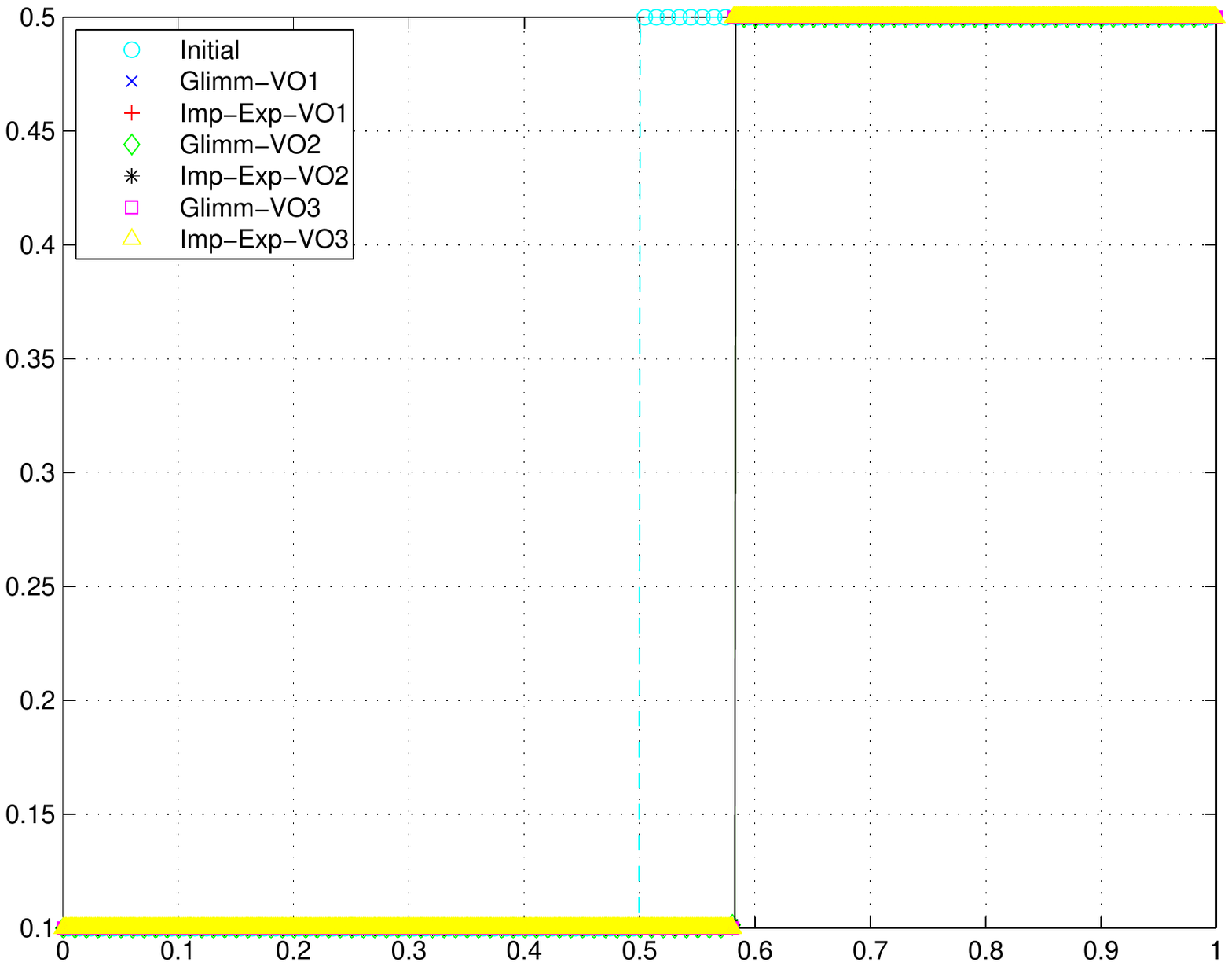}
	}
\caption{
\label{DeDe:Defig}
\textbf{Comparaison of the numerical scheme and velocity offset for the initial data \eqref{AIII}.}
We display the densities on the left and the velocities on the right, at $t=0.27$ (top), $t=0.53$ (middle) and $t=0.8$ (bottom).
Glimm scheme with pressures \eqref{pr1} (in  blue) , \eqref{pr1bis}  (in  green) and \eqref{pr2} (in pink); implicit-explicit  scheme with pressures \eqref{pr1} (in red) , \eqref{pr1bis}  (in black) and \eqref{pr2} (in yellow). 
}
\end{figure}

Figure~\ref{DeDe:Defig} represents the results for the sub-case~\eqref{AIII} for the three velocity offsets and for different times. In order to ease the comparison between the Glimm scheme and explicit-implicit scheme presented in Section~\ref{scheme}, the results for the two numerical approaches are displayed simultaneously.  According to Figure~\ref{DeDe:Defig}, all the different numerical approaches give similar results which  coincide with the solution given by~\eqref{caseAIIIsol}.

\section{Conclusion}

The model \eqref{PGCD} is intended to describe the formation and the dynamics of traffic jams, through a Lagrange multiplier that accounts 
for a density threshold.
This model can be motivated, at least formally, through asymptotic 
arguments 
from the Aw-Rascle-Zhang system with a rescaled 
velocity--offset.
It raises the question of simulating efficiently the Aw-Rascle-Zhang system 
with potentially stiff velocity offsets.
Depending 
on the values of the parameters
it can be seen either as the simulation of a 
model for traffic flows with stiff parameters 
or as a way to access the limiting  
behavior described by 
\eqref{PGCD}, alternative for instance to the approach of \cite{Maury}.
However the scaling 
induces fast propagation waves and, in turn, severe stability conditions.
In this paper, we propose several approaches to obtain asymptotically \eqref{PGCD} and we introduce an implicit--explicit
method in order to cope with the large characteristic speeds of the system.
\\

This study exhibits numerical difficulties, related to both the lack of convexity of the invariant domains of \eqref{AW1}
and the large characteristic speeds. We have proposed a time--splitting method, 
based on a decomposition of the velocity--offset and the use of the Glimm scheme which avoids the non admissible 
solutions produced by schemes based on a projection step.
Our findings bring out that  the behavior of the system \eqref{PGCD} 
can be obtained asymptotically, but the shape of the solution for intermediate 
values of the scaling parameters highly depends on the expression of the penalized velocity offset.
It means that a serious modeling work should decide what is the most appropriate model.

\section*{Acknowledgements}
We thank Fr\'ed\'eric Coquel for friendly advices and warm encouragements during the preparation of this work.

\bibliography{Biblio_Art_Traffic}

\begin{thebibliography}{10}

\bibitem{AwR2}
A.~Aw, A.~Klar, T.~Materne, and M.~Rascle.
\newblock Derivation of continuum traffic flow models from microscopic
  follow-the-leader models.
\newblock {\em SIAM J. Appl. Math.}, 63(1):259--278, 2002.

\bibitem{AwR1}
A.~Aw and M.~Rascle.
\newblock Resurrection of "second order'' models of traffic flow.
\newblock {\em SIAM J. Appl. Math.}, 60(3):916--938 (electronic), 2000.

\bibitem{BBNS}
N.~Bellomo, A.~Bellouquid, J.~Nieto, and J.~Soler.
\newblock On the multiscale modeling of vehicular traffic: From kinetic to
  hydrodynamics.
\newblock {\em Disc. Cont. Dyn. Syst.--B}, 19(7):1869--1888, 2014.

\bibitem{BeDo}
N.~Bellomo and C.~Dogb\'e.
\newblock On the modeling of traffic and crowds:a survey of models,
  speculations, and perspectives.
\newblock {\em SIAM Rev.}, 53(3):409--463, 2011.

\bibitem{BGS}
S.~Benzoni-Gavage and D.~Serre.
\newblock {\em Multidimensional Hyperbolic Partial Differential Equations:
  First-Order Systems and Applications}.
\newblock Oxford Mathematical Monographs. Oxford University Press, 2006.

\bibitem{BeBr}
F.~Berthelin and D.~Broizat.
\newblock A model for the evolution of traffic jams in multilane.
\newblock {\em Kinetic and Related Models}, 5(4):697--728, 2012.

\bibitem{BeDe}
F.~Berthelin, P.~Degond, M.~Delitala, and M.~Rascle.
\newblock A model for the formation and evolution of traffic jams.
\newblock {\em Arch. Rational Mech. Anal.}, 187:185--220, 2008.

\bibitem{BGM}
F.~Berthelin, T.~Goudon, and S.~Minjeaud.
\newblock Multifluid flows: a kinetic approach.
\newblock {\em J. Sci. Comput.}, 66(2):792--824, 2016.

\bibitem{Bouc}
F.~Bouchut.
\newblock {\em Nonlinear stability of finite volume methods for hyperbolic
  conservation laws and well--balanced schemes for sources}.
\newblock Frontiers in math. Birkh\"auser, 2004.

\bibitem{BJL}
F.~Bouchut, S.~Jin, and X.~Li.
\newblock Numerical approximations of pressureless and isothermal gas dynamics.
\newblock {\em SIAM J. Numer. Anal.}, 41(1):135--158 (electronic), 2003.

\bibitem{BrGr}
Y.~Brenier and E.~Grenier.
\newblock Sticky particles and scalar conservation laws.
\newblock {\em SIAM J. Numer. Anal.}, 35(6):2317--2328, 1998.

\bibitem{ChG}
C.~Chalons and P.~Goatin.
\newblock Transport-equilibrium schemes for computing contact discontinuities
  in traffic flow modeling.
\newblock {\em Commun. Math. Sci.}, 5(3):533--551, 2007.

\bibitem{CG08}
C.~Chalons and P.~Goatin.
\newblock Godunov scheme and sampling technique for computing phase transitions
  in traffic flow modeling.
\newblock {\em Interfaces Free Bound.}, 10(2):197--221, 2008.

\bibitem{colella}
P.~Colella.
\newblock Glimm's method for gas dynamics.
\newblock {\em SIAM J. Sci. Statist. Comput.}, 3(1):76--110, 1982.

\bibitem{daganzo}
C.~F. Daganzo.
\newblock Requiem for second-order fluid approximations of traffic flow.
\newblock {\em Transportation Research Part B: Methodological}, 29(4):277--286,
  1995.

\bibitem{DeDe}
P.~Degond and M.~Delitala.
\newblock Modelling and simulation of vehicular traffic jam formation.
\newblock {\em Kinet. Relat. Models}, 1(2):279--293, 2008.

\bibitem{DH}
P.~Degond and J.~Hua.
\newblock Self-organized hydrodynamics with congestion and path formation in
  crowds.
\newblock {\em J. Comput. Phys.}, 237:299--319, 2013.

\bibitem{DHN}
P.~Degond, J.~Hua, and L.~Navoret.
\newblock Numerical simulations of the {E}uler system with congestion
  constraint.
\newblock {\em J. Comput. Phys.}, 230(22):8057--8088, 2011.

\bibitem{DeTa}
P.~Degond and M.~Tang.
\newblock All speed scheme for the low {M}ach number limit of the isentropic
  {E}uler equations.
\newblock {\em Commun. Comput. Phys.}, 10(1):1--31, 2011.

\bibitem{gazis}
D.~C. Gazis, R.~Herman, and R.~W. Rothery.
\newblock Nonlinear follow-the-leader models of traffic flow.
\newblock {\em Operations research}, 9(4):545--567, 1961.

\bibitem{Gli}
J.~Glimm.
\newblock Solutions in the large for nonlinear hyperbolic systems of equations.
\newblock {\em Comm. Pure Applied Math.}, 18:697--715, 1965.

\bibitem{Gr}
E.~Grenier.
\newblock Existence globale pour le syst\`eme des gaz sans pression.
\newblock {\em Comptes Rendus Acad. Sci.}, pages 171--174, 1995.

\bibitem{Ju}
J.~Jung.
\newblock {\em Sch\'emas num\'eriques adapt\'es aux acc\'el\'erateurs
  multicoeurs pour les \'ecoulements bifluides}.
\newblock PhD thesis, Univ. Strasbourg, 2014.

\bibitem{AYLR}
A.-Y. Le~Roux.
\newblock Stability for some equations of gas dynamics.
\newblock {\em Math. Comput.}, 37(156):307--320, 1981.

\bibitem{LW}
M.~J. Lighthill and G.~B. Whitham.
\newblock On kinematic waves. {II}. {A} theory of traffic flow on long crowded
  roads.
\newblock {\em Proc. Roy. Soc. London. Ser. A.}, 229:317--345, 1955.

\bibitem{LiMa}
P.-L. Lions and N.~Masmoudi.
\newblock On a free boundary barotropic model.
\newblock {\em Ann. Inst. H. Poincar\'e Anal. Non Lin\'eaire}, 16(3):373--410,
  1999.

\bibitem{liu}
T.~P. Liu.
\newblock The deterministic version of the {G}limm scheme.
\newblock {\em Comm. Math. Phys.}, 57(2):135--148, 1977.

\bibitem{Maury}
B.~Maury and A.~Preux.
\newblock Pressureless {E}uler equations with maximal density constraint : a
  time-splitting scheme.
\newblock Technical report, Universit\'e Paris-Sud, 2015.
\newblock Available on \url{https://hal.archives-ouvertes.fr/hal-01224008}.

\bibitem{Nel}
P.~Nelson and A.~Sopasakis.
\newblock The {P}rigogine--{H}erman kinetic model predicts widely scattered
  traffic flow data at high concentrations.
\newblock {\em Transportation Research Part B: Methodological}, 32(8):589--604,
  1998.

\bibitem{PF}
S.~L. Paveri-Fontana.
\newblock On {B}oltzmann-like treatments for traffic flow: A critical review of
  the basic model and an alternative proposal for dilute traffic analysis.
\newblock {\em Transportation Research}, 9:225--235, 1975.

\bibitem{payne}
H.~J. Payne.
\newblock Freflo: A macroscopic simulation model of freeway traffic.
\newblock {\em Transportation Research Record}, (722), 1979.

\bibitem{Prig}
I.~Prigogine and R.~Herman.
\newblock {\em Kinetic Theory of Vehicular Traffic}.
\newblock American Elsevier Publishing, 1971.

\bibitem{PSTV}
G.~Puppo, M.~Semplice, A.~Tosin, and G.~Visconti.
\newblock Fundamental diagrams in traffic flow: the case of heterogeneous
  kinetic models.
\newblock {\em Communications in mathematical sciences}, 2015.

\bibitem{smoller}
J.~Smoller.
\newblock {\em Shock waves and reaction-diffusion equations}, volume 258 of
  {\em Grundlehren der Mathematischen Wissenschaften [Fundamental Principles of
  Mathematical Sciences]}.
\newblock Springer-Verlag, New York, second edition, 1994.

\bibitem{toro}
E.~F. Toro.
\newblock {\em Riemann solvers and numerical methods for fluid dynamics}.
\newblock Springer-Verlag, Berlin, third edition, 2009.
\newblock A practical introduction.

\bibitem{WK}
R.~Wegener and A.~Klar.
\newblock A kinetic model for vehicular traffic derived from a stochastic
  microscopic model.
\newblock {\em Transport Theory and Stat. Phys.}, 25:785--798, 1996.

\bibitem{Zhang}
H.~M. Zhang.
\newblock A non-equilibrium traffic model devoid of gas-like behavior.
\newblock {\em Transportation Research Part B: Methodological}, 36(3):275--290,
  2002.

\end{thebibliography}
\bibliographystyle{abbrv}

\end{document}